\documentclass[a4paper,12pt,twoside]{article}
\usepackage{graphicx}\usepackage{a4wide}
\usepackage{amssymb,amsmath,amsthm}
\usepackage{yhmath} 
\usepackage{mathdots}
\usepackage[only,llbracket,rrbracket,interleave]{stmaryrd}

\usepackage{array}
\usepackage{multirow}
\usepackage{tabularx}
\usepackage{extarrows}
\usepackage{booktabs}

\usepackage{subfig}
\usepackage[sort,compress]{cite}
\usepackage{enumitem}

\newcommand{\bm}[1]{\text{\boldmath $#1$\unboldmath}}
\newcommand{\abs}[1]{\lvert#1\rvert}
\newcommand{\norm}[1]{\lVert#1\rVert}

\newcommand{\Div}{{\bm{\nabla}{\cdot}}}
\newcommand{\grad}{\bm{\nabla}}

\newcommand{\RR}{\mathbb{R}}
\newcommand{\C}{\mathrm{C}}
\newcommand{\fluxE}{\mathrm{J}_S}

\newcommand{\sobo}[1][1]{\ensuremath{\mathcal{H}^{#1}}}
\newcommand{\eltwo}{\ensuremath{\mathcal{L}^2}}
\newcommand{\hDiv}[1]{\ensuremath{\mathcal{H}(\operatorname{div};{{#1}})}}

\newcommand{\nsd}  {\ensuremath{\texttt{n}_{\texttt{sd}}}}
\newcommand{\numel}{\ensuremath{\texttt{n}_{\texttt{el}}}}

\newcommand{\nface}{\ensuremath{\texttt{n}_{\texttt{fc}}}}

\newcommand{\numelfc}{\ensuremath{\texttt{n}_{\texttt{fc}}}}

\newcommand{\Insd}{\bm{I}_{\!\nsd\!}}

\newcommand{\bu}{\bm{u}}
\newcommand{\bhu}{\bm{\hat{u}}}
\newcommand{\bv}{\bm{v}}

\newcommand{\bw}{\bm{w}}
\newcommand{\bhw}{\bm{\hat{w}}}
\newcommand{\bn}{\bm{n}}

\newcommand{\bt}{\bm{t}}
\newcommand{\bx}{\bm{x}}
\newcommand{\bL}  {\bm{L}}
\newcommand{\bG}  {\bm{G}}
\newcommand{\buS}{\bu_{\star}}
\newcommand{\R}{\text{R}}
\newcommand{\Rpore}{\text{R}_{\text{p}}}
\newcommand{\Rint}{\text{R}_{\text{int}}}
\newcommand{\Rext}{\text{R}_{\text{ext}}}
\newcommand{\OmInt}{\omega_{\text{int}}}
\newcommand{\OmExt}{\omega_{\text{ext}}}

\newcommand{\buSe}{\bu_{\star,e}}

\newcommand{\bwS}{\bw_{\star}}

\newcommand{\nDegP}{\ensuremath{k}}
\newcommand{\Pk}{\ensuremath{\mathcal{P}^{\nDegP}}}

\newcommand{\PkPone}{\ensuremath{\mathcal{P}^{\nDegP+1}}}

\newcommand{\Ga}[2]{\Gamma_{\!\text{#1}}^{#2}}
\newcommand{\hGa}[1]{\hat{\Gamma}_{\!\text{#1}}}
\newcommand{\hOm}{\hat{\Omega}}
\newcommand{\Su}{\Upsilon}

\newcommand{\numelOmega}{\ensuremath{\texttt{n}_{\Omega_e}}}

\newcommand{\meshDomain}{\ensuremath{\Omega}_{\text{o}}}

\newcommand{\inSet}[2]{={#1},\dotsc,{#2}}
\newcommand{\inSetTwo}[2]{={#1},{#2}}

\newcommand{\bigSum}[3]{\bigcup_{{#1}={#2}}^{#3}}

\newcommand{\bigSumE}[3]{\sum_{{#1}={#2}}^{#3}}
\newcommand{\Nitsche}{\eta}

\newcommand{\nSu}  {\ensuremath{\texttt{n}_{\Su}}}
\newcommand{\nOmega}  {\ensuremath{\texttt{n}_{\partial\Omega}}}

\newcommand{\dirData}{\bm{u}_{\text{D}}}

\newcommand{\jump}[1]{\llbracket #1\rrbracket}
\newcommand{\mean}[1]{\{ #1\}}

\usepackage{scalerel,stackengine}
\stackMath
\newcommand\reallyhat[1]{%
\savestack{\tmpbox}{\stretchto{%
  \scaleto{%
    \scalerel*[\widthof{\ensuremath{#1}}]{\kern-.6pt\bigwedge\kern-.6pt}%
    {\rule[-\textheight/2]{1ex}{\textheight}}
  }{\textheight}%
}{0.5ex}}%
\stackon[1pt]{#1}{\tmpbox}%
}

\newenvironment{keywords}{\begin{quote}\emph{\textbf{Keywords:}}}{\end{quote}}

\theoremstyle{definition}

\newtheoremstyle{remark-AH}
                            {}{}
                            {\narrower}
                            {}
                            {\bfseries}
                            {}
                            {.5em}
                            {\thmname{#1}\thmnumber{ #2}\thmnote{ [#3]}}
\theoremstyle{remark-AH}
\newtheorem{remark}{Remark}

\usepackage{tikz}
\usetikzlibrary{fadings}
\usetikzlibrary{patterns}
\usetikzlibrary{shadows.blur}
\usetikzlibrary{shapes}

\newcommand{\kappaGlob}{\kappa^{\text{glob}}}
\newcommand{\kappaLoc}{\kappa^{\text{loc}}}

\graphicspath{ {figs/} {ImagesImmersed/} {ImagesImmersed/badlyCutCells/} {ImagesImmersed/badlyCutFaces/}{ImagesImmersed/couette/} {ImagesImmersed/microchannel/} {interfaceImages/bubbleEq/} {interfaceImages/multiBubbles/} }

\begin{document}
\title{An unfitted high-order HDG method for two-fluid Stokes flow with exact NURBS geometries}

\author{
\renewcommand{\thefootnote}{\arabic{footnote}}
			Stefano Piccardo\footnotemark[1]\textsuperscript{ \ ,}\footnotemark[3]\textsuperscript{ \ ,}\footnotemark[4] , \ 
			Matteo Giacomini\footnotemark[1]\textsuperscript{ \ ,}\footnotemark[2]\textsuperscript{ \ ,}*  \  and
			Antonio Huerta\footnotemark[1]\textsuperscript{ \ ,}\footnotemark[2]
}

\date{}
\maketitle

\renewcommand{\thefootnote}{\arabic{footnote}}

\footnotetext[1]{Laboratori de C\`alcul Num\`eric (LaC\`aN), ETS de Ingenier\'ia de Caminos, Canales y Puertos, Universitat Polit\`ecnica de Catalunya, Barcelona, Spain.}
\footnotetext[2]{Centre Internacional de M\`etodes Num\`erics en Enginyeria (CIMNE), Barcelona, Spain.}
\footnotetext[3]{CERMICS, Ecole des Ponts, 77455 Marne-la-Vallée, France.}
\footnotetext[4]{INRIA, 2 rue Simone Iff, 75589 Paris, France.
\vspace{5pt}\\
* Corresponding author: Matteo Giacomini. \textit{E-mail:} \texttt{matteo.giacomini@upc.edu}
}

\begin{abstract}
A high-order, degree-adaptive hybridizable discontinuous Galerkin (HDG) method is presented for two-fluid incompressible Stokes flows, with boundaries and interfaces described using NURBS.
The NURBS curves are embedded in a fixed Cartesian grid, yielding an unfitted HDG scheme capable of treating the exact geometry of the boundaries/interfaces, circumventing the need for fitted, high-order, curved meshes.
The framework of the NURBS-enhanced finite element method (NEFEM) is employed for accurate quadrature along immersed NURBS and in elements cut by NURBS curves.
A Nitsche's formulation is used to enforce Dirichlet conditions on embedded surfaces, yielding unknowns only on the mesh skeleton as in standard HDG, without introducing any additional degree of freedom on non-matching boundaries/interfaces.
The resulting unfitted HDG-NEFEM method combines non-conforming meshes, exact NURBS geometry and high-order approximations to provide high-fidelity results on coarse meshes, independent of the geometric features of the domain.
Numerical examples illustrate the optimal accuracy and robustness of the method, even in the presence of badly cut cells or faces, and its suitability to simulate microfluidic systems from CAD geometries.
\end{abstract}


\newpage
%
\section{Introduction}
The simulation of multi-fluid systems and multi-phase flows entails a series of numerical challenges related to the different flow features that may arise from the interaction of different materials and phases\ \cite{Reusken-GR-11,Diez-ZD-09,HEIB13}. In this context, high-order methods\ \cite{Wang2013} have received increasing attention in recent years, given their suitability to devise high-fidelity solvers for complex flow problems\ \cite{Riviere-FKR-18,CostaSole-CRS-19}.

Despite their proven superiority in terms of accuracy and computational efficiency\ \cite{Cockburn-KSC:11,AA-HARP:13,May-WBMS-14}, the adoption of high-order methods (e.g., discontinuous Galerkin --DG--) is still limited outside of the academic environment, mainly due to their reduced robustness with respect to low-order methods and the challenges to automatically generate high-order, curved meshes. The latter issue is particularly critical since the error introduced by the geometric approximation of a computer aided design (CAD) model can be responsible for the appearance of non-physical effects in the numerical solution, as reported in the literature\ \cite{Bassi-BR-97,Sevilla-SFH-08}.
To remedy this issue, the NURBS-enhanced finite element method (NEFEM) was proposed in\ \cite{Sevilla-SFH-08-IJNME,Sevilla-SFH-11}, offering a seamless integration of the standard finite element method (FEM) with a CAD representation of the boundary during simulation.  Whilst NEFEM provides a flexible computing environment, allowing for a simple treatment of trimmed and degenerate NURBS,  without the need for complex, manual de-featuring procedures, the automatic generation of meshes suitable for NEFEM represents an open field of investigation, with some promising recent contributions\ \cite{Sevilla-SRH-16,Sevilla-ZLSHM-24}.

An alternative approach to ease the difficulties of geometry treatment is represented by immersed boundary methods\ \cite{Iaccarino-MI-05}.
Many \emph{unfitted}, \emph{embedded} or \emph{immersed} methods have been proposed in the literature, spanning the unfitted finite element method\ \cite{BarEl:87} coupled with Nitsche's method\ \cite{HanHan:02,HanHan:04}, the fictitious domain method\ \cite{GlowTP:94,GlowTP2:94},  the cutFEM\ \cite{BurHan:10,BurHan:12,Burman-BCHLM:15} and cutDG\ \cite{Gurkan-GM-19} methods,  the finite cell method\ \cite{Rank-PDR-07}, the aggregated FEM\ \cite{Badia-BVM-18} and the shifted boundary method\ \cite{Main:18}.
It is well known that unfitted methods tend to suffer in the presence of small cuts in the mesh, possibly leading to ill-conditioned systems\ \cite{Brummelen-PVZB-17}. To remedy this issue, several techniques have been proposed in the literature, including ghost penalty\ \cite{Burman:10}, cell agglomeration\ \cite{JohLa:13} and element extension\ \cite{Navarro:21}. For a detailed discussion about these topics, interested readers are referred to\ \cite{Hansbo-BCHLM-15,Badia-PVBLB-23}.

Whilst well-established for low-order functional and geometric approximations,  unfitted methods have been less studied in a high-order context featuring exact geometry.  
Extensions towards high order have been proposed for the finite cell method\ \cite{Schillinger-SR-15},  the unfitted FEM using high-order geometrical maps\ \cite{Lehre2:16} and the unfitted DG method with cell agglomeration\ \cite{Kummer-MKKO-17}. More recently,  high-order versions were presented for ghost penalty\ \cite{Larson-LZ-20},  shifted boundary method\ \cite{Scovazzi-ACS-22} and aggregated FEM\ \cite{Badia-BNV-22}.
Concerning exact geometry treatment, unfitted approaches employing B-Splines and NURBS were discussed in\ \cite{Rank-SR-11,Rank-RSBVR-13} using the finite cell method, in\ \cite{Legrain-13} via NEFEM, in\ \cite{Bazilevs-KHSEABSH-15,Reali-HVQCARB-19} based on the immersogeometric framework and in\ \cite{Sevilla-MSZRT-15,Navarro:21} for the Cartesian grid FEM and DG methods.

The goal of this work is to devise an unfitted method merging the advantages of high-order functional approximations and exact treatment of the NURBS geometry. 
Stemming from the high-order hybridizable discontinuous Galerkin (HDG) rationale proposed in\ \cite{Jay-CGL:09}, a high-fidelity solver is devised for Stokes flows featuring immersed boundaries, as well as interface problems with two fluids.  Owing to the hybridization procedure, this allows to significantly reduce the globally coupled degrees of freedom of the problem, eliminating all element-based unknowns\ \cite{AA-HARP:13}.
In this work, boundaries and interfaces are represented exactly by means of NURBS\ \cite{Piegl-Tiller}, which are incorporated in the HDG solver via the NEFEM paradigm\ \cite{SevHu:18}.
The main novelties of the proposed approach are: (i) the flexibility of treating exact NURBS geometry within non-matching grids, thus eliminating the difficulties of high-order mesh generation; (ii) the suitability of HDG to construct high-order approximations with immersed boundaries, without the need to introduce additional unknowns along the contours/interfaces; (iii) the capability of devising a degree-adaptive procedure, exploiting the exact representation of the geometry via NEFEM and the superconvergent properties of the HDG primal variable.
The resulting unfitted high-order HDG-NEFEM method has close connections with other techniques proposed in recent years. 
Whilst it shares the definition of face unknowns with the unfitted HDG method\ \cite{Dong-DWXW-16,Nguyen-MNS-22}, the extended HDG method\ \cite{Gurkan-GKF:17,Gurkan-GKF:19}, the hybridized cutFEM\ \cite{Burman-BEHLL-19} and the unfitted hybrid high order method\ \cite{Ern-BE-18,BuCDE:21,BuDeE:21,PiErn:23}, the uniqueness of the proposed approach lies in the definition of unfitted exact geometries via NURBS.
%
%
Moreover, this work inherits from\ \cite{Navarro:21} the element extension strategy employed to handle badly cut cells,  whereas, differently from it, it exploits hybridization to reduce the number of globally coupled degrees of freedom in the system.

The rest of this article is organized as follows.
After introducing the two-fluid Stokes interface problem in\ Section\ \ref{sc:problem}, Section\ \ref{sc:problemDisc} describes the unfitted HDG formulation, accounting for both immersed boundaries and immersed interfaces, and a degree-adaptive procedure.
In section\ \ref{sec:quadDef}, technical details are discussed concerning the treatment of unfitted interfaces, the NEFEM numerical quadrature in the presence of elements cut by NURBS and the element extension strategy to handle badly cut cells.
The numerical results are reported in Section\ \ref{sc:numSimu}: the robustness and accuracy of the method are assessed by means of a series of benchmarks with unfitted boundaries and unfitted interfaces and two applications to microfluidic systems.
Finally, Section\ \ref{sc:Conclusions} summarizes the results of this work and two appendices provide technical details for the use of NURBS in the presented unfitted HDG-NEFEM framework.

%
\section{Problem statement}
\label{sc:problem}
Let $\Omega\subset \RR^{\nsd}$ be an open bounded domain with boundary $\partial\Omega$ and $\nsd = 2$ be the number of spatial dimensions. The boundary $\partial\Omega$ is composed of two disjoint parts, the Dirichlet portion $\Ga{D}{}$,  and the Neumann portion $\Ga{N}{}$.  Formally, $\partial\Omega = \overline{\Ga{D}{}}\cup\overline{\Ga{N}{}}$ such that $\Ga{D}{}\cap\Ga{N}{} = \emptyset$.

Suppose also that $\Omega$ is split by a fixed interface $\Su$ into two disjoint subdomains, each occupied by an immiscible incompressible Stokes fluid. That is, $\overline{\Omega} =  \overline{\Omega^1}\cup\overline{\Omega^2}$, $\Omega^1\cap\Omega^2=\emptyset$, and $\Su = \overline{\Omega^1}\cap\overline{\Omega^2}$. Note that, for each fluid, the boundary is $\partial\Omega^i=\Ga{D}{i} \cup \Ga{N}{i}\cup\Su$ with $\Ga{D}{i}\cap\Su = \emptyset$, and $\Ga{N}{i}\cap\Su = \emptyset$, for $i=1,2$. Moreover, define the domain $\hat{\Omega}$ as $\hat{\Omega} := \Omega^1\cup\Omega^2$ and, similarly, $\hGa{D} := \Ga{D}{1}\cup\Ga{D}{2}$ and $\hGa{N} := \Ga{N}{1}\cup\Ga{N}{2}$. The interface $\Su$ and the boundary $\partial\Omega$ are represented using NURBS and are composed by $\nSu$ and $\nOmega$ NURBS curves, respectively. That is, 
\begin{equation}\label{eq:nurbsDef}
   \Su := \bigSum{j}{1}{\nSu} \mathcal{\pmb{C}}^j_{\Su}([0,1]),\qquad \partial\Omega := \bigSum{j}{1}{\nOmega} \mathcal{\pmb{C}}^j_{\partial\Omega}([0,1]),
\end{equation}
where $\mathcal{\pmb{C}}^j_{\diamond} : \lambda \mapsto \mathcal{\pmb{C}}^j_{\diamond}(\lambda)$ is a generic NURBS curve\ \cite{Piegl-Tiller,Rogers} defined in the parametric domain $\lambda\in[0,1]$.

The problem aims to find the unknown velocity and pressure fields, $(\bu,p)\in[\sobo(\hOm)]^{\nsd}\times\eltwo(\hOm)$, whose restrictions to each subdomain $\Omega^i$ are $(\bu,p)\vert_{\Omega^i}=(\bu^i,p^i)\in[\sobo(\Omega^i)]^{\nsd}\times\eltwo(\Omega^i)$, $i=1,2$, and where 
$
 [\sobo(\hOm)]^{\nsd} := \bigl\{ \bv\in [\eltwo(\hOm)]^{\nsd} \mid \bv\vert_{\Omega^i}\in[\sobo(\Omega^i)]^{\nsd},\; \forall i=1,2 \bigr\} 
$. 
Given the velocity profile $\dirData\in[\sobo[\frac{1}{2}](\hGa{D}{})]^{\nsd}$ imposed on the Dirichlet portion of the boundary and the Neumann datum $\bt\in[\eltwo(\hGa{N}{})]^{\nsd}$ applied on $\hGa{N}{}$,
the strong form of the problem can be written as find  $(\bu,p)\in[\sobo(\hOm)]^{\nsd}\times\eltwo(\hOm)$ such that
\begin{equation}\label{eq:probstat}
\begin{split}
 &\left\{
 \begin{aligned}
 -\grad\cdot(\mu\grad\bu - p \Insd) &= \bm{s}                    &&\text{in $\hOm$,}\\
   \grad\cdot\bu &= 0                                                         &&\text{in $\hOm$,}\\
   \bu &= \dirData                                                                  &&\text{on $\hGa{D}$,}\\
   \bigl(\mu\grad\bu - p \Insd \bigr)\,\bn &= \bt   &&\text{on $\hGa{N}$,}\\
 \end{aligned}\right. 
 \\[1ex]
 &\left\{\begin{aligned}
     \bu^1&= \bu^2 \\
   \jump{(\mu\grad\bu {-} p \Insd )\bn}
    &= \gamma(\Div\bn^1)\bn^1{-}(\Insd{-}\bn^1{\otimes}\bn^1)\grad\gamma
 \end{aligned}\right.\quad \text{on $\!\Su$,}
 \end{split}
\end{equation}
where $\mu=\mu^i>0$ in $\Omega^i$, $i=1,2$, is the piecewise dynamic viscosity of each fluid, assumed constant in $\Omega^i$, $\bm{s}\in[\eltwo(\hOm)]^{\nsd}$ is the volumetric source term, $\bn^i$ is the unit normal vector exiting domain $\Omega^i$, and $\gamma$ is the surface tension coefficient. Finally, the last two equations impose the continuity of velocity and equilibrium of forces along the interface $\Su$.

Note also the use of the \emph{jump} operator, $\jump{\odot}$, which follows the definition introduced in \cite{AdM-MFH:08}. Along any $(\nsd{-}1)$-dimensional manifold the \emph{jump} operator sums the values of a generic quantity $\odot$ from the left and from the right, namely,
\begin{equation}\label{eq:jump}
  \jump{\odot} := \odot_l + \odot_r ,
\end{equation}
in this case, since the discontinuity is along $\Su$, then $l{=}1$ and $r{=}2$ or viceversa. Note that the above definition of the jump operator always involves the outward unit normal to a surface. 

\begin{remark}[Surface tension]
Note that in the last equation in~\eqref{eq:probstat} the choice of $\bn^1$ in the right-hand-side of the equation is arbitrary because
\begin{equation*}
 \gamma(\Div\bn^1)\bn^1{-}(\Insd{-}\bn^1{\otimes}\bn^1)\grad\gamma = \gamma(\Div\bn^2)\bn^2{-}(\Insd{-}\bn^2{\otimes}\bn^2)\grad\gamma .
\end{equation*}
Moreover, in most cases, the surface tension, $\gamma$, is assumed constant, thus, $\grad\gamma=\bm{0}$, and, consequently, the interface shear stress is continuous. Without loss of generality, this assumption is retained in the rest of the paper.
\end{remark}

\begin{remark}[One-fluid problem]
When $\Su = \emptyset$, the interface conditions disappear, and one-fluid problem is solved in domain $\Omega^1 \equiv \Omega$.
\end{remark}

\begin{remark}[Uniqueness of pressure]
When $\hGa{N} = \emptyset$, the Stokes interface problem is solvable up to a global additive constant on the pressure, which we fix by imposing
\begin{equation*}
 \int_{\hOm} p \, \mathrm{d}{\hOm} =  \int_{\hOm} p^{\text{ref}}  \mathrm{d}{\hOm},
\end{equation*}
where $p^{\text{ref}}$ is a reference pressure value, typically equal to zero. At the discrete level, see, e.g.,~\cite[Chapter 17]{Quarteroni-17},  an alternative strategy to eliminate this indeterminacy relies on prescribing the pressure at one point, as commonly performed by finite volume methods\ \cite{Jasak-PhD-96,RS-SGH:2018_FCFV1}.
%
\label{rm:pressureC0}
\end{remark}

\begin{remark}[Multi-fluid interface problem]
Formulation~\eqref{eq:probstat} can be readily extended to address multi-fluid problems, where $\Omega$ is partitioned into a generic number of connected, open, bounded sets, each occupied by an immiscible, incompressible Stokes fluid. 
\end{remark}

%
\section{Unfitted HDG formulation}
\label{sc:problemDisc}
Let $\meshDomain$ be a shape-regular polyhedral domain containing $\Omega$, i.e., $\Omega \subseteq \meshDomain $. Then, consider a shape-regular mesh composed of $\numel$ disjoint (open) subdomains $\Omega_e$ such that $\meshDomain$ is exactly covered, i.e.,
\begin{equation*}
 \overline{\meshDomain} := \bigSum{e}{1}{\numel} \overline{\Omega_e}. 
\end{equation*}
%
The internal boundaries $\partial\Omega_e$ of the mesh elements $\Omega_e, \ e=1,\ldots,\numel$ define the \emph{internal skeleton}, $\Gamma$, as
\begin{equation*}
 \Gamma := \Bigg[ \bigSum{e}{1}{\numel} \partial\Omega_e \Bigg]\setminus\bigl(\hGa{D}\cup\hGa{N}\cup\partial\meshDomain\bigr). 
\end{equation*} 
We denote by $\Gamma_{\!\! f}$ a generic internal face of $\Gamma$, and by $\nface$ the total number of internal faces, i.e., $\Gamma  = \bigSum{f}{1}{\nface}\Gamma_{\!\! f}$. 
It is important to recall that the interface, $\Su$, and the external boundary, $\hGa{D}$ and $\hGa{N}$, do not need to align with the mesh.
Moreover, for all $e\inSet{1}{\numel}$,  define $\Omega_e^i$ the region of $\Omega_e$ that belongs to the fluid indexed by $i$, that is,
\begin{equation*}
    \Omega_e^{i} := \Omega_e \cap \Omega^{i}, \qquad \forall i=1,2.
\end{equation*}
Similarly, define $\partial\Omega_e^{i} := \partial\Omega_e \cap(\overline{ \Omega^{i}} \setminus \Su)$ as the portion of the boundary $\partial\Omega_e$ belonging to fluid $i$. 

On this broken domain, problem~\eqref{eq:probstat} can be written in mixed form as follows:
\begin{subequations}\label{eq:brokenprbsta}
\begin{align}
 &\left\{
 \begin{aligned}
   \bL^i_e + \sqrt{\mu^i}\grad\bu^i_e &= \bm{0}                                             &&\text{in $\Omega_e^i$,}\\ 
  \grad\cdot(\sqrt{\mu^i}\bL^i_e) + \grad p^i_e &= \bm{s}                              &&\text{in $\Omega_e^i$,}\\
  \grad\cdot\bu^i_e &= 0                                                                               &&\text{in $\Omega_e^i$,}\\
   \bu^i_e &= \dirData                                                                                      &&\text{on $\Ga{D}{i}\cap\overline{\Omega_e}$,}\\
   \bigl(\sqrt{\mu^i}\bL^i_e + p^i_e \Insd \bigr)\,\bn^i_e &= -\bt  &&\text{on $\Ga{N}{i}\cap\overline{\Omega_e}$,}\\
 \end{aligned}\right.  
 && \text{for } e\inSet{1}{\numel} \atop {\hspace{-28pt} \text{ and } i=1,2},
 \label{eq:PSLocal}
 \\[1ex]
&\left\{\begin{aligned}
 \bu^1_e&=\bu^2_e \\
 \jump{(\sqrt{\mu}\bL_e {+} p_e \Insd )\bn_e}
    &= 
    {-}\gamma(\Div\bn^1_e)\bn^1_e
 \end{aligned}\right. \text{ on $\Su{\cap}\Omega_e$,} 
 &&\text{for $e=1,\dotsc,\numel$,}
 \label{eq:PSInter}
 \\[1ex]
&\left\{
  \begin{aligned}
    \jump{\bu^i \otimes \bn^i}                     &= \bm{0}  \\
    \jump{(\sqrt{\mu^i} \bL^i + p^i \Insd)\,\bn^i} &= \bm{0}  \\
 \end{aligned}\right. \quad\text{on $\Gamma\cap\Omega^i$,} 
 &&\text{for $i=1,2$,}
 \label{eq:PSTrans}
 \end{align}
\end{subequations}
where $\bL:= - \sqrt{\mu} \grad\bu$ is the mixed variable so that $\bL^{i}_e= - \sqrt{\mu^{i}} \grad\bu^{i}_e$ is its restriction to $\Omega^{i}_e$, for $e\inSet{1}{\numel}$ and $i=1,2$. Note that the \emph{jump} operator acts over the interface, $\Su$, in~\eqref{eq:PSInter} whereas this operator acts along $\Gamma\cap\Omega^i$ for  $i=1,2$ in~\eqref{eq:PSTrans}.

The interface conditions stated in~\eqref{eq:PSInter} are those presented in~\eqref{eq:probstat}, rewritten element-by-element to ensure the enforcement of the desired continuity conditions. 

Equations~\eqref{eq:PSLocal} and~\eqref{eq:PSTrans} are, respectively, the so-called \emph{local problems} and \emph{transmission conditions} common in HDG, see \cite{MG-GSH:20}. It is worth noting that here, following \cite{RS-SH:16}, the HDG formulation imposes the Neumann boundary condition in the local problem instead of the global one as it is usually done in HDG. This leads to a marginally smaller discrete global problem but, more importantly, avoids the definition of hybrid unknowns along the Neumann boundary (regardless of having a fitted or unfitted Neumann boundary). The transmission conditions impose, respectively, continuity of velocity and normal flux across the skeleton of the mesh. Finally, in order to solve the local problems, HDG methods introduce the so-called \emph{hybrid velocity}: this is a face-based variable $\bhu^i$ that represents the trace of the solution on each face of $\Gamma\cap\Omega^i$ so that each local (element-by-element) problem described by~\eqref{eq:PSLocal} is completed by adding 
%
\begin{equation}\label{eq:DirichletHybrid}  
  \bu_e^i = \bhu^i \text{ on $\Gamma \cap \partial\Omega_e^i$, for $e\inSet{1}{\numel}$ and $i=1,2$.}
\end{equation}

To simplify the presentation,  consider three types of elements: (i) \emph{standard HDG elements}, which are those not cut by $\Su$ or $\partial\Omega$, that is, $\Su\cap \Omega_e = \emptyset$ and $\partial\Omega\cap \Omega_e = \emptyset$; (ii)  \emph{immersed boundary elements} cut by $\partial\Omega$ but not by $\Su$, that is, $\Su\cap \Omega_e = \emptyset$ and $\partial\Omega\cap \Omega_e \neq \emptyset$; and (iii) \emph{interface elements} cut by $\Su$ but not by $\partial\Omega$, that is, $\Su\cap \Omega_e \neq \emptyset$ and $\partial\Omega\cap \Omega_e = \emptyset$. 
Note that all the defined elements, (i), (ii), and (iii), also encompass the case in which the external boundary is aligned with the mesh skeleton, which means $\partial\Omega_e\cap\partial\Omega\neq \emptyset$.
%
\begin{figure}[t]
    \centering
\tikzset{every picture/.style={line width=0.75pt}} 

\begin{tikzpicture}[x=0.75pt,y=0.75pt,yscale=-1,xscale=1]

\draw  [draw opacity=0][fill={rgb, 255:red, 126; green, 211; blue, 33 }  ,fill opacity=1 ] (158,146) -- (196,146) -- (196,184) -- (158,184) -- cycle ;
\draw  [draw opacity=0][fill={rgb, 255:red, 245; green, 166; blue, 35 }  ,fill opacity=1 ] (310,184) -- (348,184) -- (348,222) -- (310,222) -- cycle ;
\draw [line width=0.75]    (167,102) -- (145.12,80.12) ;
\draw [shift={(143,78)}, rotate = 45] [fill={rgb, 255:red, 0; green, 0; blue, 0 }  ][line width=0.08]  [draw opacity=0] (8.93,-4.29) -- (0,0) -- (8.93,4.29) -- cycle    ;
\draw [color={rgb, 255:red, 74; green, 144; blue, 226 }  ,draw opacity=1 ][line width=2.25]    (261,199) .. controls (133,196) and (131,119) .. (180,91) .. controls (221,66) and (317,74) .. (262,101) .. controls (207,128) and (337,207) .. (261,199) -- cycle ;
\draw  [line width=1.5]  (220.85,85) -- (230.9,77.06) -- (218.45,72.29) ;
\draw [line width=0.75]    (150,225) -- (126.96,251.73) ;
\draw [shift={(125,254)}, rotate = 310.76] [fill={rgb, 255:red, 0; green, 0; blue, 0 }  ][line width=0.08]  [draw opacity=0] (8.93,-4.29) -- (0,0) -- (8.93,4.29) -- cycle    ;
\draw  [draw opacity=0] (120,32) -- (349,32) -- (349,261) -- (120,261) -- cycle ; \draw  [color={rgb, 255:red, 128; green, 128; blue, 128 }  ,draw opacity=1 ] (120,32) -- (120,261)(158,32) -- (158,261)(196,32) -- (196,261)(234,32) -- (234,261)(272,32) -- (272,261)(310,32) -- (310,261)(348,32) -- (348,261) ; \draw  [color={rgb, 255:red, 128; green, 128; blue, 128 }  ,draw opacity=1 ] (120,32) -- (349,32)(120,70) -- (349,70)(120,108) -- (349,108)(120,146) -- (349,146)(120,184) -- (349,184)(120,222) -- (349,222)(120,260) -- (349,260) ; \draw  [color={rgb, 255:red, 128; green, 128; blue, 128 }  ,draw opacity=1 ]  ;
\draw  [line width=2.25]  (120,146) .. controls (120,83.04) and (171.04,32) .. (234,32) .. controls (296.96,32) and (348,83.04) .. (348,146) .. controls (348,208.96) and (296.96,260) .. (234,260) .. controls (171.04,260) and (120,208.96) .. (120,146) -- cycle ;
\draw [line width=0.75]    (201,188) -- (210.94,161.8) ;
\draw [shift={(212,159)}, rotate = 110.77] [fill={rgb, 255:red, 0; green, 0; blue, 0 }  ][line width=0.08]  [draw opacity=0] (8.93,-4.29) -- (0,0) -- (8.93,4.29) -- cycle    ;
\draw  [draw opacity=0][fill={rgb, 255:red, 208; green, 2; blue, 27 }  ,fill opacity=1 ] (272,108) -- (310,108) -- (310,146) -- (272,146) -- cycle ;
\draw  [draw opacity=0][fill={rgb, 255:red, 208; green, 2; blue, 27 }  ,fill opacity=1 ] (196,108) -- (234,108) -- (234,146) -- (196,146) -- cycle ;

\draw (166,115) node [anchor=north west][inner sep=0.75pt]  [font=\large] [align=left] {$\displaystyle \Omega ^{1}$};
\draw (244,226) node [anchor=north west][inner sep=0.75pt]  [font=\large] [align=left] {$\displaystyle \Omega ^{2}$};
\draw (216,51) node [anchor=north west][inner sep=0.75pt]  [font=\large,color={rgb, 255:red, 0; green, 0; blue, 0 }  ,opacity=1 ] [align=left] {$\displaystyle \Upsilon $};
\draw (350,35) node [anchor=north west][inner sep=0.75pt]  [font=\large] [align=left] {$\displaystyle \meshDomain $};
\draw (315,230) node [anchor=north west][inner sep=0.75pt]  [font=\large] [align=left] {$\displaystyle \partial \Omega$};
\draw (160,73) node [anchor=north west][inner sep=0.75pt]  [font=\large] [align=left] {$\displaystyle n^{1}$};
\draw (136,237) node [anchor=north west][inner sep=0.75pt]  [font=\large] [align=left] {$\displaystyle n^{2}$};
\draw (213,158) node [anchor=north west][inner sep=0.75pt]  [font=\large] [align=left] {$\displaystyle n^{2}$};

\end{tikzpicture}
    \caption[Schematic representation of the geometry and its discretization.]{Schematic representation of the geometry and its discretization. Computational domain $\meshDomain$ partitioned into $6\times6$ square elements. The external physical boundary $\partial\Omega$ and the interface $\Su$ not aligned with the mesh skeleton. Standard HDG elements in red, immersed boundary element (cut by $\partial\Omega$) in orange, interface element (cut by $\Su$) in green.}
    \label{fg:cut_cell_ex}
\end{figure}
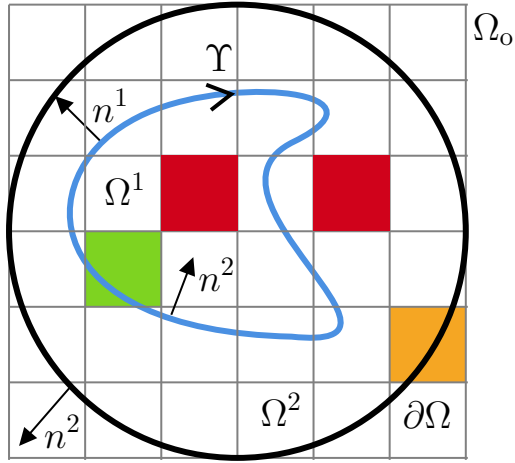

The method is henceforth presented assuming that each element is split, at most, into two different regions. Figure\ \ref{fg:cut_cell_ex} illustrates one possible configuration of such a case. Note that this assumption is introduced exclusively to simplify the readability of the proposed unfitted HDG method and does not represent a limitation of the presented methodology. Indeed,  the numerical results in Section\ \ref{sc:microMix} showcase the capability of the method to handle mesh elements subdivided in more than two subregions. Finally, the technical steps required to extend the discussed formulation to the case of multiple subregions are presented in Appendix\ \ref{ap:intersect}.
%
%

\begin{figure}[!h]
    \centering
    \subfloat[Standard HDG element.]{
\tikzset{every picture/.style={line width=0.75pt}} 

\begin{tikzpicture}[x=0.75pt,y=0.75pt,yscale=-0.6,xscale=0.6]

\draw   (151,100) -- (351,100) -- (351,300) -- (151,300) -- cycle ;
\draw    (151,351) -- (350,351) ;
\draw    (101,301) -- (101,100) ;
\draw  [color={rgb, 255:red, 0; green, 0; blue, 0 }  ,draw opacity=1 ][fill={rgb, 255:red, 208; green, 2; blue, 27 }  ,fill opacity=1 ] (242,351) .. controls (242,346.31) and (245.81,342.5) .. (250.5,342.5) .. controls (255.19,342.5) and (259,346.31) .. (259,351) .. controls (259,355.69) and (255.19,359.5) .. (250.5,359.5) .. controls (245.81,359.5) and (242,355.69) .. (242,351) -- cycle ;
\draw  [color={rgb, 255:red, 0; green, 0; blue, 0 }  ,draw opacity=1 ][fill={rgb, 255:red, 208; green, 2; blue, 27 }  ,fill opacity=1 ] (341.5,351) .. controls (341.5,346.31) and (345.31,342.5) .. (350,342.5) .. controls (354.69,342.5) and (358.5,346.31) .. (358.5,351) .. controls (358.5,355.69) and (354.69,359.5) .. (350,359.5) .. controls (345.31,359.5) and (341.5,355.69) .. (341.5,351) -- cycle ;
\draw  [color={rgb, 255:red, 0; green, 0; blue, 0 }  ,draw opacity=1 ][fill={rgb, 255:red, 208; green, 2; blue, 27 }  ,fill opacity=1 ] (142.5,351) .. controls (142.5,346.31) and (146.31,342.5) .. (151,342.5) .. controls (155.69,342.5) and (159.5,346.31) .. (159.5,351) .. controls (159.5,355.69) and (155.69,359.5) .. (151,359.5) .. controls (146.31,359.5) and (142.5,355.69) .. (142.5,351) -- cycle ;
\draw  [color={rgb, 255:red, 0; green, 0; blue, 0 }  ,draw opacity=1 ][fill={rgb, 255:red, 74; green, 144; blue, 226 }  ,fill opacity=1 ] (242.5,100) .. controls (242.5,95.31) and (246.31,91.5) .. (251,91.5) .. controls (255.69,91.5) and (259.5,95.31) .. (259.5,100) .. controls (259.5,104.69) and (255.69,108.5) .. (251,108.5) .. controls (246.31,108.5) and (242.5,104.69) .. (242.5,100) -- cycle ;
\draw  [color={rgb, 255:red, 0; green, 0; blue, 0 }  ,draw opacity=1 ][fill={rgb, 255:red, 74; green, 144; blue, 226 }  ,fill opacity=1 ] (342.5,100) .. controls (342.5,95.31) and (346.31,91.5) .. (351,91.5) .. controls (355.69,91.5) and (359.5,95.31) .. (359.5,100) .. controls (359.5,104.69) and (355.69,108.5) .. (351,108.5) .. controls (346.31,108.5) and (342.5,104.69) .. (342.5,100) -- cycle ;
\draw  [color={rgb, 255:red, 0; green, 0; blue, 0 }  ,draw opacity=1 ][fill={rgb, 255:red, 74; green, 144; blue, 226 }  ,fill opacity=1 ] (142.5,200) .. controls (142.5,195.31) and (146.31,191.5) .. (151,191.5) .. controls (155.69,191.5) and (159.5,195.31) .. (159.5,200) .. controls (159.5,204.69) and (155.69,208.5) .. (151,208.5) .. controls (146.31,208.5) and (142.5,204.69) .. (142.5,200) -- cycle ;
\draw  [color={rgb, 255:red, 0; green, 0; blue, 0 }  ,draw opacity=1 ][fill={rgb, 255:red, 74; green, 144; blue, 226 }  ,fill opacity=1 ] (242.5,200) .. controls (242.5,195.31) and (246.31,191.5) .. (251,191.5) .. controls (255.69,191.5) and (259.5,195.31) .. (259.5,200) .. controls (259.5,204.69) and (255.69,208.5) .. (251,208.5) .. controls (246.31,208.5) and (242.5,204.69) .. (242.5,200) -- cycle ;
\draw  [color={rgb, 255:red, 0; green, 0; blue, 0 }  ,draw opacity=1 ][fill={rgb, 255:red, 74; green, 144; blue, 226 }  ,fill opacity=1 ] (342.5,300) .. controls (342.5,295.31) and (346.31,291.5) .. (351,291.5) .. controls (355.69,291.5) and (359.5,295.31) .. (359.5,300) .. controls (359.5,304.69) and (355.69,308.5) .. (351,308.5) .. controls (346.31,308.5) and (342.5,304.69) .. (342.5,300) -- cycle ;
\draw  [color={rgb, 255:red, 0; green, 0; blue, 0 }  ,draw opacity=1 ][fill={rgb, 255:red, 74; green, 144; blue, 226 }  ,fill opacity=1 ] (142.5,300) .. controls (142.5,295.31) and (146.31,291.5) .. (151,291.5) .. controls (155.69,291.5) and (159.5,295.31) .. (159.5,300) .. controls (159.5,304.69) and (155.69,308.5) .. (151,308.5) .. controls (146.31,308.5) and (142.5,304.69) .. (142.5,300) -- cycle ;
\draw  [color={rgb, 255:red, 0; green, 0; blue, 0 }  ,draw opacity=1 ][fill={rgb, 255:red, 74; green, 144; blue, 226 }  ,fill opacity=1 ] (142.5,100) .. controls (142.5,95.31) and (146.31,91.5) .. (151,91.5) .. controls (155.69,91.5) and (159.5,95.31) .. (159.5,100) .. controls (159.5,104.69) and (155.69,108.5) .. (151,108.5) .. controls (146.31,108.5) and (142.5,104.69) .. (142.5,100) -- cycle ;
\draw  [color={rgb, 255:red, 0; green, 0; blue, 0 }  ,draw opacity=1 ][fill={rgb, 255:red, 74; green, 144; blue, 226 }  ,fill opacity=1 ] (342.5,200) .. controls (342.5,195.31) and (346.31,191.5) .. (351,191.5) .. controls (355.69,191.5) and (359.5,195.31) .. (359.5,200) .. controls (359.5,204.69) and (355.69,208.5) .. (351,208.5) .. controls (346.31,208.5) and (342.5,204.69) .. (342.5,200) -- cycle ;
\draw  [color={rgb, 255:red, 0; green, 0; blue, 0 }  ,draw opacity=1 ][fill={rgb, 255:red, 74; green, 144; blue, 226 }  ,fill opacity=1 ] (242.5,300) .. controls (242.5,295.31) and (246.31,291.5) .. (251,291.5) .. controls (255.69,291.5) and (259.5,295.31) .. (259.5,300) .. controls (259.5,304.69) and (255.69,308.5) .. (251,308.5) .. controls (246.31,308.5) and (242.5,304.69) .. (242.5,300) -- cycle ;
\draw    (153,51) -- (352,51) ;
\draw  [color={rgb, 255:red, 0; green, 0; blue, 0 }  ,draw opacity=1 ][fill={rgb, 255:red, 208; green, 2; blue, 27 }  ,fill opacity=1 ] (92.5,200.5) .. controls (92.5,195.81) and (96.31,192) .. (101,192) .. controls (105.69,192) and (109.5,195.81) .. (109.5,200.5) .. controls (109.5,205.19) and (105.69,209) .. (101,209) .. controls (96.31,209) and (92.5,205.19) .. (92.5,200.5) -- cycle ;
\draw  [color={rgb, 255:red, 0; green, 0; blue, 0 }  ,draw opacity=1 ][fill={rgb, 255:red, 208; green, 2; blue, 27 }  ,fill opacity=1 ] (92.5,301) .. controls (92.5,296.31) and (96.31,292.5) .. (101,292.5) .. controls (105.69,292.5) and (109.5,296.31) .. (109.5,301) .. controls (109.5,305.69) and (105.69,309.5) .. (101,309.5) .. controls (96.31,309.5) and (92.5,305.69) .. (92.5,301) -- cycle ;
\draw  [color={rgb, 255:red, 0; green, 0; blue, 0 }  ,draw opacity=1 ][fill={rgb, 255:red, 208; green, 2; blue, 27 }  ,fill opacity=1 ] (92.5,100) .. controls (92.5,95.31) and (96.31,91.5) .. (101,91.5) .. controls (105.69,91.5) and (109.5,95.31) .. (109.5,100) .. controls (109.5,104.69) and (105.69,108.5) .. (101,108.5) .. controls (96.31,108.5) and (92.5,104.69) .. (92.5,100) -- cycle ;
\draw  [color={rgb, 255:red, 0; green, 0; blue, 0 }  ,draw opacity=1 ][fill={rgb, 255:red, 208; green, 2; blue, 27 }  ,fill opacity=1 ] (244,51) .. controls (244,46.31) and (247.81,42.5) .. (252.5,42.5) .. controls (257.19,42.5) and (261,46.31) .. (261,51) .. controls (261,55.69) and (257.19,59.5) .. (252.5,59.5) .. controls (247.81,59.5) and (244,55.69) .. (244,51) -- cycle ;
\draw  [color={rgb, 255:red, 0; green, 0; blue, 0 }  ,draw opacity=1 ][fill={rgb, 255:red, 208; green, 2; blue, 27 }  ,fill opacity=1 ] (144.5,51) .. controls (144.5,46.31) and (148.31,42.5) .. (153,42.5) .. controls (157.69,42.5) and (161.5,46.31) .. (161.5,51) .. controls (161.5,55.69) and (157.69,59.5) .. (153,59.5) .. controls (148.31,59.5) and (144.5,55.69) .. (144.5,51) -- cycle ;
\draw  [color={rgb, 255:red, 0; green, 0; blue, 0 }  ,draw opacity=1 ][fill={rgb, 255:red, 208; green, 2; blue, 27 }  ,fill opacity=1 ] (343.5,51) .. controls (343.5,46.31) and (347.31,42.5) .. (352,42.5) .. controls (356.69,42.5) and (360.5,46.31) .. (360.5,51) .. controls (360.5,55.69) and (356.69,59.5) .. (352,59.5) .. controls (347.31,59.5) and (343.5,55.69) .. (343.5,51) -- cycle ;
\draw    (400,302) -- (400,101) ;
\draw  [color={rgb, 255:red, 0; green, 0; blue, 0 }  ,draw opacity=1 ][fill={rgb, 255:red, 208; green, 2; blue, 27 }  ,fill opacity=1 ] (391.5,101) .. controls (391.5,96.31) and (395.31,92.5) .. (400,92.5) .. controls (404.69,92.5) and (408.5,96.31) .. (408.5,101) .. controls (408.5,105.69) and (404.69,109.5) .. (400,109.5) .. controls (395.31,109.5) and (391.5,105.69) .. (391.5,101) -- cycle ;
\draw  [color={rgb, 255:red, 0; green, 0; blue, 0 }  ,draw opacity=1 ][fill={rgb, 255:red, 208; green, 2; blue, 27 }  ,fill opacity=1 ] (391.5,302) .. controls (391.5,297.31) and (395.31,293.5) .. (400,293.5) .. controls (404.69,293.5) and (408.5,297.31) .. (408.5,302) .. controls (408.5,306.69) and (404.69,310.5) .. (400,310.5) .. controls (395.31,310.5) and (391.5,306.69) .. (391.5,302) -- cycle ;
\draw  [color={rgb, 255:red, 0; green, 0; blue, 0 }  ,draw opacity=1 ][fill={rgb, 255:red, 208; green, 2; blue, 27 }  ,fill opacity=1 ] (391.5,201.5) .. controls (391.5,196.81) and (395.31,193) .. (400,193) .. controls (404.69,193) and (408.5,196.81) .. (408.5,201.5) .. controls (408.5,206.19) and (404.69,210) .. (400,210) .. controls (395.31,210) and (391.5,206.19) .. (391.5,201.5) -- cycle ;

\end{tikzpicture}
    \label{fig:elementsBasesA}
    }\qquad
    \subfloat[Immersed boundary element.]{

\tikzset{every picture/.style={line width=0.75pt}} 

\begin{tikzpicture}[x=0.75pt,y=0.75pt,yscale=-0.6,xscale=0.6]

\draw   (151,100) -- (351,100) -- (351,300) -- (151,300) -- cycle ;
\draw    (151,351) -- (350,351) ;
\draw    (101,301) -- (101,100) ;
\draw  [color={rgb, 255:red, 0; green, 0; blue, 0 }  ,draw opacity=1 ][fill={rgb, 255:red, 208; green, 2; blue, 27 }  ,fill opacity=1 ] (242,351) .. controls (242,346.31) and (245.81,342.5) .. (250.5,342.5) .. controls (255.19,342.5) and (259,346.31) .. (259,351) .. controls (259,355.69) and (255.19,359.5) .. (250.5,359.5) .. controls (245.81,359.5) and (242,355.69) .. (242,351) -- cycle ;
\draw  [color={rgb, 255:red, 0; green, 0; blue, 0 }  ,draw opacity=1 ][fill={rgb, 255:red, 208; green, 2; blue, 27 }  ,fill opacity=1 ] (341.5,351) .. controls (341.5,346.31) and (345.31,342.5) .. (350,342.5) .. controls (354.69,342.5) and (358.5,346.31) .. (358.5,351) .. controls (358.5,355.69) and (354.69,359.5) .. (350,359.5) .. controls (345.31,359.5) and (341.5,355.69) .. (341.5,351) -- cycle ;
\draw  [color={rgb, 255:red, 0; green, 0; blue, 0 }  ,draw opacity=1 ][fill={rgb, 255:red, 208; green, 2; blue, 27 }  ,fill opacity=1 ] (142.5,351) .. controls (142.5,346.31) and (146.31,342.5) .. (151,342.5) .. controls (155.69,342.5) and (159.5,346.31) .. (159.5,351) .. controls (159.5,355.69) and (155.69,359.5) .. (151,359.5) .. controls (146.31,359.5) and (142.5,355.69) .. (142.5,351) -- cycle ;
\draw  [color={rgb, 255:red, 0; green, 0; blue, 0 }  ,draw opacity=1 ][fill={rgb, 255:red, 74; green, 144; blue, 226 }  ,fill opacity=1 ] (242.5,100) .. controls (242.5,95.31) and (246.31,91.5) .. (251,91.5) .. controls (255.69,91.5) and (259.5,95.31) .. (259.5,100) .. controls (259.5,104.69) and (255.69,108.5) .. (251,108.5) .. controls (246.31,108.5) and (242.5,104.69) .. (242.5,100) -- cycle ;
\draw  [color={rgb, 255:red, 0; green, 0; blue, 0 }  ,draw opacity=1 ][fill={rgb, 255:red, 74; green, 144; blue, 226 }  ,fill opacity=1 ] (342.5,100) .. controls (342.5,95.31) and (346.31,91.5) .. (351,91.5) .. controls (355.69,91.5) and (359.5,95.31) .. (359.5,100) .. controls (359.5,104.69) and (355.69,108.5) .. (351,108.5) .. controls (346.31,108.5) and (342.5,104.69) .. (342.5,100) -- cycle ;
\draw  [color={rgb, 255:red, 0; green, 0; blue, 0 }  ,draw opacity=1 ][fill={rgb, 255:red, 74; green, 144; blue, 226 }  ,fill opacity=1 ] (142.5,200) .. controls (142.5,195.31) and (146.31,191.5) .. (151,191.5) .. controls (155.69,191.5) and (159.5,195.31) .. (159.5,200) .. controls (159.5,204.69) and (155.69,208.5) .. (151,208.5) .. controls (146.31,208.5) and (142.5,204.69) .. (142.5,200) -- cycle ;
\draw  [color={rgb, 255:red, 0; green, 0; blue, 0 }  ,draw opacity=1 ][fill={rgb, 255:red, 74; green, 144; blue, 226 }  ,fill opacity=1 ] (242.5,200) .. controls (242.5,195.31) and (246.31,191.5) .. (251,191.5) .. controls (255.69,191.5) and (259.5,195.31) .. (259.5,200) .. controls (259.5,204.69) and (255.69,208.5) .. (251,208.5) .. controls (246.31,208.5) and (242.5,204.69) .. (242.5,200) -- cycle ;
\draw  [color={rgb, 255:red, 0; green, 0; blue, 0 }  ,draw opacity=1 ][fill={rgb, 255:red, 74; green, 144; blue, 226 }  ,fill opacity=1 ] (342.5,300) .. controls (342.5,295.31) and (346.31,291.5) .. (351,291.5) .. controls (355.69,291.5) and (359.5,295.31) .. (359.5,300) .. controls (359.5,304.69) and (355.69,308.5) .. (351,308.5) .. controls (346.31,308.5) and (342.5,304.69) .. (342.5,300) -- cycle ;
\draw  [color={rgb, 255:red, 0; green, 0; blue, 0 }  ,draw opacity=1 ][fill={rgb, 255:red, 74; green, 144; blue, 226 }  ,fill opacity=1 ] (142.5,300) .. controls (142.5,295.31) and (146.31,291.5) .. (151,291.5) .. controls (155.69,291.5) and (159.5,295.31) .. (159.5,300) .. controls (159.5,304.69) and (155.69,308.5) .. (151,308.5) .. controls (146.31,308.5) and (142.5,304.69) .. (142.5,300) -- cycle ;
\draw  [color={rgb, 255:red, 0; green, 0; blue, 0 }  ,draw opacity=1 ][fill={rgb, 255:red, 74; green, 144; blue, 226 }  ,fill opacity=1 ] (142.5,100) .. controls (142.5,95.31) and (146.31,91.5) .. (151,91.5) .. controls (155.69,91.5) and (159.5,95.31) .. (159.5,100) .. controls (159.5,104.69) and (155.69,108.5) .. (151,108.5) .. controls (146.31,108.5) and (142.5,104.69) .. (142.5,100) -- cycle ;
\draw  [color={rgb, 255:red, 0; green, 0; blue, 0 }  ,draw opacity=1 ][fill={rgb, 255:red, 74; green, 144; blue, 226 }  ,fill opacity=1 ] (342.5,200) .. controls (342.5,195.31) and (346.31,191.5) .. (351,191.5) .. controls (355.69,191.5) and (359.5,195.31) .. (359.5,200) .. controls (359.5,204.69) and (355.69,208.5) .. (351,208.5) .. controls (346.31,208.5) and (342.5,204.69) .. (342.5,200) -- cycle ;
\draw  [color={rgb, 255:red, 0; green, 0; blue, 0 }  ,draw opacity=1 ][fill={rgb, 255:red, 74; green, 144; blue, 226 }  ,fill opacity=1 ] (242.5,300) .. controls (242.5,295.31) and (246.31,291.5) .. (251,291.5) .. controls (255.69,291.5) and (259.5,295.31) .. (259.5,300) .. controls (259.5,304.69) and (255.69,308.5) .. (251,308.5) .. controls (246.31,308.5) and (242.5,304.69) .. (242.5,300) -- cycle ;
\draw  [color={rgb, 255:red, 0; green, 0; blue, 0 }  ,draw opacity=1 ][fill={rgb, 255:red, 208; green, 2; blue, 27 }  ,fill opacity=1 ] (92.5,200.5) .. controls (92.5,195.81) and (96.31,192) .. (101,192) .. controls (105.69,192) and (109.5,195.81) .. (109.5,200.5) .. controls (109.5,205.19) and (105.69,209) .. (101,209) .. controls (96.31,209) and (92.5,205.19) .. (92.5,200.5) -- cycle ;
\draw  [color={rgb, 255:red, 0; green, 0; blue, 0 }  ,draw opacity=1 ][fill={rgb, 255:red, 208; green, 2; blue, 27 }  ,fill opacity=1 ] (92.92,300.47) .. controls (92.92,295.77) and (96.73,291.97) .. (101.42,291.97) .. controls (106.12,291.97) and (109.92,295.77) .. (109.92,300.47) .. controls (109.92,305.16) and (106.12,308.97) .. (101.42,308.97) .. controls (96.73,308.97) and (92.92,305.16) .. (92.92,300.47) -- cycle ;
\draw  [color={rgb, 255:red, 0; green, 0; blue, 0 }  ,draw opacity=1 ][fill={rgb, 255:red, 208; green, 2; blue, 27 }  ,fill opacity=1 ] (92.5,100) .. controls (92.5,95.31) and (96.31,91.5) .. (101,91.5) .. controls (105.69,91.5) and (109.5,95.31) .. (109.5,100) .. controls (109.5,104.69) and (105.69,108.5) .. (101,108.5) .. controls (96.31,108.5) and (92.5,104.69) .. (92.5,100) -- cycle ;
\draw    (400,302) -- (400,101) ;
\draw  [color={rgb, 255:red, 0; green, 0; blue, 0 }  ,draw opacity=1 ][fill={rgb, 255:red, 208; green, 2; blue, 27 }  ,fill opacity=1 ] (390.91,102.03) .. controls (390.91,97.34) and (394.72,93.53) .. (399.41,93.53) .. controls (404.11,93.53) and (407.91,97.34) .. (407.91,102.03) .. controls (407.91,106.73) and (404.11,110.53) .. (399.41,110.53) .. controls (394.72,110.53) and (390.91,106.73) .. (390.91,102.03) -- cycle ;
\draw  [color={rgb, 255:red, 0; green, 0; blue, 0 }  ,draw opacity=1 ][fill={rgb, 255:red, 208; green, 2; blue, 27 }  ,fill opacity=1 ] (391.5,302) .. controls (391.5,297.31) and (395.31,293.5) .. (400,293.5) .. controls (404.69,293.5) and (408.5,297.31) .. (408.5,302) .. controls (408.5,306.69) and (404.69,310.5) .. (400,310.5) .. controls (395.31,310.5) and (391.5,306.69) .. (391.5,302) -- cycle ;
\draw  [color={rgb, 255:red, 0; green, 0; blue, 0 }  ,draw opacity=1 ][fill={rgb, 255:red, 208; green, 2; blue, 27 }  ,fill opacity=1 ] (391.5,201.5) .. controls (391.5,196.81) and (395.31,193) .. (400,193) .. controls (404.69,193) and (408.5,196.81) .. (408.5,201.5) .. controls (408.5,206.19) and (404.69,210) .. (400,210) .. controls (395.31,210) and (391.5,206.19) .. (391.5,201.5) -- cycle ;
\draw [color={rgb, 255:red, 0; green, 0; blue, 0 }  ,draw opacity=1 ][line width=2.25]    (119,285) .. controls (180,222) and (181.5,179) .. (210.5,162) .. controls (239.5,145) and (334.5,157.5) .. (374,124) ;

\draw (268,120) node [anchor=north west][inner sep=0.75pt]  [font=\large]  {$\partial \Omega $};

\end{tikzpicture}

    \label{fig:elementsBasesB}
    }\\ 

    \subfloat[Interface element.]{

\tikzset{every picture/.style={line width=0.75pt}} 

\begin{tikzpicture}[x=0.75pt,y=0.75pt,yscale=-0.6,xscale=0.6]

\draw  [fill={rgb, 255:red, 248; green, 231; blue, 28 }  ,fill opacity=1 ] (106.19,187.33) -- (114.53,196.33) -- (108.91,201.53) -- (114.11,207.15) -- (106.63,214.08) -- (101.42,208.47) -- (95.81,213.67) -- (87.47,204.67) -- (93.09,199.47) -- (87.89,193.85) -- (95.37,186.92) -- (100.58,192.53) -- cycle ;
\draw  [fill={rgb, 255:red, 248; green, 231; blue, 28 }  ,fill opacity=1 ] (106.19,86.83) -- (114.53,95.83) -- (108.91,101.03) -- (114.11,106.65) -- (106.63,113.58) -- (101.42,107.97) -- (95.81,113.17) -- (87.47,104.17) -- (93.09,98.97) -- (87.89,93.35) -- (95.37,86.42) -- (100.58,92.03) -- cycle ;
\draw  [fill={rgb, 255:red, 248; green, 231; blue, 28 }  ,fill opacity=1 ] (404.6,88.86) -- (412.94,97.86) -- (407.33,103.06) -- (412.53,108.68) -- (405.04,115.62) -- (399.84,110) -- (394.22,115.2) -- (385.88,106.2) -- (391.5,101) -- (386.3,95.38) -- (393.79,88.45) -- (398.99,94.06) -- cycle ;
\draw  [fill={rgb, 255:red, 248; green, 231; blue, 28 }  ,fill opacity=1 ] (106.62,287.3) -- (114.95,296.3) -- (109.34,301.5) -- (114.54,307.12) -- (107.05,314.05) -- (101.85,308.44) -- (96.23,313.64) -- (87.9,304.64) -- (93.51,299.44) -- (88.31,293.82) -- (95.8,286.88) -- (101,292.5) -- cycle ;
\draw  [fill={rgb, 255:red, 248; green, 231; blue, 28 }  ,fill opacity=1 ] (405.19,288.83) -- (413.53,297.83) -- (407.91,303.03) -- (413.11,308.65) -- (405.63,315.58) -- (400.42,309.97) -- (394.81,315.17) -- (386.47,306.17) -- (392.09,300.97) -- (386.89,295.35) -- (394.37,288.42) -- (399.58,294.03) -- cycle ;
\draw  [fill={rgb, 255:red, 248; green, 231; blue, 28 }  ,fill opacity=1 ] (405.19,188.33) -- (413.53,197.33) -- (407.91,202.53) -- (413.11,208.15) -- (405.63,215.08) -- (400.42,209.47) -- (394.81,214.67) -- (386.47,205.67) -- (392.09,200.47) -- (386.89,194.85) -- (394.37,187.92) -- (399.58,193.53) -- cycle ;
\draw  [fill={rgb, 255:red, 184; green, 233; blue, 134 }  ,fill opacity=1 ] (356.19,186.83) -- (364.53,195.83) -- (358.91,201.03) -- (364.11,206.65) -- (356.63,213.58) -- (351.42,207.97) -- (345.81,213.17) -- (337.47,204.17) -- (343.09,198.97) -- (337.89,193.35) -- (345.37,186.42) -- (350.58,192.03) -- cycle ;
\draw  [fill={rgb, 255:red, 184; green, 233; blue, 134 }  ,fill opacity=1 ] (156.19,186.83) -- (164.53,195.83) -- (158.91,201.03) -- (164.11,206.65) -- (156.63,213.58) -- (151.42,207.97) -- (145.81,213.17) -- (137.47,204.17) -- (143.09,198.97) -- (137.89,193.35) -- (145.37,186.42) -- (150.58,192.03) -- cycle ;
\draw  [fill={rgb, 255:red, 184; green, 233; blue, 134 }  ,fill opacity=1 ] (356.19,286.83) -- (364.53,295.83) -- (358.91,301.03) -- (364.11,306.65) -- (356.63,313.58) -- (351.42,307.97) -- (345.81,313.17) -- (337.47,304.17) -- (343.09,298.97) -- (337.89,293.35) -- (345.37,286.42) -- (350.58,292.03) -- cycle ;
\draw  [fill={rgb, 255:red, 184; green, 233; blue, 134 }  ,fill opacity=1 ] (156.19,86.83) -- (164.53,95.83) -- (158.91,101.03) -- (164.11,106.65) -- (156.63,113.58) -- (151.42,107.97) -- (145.81,113.17) -- (137.47,104.17) -- (143.09,98.97) -- (137.89,93.35) -- (145.37,86.42) -- (150.58,92.03) -- cycle ;
\draw  [fill={rgb, 255:red, 184; green, 233; blue, 134 }  ,fill opacity=1 ] (256.19,86.83) -- (264.53,95.83) -- (258.91,101.03) -- (264.11,106.65) -- (256.63,113.58) -- (251.42,107.97) -- (245.81,113.17) -- (237.47,104.17) -- (243.09,98.97) -- (237.89,93.35) -- (245.37,86.42) -- (250.58,92.03) -- cycle ;
\draw  [fill={rgb, 255:red, 184; green, 233; blue, 134 }  ,fill opacity=1 ] (156.19,286.83) -- (164.53,295.83) -- (158.91,301.03) -- (164.11,306.65) -- (156.63,313.58) -- (151.42,307.97) -- (145.81,313.17) -- (137.47,304.17) -- (143.09,298.97) -- (137.89,293.35) -- (145.37,286.42) -- (150.58,292.03) -- cycle ;
\draw  [fill={rgb, 255:red, 184; green, 233; blue, 134 }  ,fill opacity=1 ] (256.96,187.66) -- (263.76,195) -- (258.25,200.11) -- (263.36,205.62) -- (255.66,212.75) -- (250.55,207.24) -- (245.04,212.34) -- (238.24,205) -- (243.75,199.89) -- (238.64,194.38) -- (246.34,187.25) -- (251.45,192.76) -- cycle ;
\draw  [fill={rgb, 255:red, 184; green, 233; blue, 134 }  ,fill opacity=1 ] (256.19,286.83) -- (264.53,295.83) -- (258.91,301.03) -- (264.11,306.65) -- (256.63,313.58) -- (251.42,307.97) -- (245.81,313.17) -- (237.47,304.17) -- (243.09,298.97) -- (237.89,293.35) -- (245.37,286.42) -- (250.58,292.03) -- cycle ;
\draw  [fill={rgb, 255:red, 184; green, 233; blue, 134 }  ,fill opacity=1 ] (356.19,86.83) -- (364.53,95.83) -- (358.91,101.03) -- (364.11,106.65) -- (356.63,113.58) -- (351.42,107.97) -- (345.81,113.17) -- (337.47,104.17) -- (343.09,98.97) -- (337.89,93.35) -- (345.37,86.42) -- (350.58,92.03) -- cycle ;
\draw   (151,100) -- (351,100) -- (351,300) -- (151,300) -- cycle ;
\draw    (151,351) -- (350,351) ;
\draw    (101,301) -- (101,100) ;
\draw  [color={rgb, 255:red, 0; green, 0; blue, 0 }  ,draw opacity=1 ][fill={rgb, 255:red, 208; green, 2; blue, 27 }  ,fill opacity=1 ] (242,351) .. controls (242,346.31) and (245.81,342.5) .. (250.5,342.5) .. controls (255.19,342.5) and (259,346.31) .. (259,351) .. controls (259,355.69) and (255.19,359.5) .. (250.5,359.5) .. controls (245.81,359.5) and (242,355.69) .. (242,351) -- cycle ;
\draw  [color={rgb, 255:red, 0; green, 0; blue, 0 }  ,draw opacity=1 ][fill={rgb, 255:red, 208; green, 2; blue, 27 }  ,fill opacity=1 ] (341.5,351) .. controls (341.5,346.31) and (345.31,342.5) .. (350,342.5) .. controls (354.69,342.5) and (358.5,346.31) .. (358.5,351) .. controls (358.5,355.69) and (354.69,359.5) .. (350,359.5) .. controls (345.31,359.5) and (341.5,355.69) .. (341.5,351) -- cycle ;
\draw  [color={rgb, 255:red, 0; green, 0; blue, 0 }  ,draw opacity=1 ][fill={rgb, 255:red, 208; green, 2; blue, 27 }  ,fill opacity=1 ] (142.5,351) .. controls (142.5,346.31) and (146.31,342.5) .. (151,342.5) .. controls (155.69,342.5) and (159.5,346.31) .. (159.5,351) .. controls (159.5,355.69) and (155.69,359.5) .. (151,359.5) .. controls (146.31,359.5) and (142.5,355.69) .. (142.5,351) -- cycle ;
\draw  [color={rgb, 255:red, 0; green, 0; blue, 0 }  ,draw opacity=1 ][fill={rgb, 255:red, 74; green, 144; blue, 226 }  ,fill opacity=1 ] (242.5,100) .. controls (242.5,95.31) and (246.31,91.5) .. (251,91.5) .. controls (255.69,91.5) and (259.5,95.31) .. (259.5,100) .. controls (259.5,104.69) and (255.69,108.5) .. (251,108.5) .. controls (246.31,108.5) and (242.5,104.69) .. (242.5,100) -- cycle ;
\draw  [color={rgb, 255:red, 0; green, 0; blue, 0 }  ,draw opacity=1 ][fill={rgb, 255:red, 74; green, 144; blue, 226 }  ,fill opacity=1 ] (342.5,100) .. controls (342.5,95.31) and (346.31,91.5) .. (351,91.5) .. controls (355.69,91.5) and (359.5,95.31) .. (359.5,100) .. controls (359.5,104.69) and (355.69,108.5) .. (351,108.5) .. controls (346.31,108.5) and (342.5,104.69) .. (342.5,100) -- cycle ;
\draw  [color={rgb, 255:red, 0; green, 0; blue, 0 }  ,draw opacity=1 ][fill={rgb, 255:red, 74; green, 144; blue, 226 }  ,fill opacity=1 ] (142.5,200) .. controls (142.5,195.31) and (146.31,191.5) .. (151,191.5) .. controls (155.69,191.5) and (159.5,195.31) .. (159.5,200) .. controls (159.5,204.69) and (155.69,208.5) .. (151,208.5) .. controls (146.31,208.5) and (142.5,204.69) .. (142.5,200) -- cycle ;
\draw  [color={rgb, 255:red, 0; green, 0; blue, 0 }  ,draw opacity=1 ][fill={rgb, 255:red, 74; green, 144; blue, 226 }  ,fill opacity=1 ] (242.5,200) .. controls (242.5,195.31) and (246.31,191.5) .. (251,191.5) .. controls (255.69,191.5) and (259.5,195.31) .. (259.5,200) .. controls (259.5,204.69) and (255.69,208.5) .. (251,208.5) .. controls (246.31,208.5) and (242.5,204.69) .. (242.5,200) -- cycle ;
\draw  [color={rgb, 255:red, 0; green, 0; blue, 0 }  ,draw opacity=1 ][fill={rgb, 255:red, 74; green, 144; blue, 226 }  ,fill opacity=1 ] (342.5,300) .. controls (342.5,295.31) and (346.31,291.5) .. (351,291.5) .. controls (355.69,291.5) and (359.5,295.31) .. (359.5,300) .. controls (359.5,304.69) and (355.69,308.5) .. (351,308.5) .. controls (346.31,308.5) and (342.5,304.69) .. (342.5,300) -- cycle ;
\draw  [color={rgb, 255:red, 0; green, 0; blue, 0 }  ,draw opacity=1 ][fill={rgb, 255:red, 74; green, 144; blue, 226 }  ,fill opacity=1 ] (142.5,300) .. controls (142.5,295.31) and (146.31,291.5) .. (151,291.5) .. controls (155.69,291.5) and (159.5,295.31) .. (159.5,300) .. controls (159.5,304.69) and (155.69,308.5) .. (151,308.5) .. controls (146.31,308.5) and (142.5,304.69) .. (142.5,300) -- cycle ;
\draw  [color={rgb, 255:red, 0; green, 0; blue, 0 }  ,draw opacity=1 ][fill={rgb, 255:red, 74; green, 144; blue, 226 }  ,fill opacity=1 ] (142.5,100) .. controls (142.5,95.31) and (146.31,91.5) .. (151,91.5) .. controls (155.69,91.5) and (159.5,95.31) .. (159.5,100) .. controls (159.5,104.69) and (155.69,108.5) .. (151,108.5) .. controls (146.31,108.5) and (142.5,104.69) .. (142.5,100) -- cycle ;
\draw  [color={rgb, 255:red, 0; green, 0; blue, 0 }  ,draw opacity=1 ][fill={rgb, 255:red, 74; green, 144; blue, 226 }  ,fill opacity=1 ] (342.5,200) .. controls (342.5,195.31) and (346.31,191.5) .. (351,191.5) .. controls (355.69,191.5) and (359.5,195.31) .. (359.5,200) .. controls (359.5,204.69) and (355.69,208.5) .. (351,208.5) .. controls (346.31,208.5) and (342.5,204.69) .. (342.5,200) -- cycle ;
\draw  [color={rgb, 255:red, 0; green, 0; blue, 0 }  ,draw opacity=1 ][fill={rgb, 255:red, 74; green, 144; blue, 226 }  ,fill opacity=1 ] (242.5,300) .. controls (242.5,295.31) and (246.31,291.5) .. (251,291.5) .. controls (255.69,291.5) and (259.5,295.31) .. (259.5,300) .. controls (259.5,304.69) and (255.69,308.5) .. (251,308.5) .. controls (246.31,308.5) and (242.5,304.69) .. (242.5,300) -- cycle ;
\draw    (153,51) -- (352,51) ;
\draw  [color={rgb, 255:red, 0; green, 0; blue, 0 }  ,draw opacity=1 ][fill={rgb, 255:red, 208; green, 2; blue, 27 }  ,fill opacity=1 ] (92.5,200.5) .. controls (92.5,195.81) and (96.31,192) .. (101,192) .. controls (105.69,192) and (109.5,195.81) .. (109.5,200.5) .. controls (109.5,205.19) and (105.69,209) .. (101,209) .. controls (96.31,209) and (92.5,205.19) .. (92.5,200.5) -- cycle ;
\draw  [color={rgb, 255:red, 0; green, 0; blue, 0 }  ,draw opacity=1 ][fill={rgb, 255:red, 208; green, 2; blue, 27 }  ,fill opacity=1 ] (92.92,300.47) .. controls (92.92,295.77) and (96.73,291.97) .. (101.42,291.97) .. controls (106.12,291.97) and (109.92,295.77) .. (109.92,300.47) .. controls (109.92,305.16) and (106.12,308.97) .. (101.42,308.97) .. controls (96.73,308.97) and (92.92,305.16) .. (92.92,300.47) -- cycle ;
\draw  [color={rgb, 255:red, 0; green, 0; blue, 0 }  ,draw opacity=1 ][fill={rgb, 255:red, 208; green, 2; blue, 27 }  ,fill opacity=1 ] (92.5,100) .. controls (92.5,95.31) and (96.31,91.5) .. (101,91.5) .. controls (105.69,91.5) and (109.5,95.31) .. (109.5,100) .. controls (109.5,104.69) and (105.69,108.5) .. (101,108.5) .. controls (96.31,108.5) and (92.5,104.69) .. (92.5,100) -- cycle ;
\draw  [color={rgb, 255:red, 0; green, 0; blue, 0 }  ,draw opacity=1 ][fill={rgb, 255:red, 208; green, 2; blue, 27 }  ,fill opacity=1 ] (244,51) .. controls (244,46.31) and (247.81,42.5) .. (252.5,42.5) .. controls (257.19,42.5) and (261,46.31) .. (261,51) .. controls (261,55.69) and (257.19,59.5) .. (252.5,59.5) .. controls (247.81,59.5) and (244,55.69) .. (244,51) -- cycle ;
\draw  [color={rgb, 255:red, 0; green, 0; blue, 0 }  ,draw opacity=1 ][fill={rgb, 255:red, 208; green, 2; blue, 27 }  ,fill opacity=1 ] (144.5,51) .. controls (144.5,46.31) and (148.31,42.5) .. (153,42.5) .. controls (157.69,42.5) and (161.5,46.31) .. (161.5,51) .. controls (161.5,55.69) and (157.69,59.5) .. (153,59.5) .. controls (148.31,59.5) and (144.5,55.69) .. (144.5,51) -- cycle ;
\draw  [color={rgb, 255:red, 0; green, 0; blue, 0 }  ,draw opacity=1 ][fill={rgb, 255:red, 208; green, 2; blue, 27 }  ,fill opacity=1 ] (343.5,51) .. controls (343.5,46.31) and (347.31,42.5) .. (352,42.5) .. controls (356.69,42.5) and (360.5,46.31) .. (360.5,51) .. controls (360.5,55.69) and (356.69,59.5) .. (352,59.5) .. controls (347.31,59.5) and (343.5,55.69) .. (343.5,51) -- cycle ;
\draw    (400,302) -- (400,101) ;
\draw  [color={rgb, 255:red, 0; green, 0; blue, 0 }  ,draw opacity=1 ][fill={rgb, 255:red, 208; green, 2; blue, 27 }  ,fill opacity=1 ] (390.91,102.03) .. controls (390.91,97.34) and (394.72,93.53) .. (399.41,93.53) .. controls (404.11,93.53) and (407.91,97.34) .. (407.91,102.03) .. controls (407.91,106.73) and (404.11,110.53) .. (399.41,110.53) .. controls (394.72,110.53) and (390.91,106.73) .. (390.91,102.03) -- cycle ;
\draw  [color={rgb, 255:red, 0; green, 0; blue, 0 }  ,draw opacity=1 ][fill={rgb, 255:red, 208; green, 2; blue, 27 }  ,fill opacity=1 ] (391.5,302) .. controls (391.5,297.31) and (395.31,293.5) .. (400,293.5) .. controls (404.69,293.5) and (408.5,297.31) .. (408.5,302) .. controls (408.5,306.69) and (404.69,310.5) .. (400,310.5) .. controls (395.31,310.5) and (391.5,306.69) .. (391.5,302) -- cycle ;
\draw  [color={rgb, 255:red, 0; green, 0; blue, 0 }  ,draw opacity=1 ][fill={rgb, 255:red, 208; green, 2; blue, 27 }  ,fill opacity=1 ] (391.5,201.5) .. controls (391.5,196.81) and (395.31,193) .. (400,193) .. controls (404.69,193) and (408.5,196.81) .. (408.5,201.5) .. controls (408.5,206.19) and (404.69,210) .. (400,210) .. controls (395.31,210) and (391.5,206.19) .. (391.5,201.5) -- cycle ;
\draw [color={rgb, 255:red, 0; green, 0; blue, 0 }  ,draw opacity=1 ][line width=2.25]    (119,285) .. controls (180,222) and (181.5,179) .. (210.5,162) .. controls (239.5,145) and (334.5,157.5) .. (374,124) ;

\draw (268,120) node [anchor=north west][inner sep=0.75pt]  [font=\large]  {$\Upsilon $};

\end{tikzpicture}
\label{fig:elementsBasesC}
    }
    \caption{Schematic representation of the unknowns for the different element types.  (a) Standard HDG element. (b) Immersed boundary element. (c) Interface element. Blue bullets represent local unknowns, red bullets hybrid unknowns, green crosses the doubled local unknowns employed in interface elements, and yellow crosses the doubled hybrid unknowns employed in interface cut faces.}
    \label{fig:elementsBases}
\end{figure}
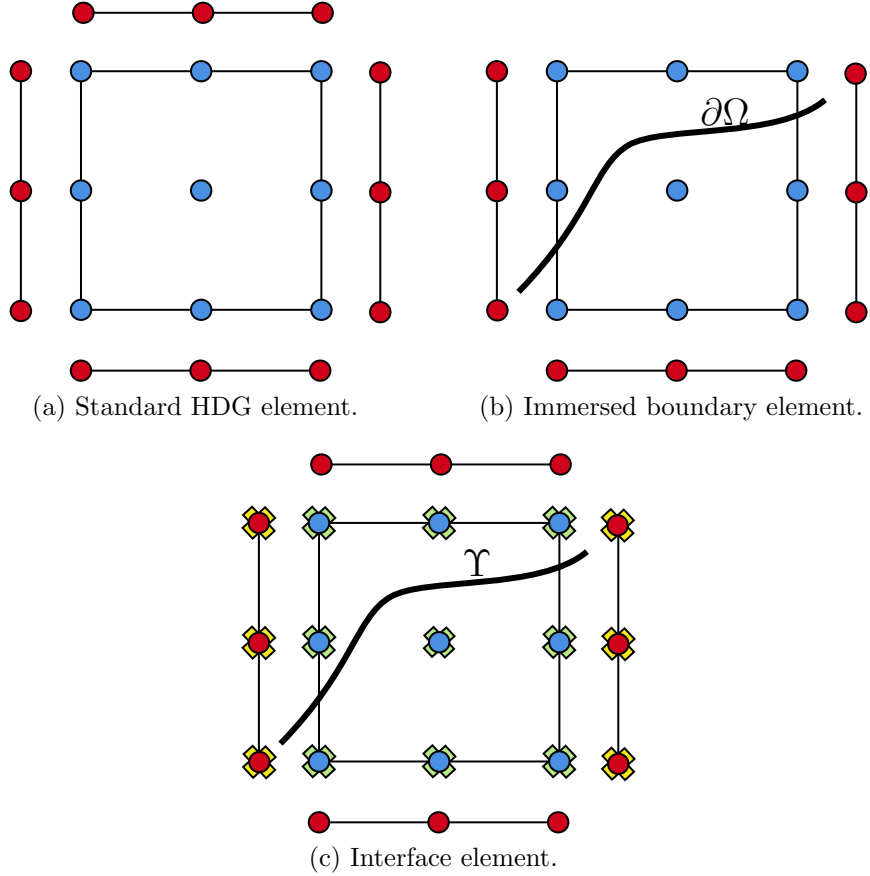

As usual in HDG, equal-order polynomial approximations are used for velocity, pressure and the mixed variable. Of course, these approximations only exist for those elements occupied (totally or partially) by fluid, that is, for those $\Omega_e$ such that $\Omega_e\cap\Omega^i\neq\emptyset$, for at least one fluid $i$. In addition, for the interface elements these approximations are duplicated to determine velocity, pressure and the mixed variable in each fluid. Similarly, the hybrid variable, the trace of the velocity, $\bhu$, whose polynomial order is chosen accordingly to the highest one in the neighboring elements, is defined only on the faces (totally or partially) in the fluid, that is, when $\Gamma_{\!\! f}\cap\Omega^i\neq\emptyset$, for at least one fluid $i$. Analogously, the unknowns are duplicated for those faces intersected by the interface $\Su$. Figure\ \ref{fig:elementsBases} displays a schematic representation of these approximations for the three element types under analysis. 
%

Henceforth, consider rectangular computational domains $\meshDomain$, partitioned in rectangular-shaped elements $\Omega_e$. 
In this specific context, for the element-based unknowns, Lagrange basis functions are employed as common in the case of fitted HDG methods. 
For the face-based unknowns, two possible bases are compared: the Lagrange basis, $\Pk_{\text{La}}$, and the Legendre basis, $\Pk_{\text{Le}}$. 
In Section\ \ref{sc:conditioning}, a comparative analysis of both bases shows the superior performance of Legendre face-based unknowns in terms of mitigating the conditioning of the global problem.

%
\subsection{HDG local problems} 
\label{sc:problemDiscLoc}
Consistently with the three element types identified (standard, immersed and interface), three types of local problems must be solved. In the first case, the strong form of the local problem for \emph{standard HDG elements} is already detailed in equation~\eqref{eq:PSLocal}, supplemented, as noted earlier, with the hybrid velocity $\bhu^i$ along $\Gamma \cap \Omega^i$ (see \cite{MG-GSH:20} for more details). Recall that the numerical trace of the diffusive flux must be defined. Following \cite{Jay-CGL:09,MG-GSH:20} it is prescribed as 
%
\begin{equation}\label{eq:difflux}
     (\reallyhat{\sqrt{\mu^i} \bL_e^i {+} p_e^i \Insd})\,\bn_e :=  ({\sqrt{\mu^i} \bL_e^i {+} p_e^i \Insd})\,\bn_e + \begin{cases} 
     \tau (\bu_e^i - \dirData), &\text{on $\hGa{D}\cap\partial\Omega_e^i$,}\\
     \tau (\bu_e^i - \bhu^i),    &\text{on $\Gamma\cap\partial\Omega_e^i$,}
     \end{cases} 
\end{equation}
where $\tau$ is the stabilization parameter. Consistently with\ \cite{GiaSeH:21}, this parameter is defined as $ \tau:= \C \mu / \ell$ where $\mu := \max(\mu^1,\mu^2)$, $\ell$ is a characteristic length of the domain, and $\C$ is a scaling factor related to the element type, which is henceforth set equal to $3$. See \cite{MG-GSH:20} for more details.  

\begin{remark}[Stabilization for high viscosity constrast]
It is known that the definition of jump and average operators with weighting coefficients depending on the viscosity can improve the accuracy and the robustness of the resulting numerical schemes in the presence of heterogeneous materials\ \cite{Ern-ESZ-09,Dolbow-AHD-12,Giacomini-18}, especially in the case of high contrast.
The numerical tests presented in Section\ \ref{sc:numSimu} feature moderate viscosity contrasts for which the stabilization coefficient introduced above provides stable and accurate results. 
The suitability of the proposed stabilization for higher contrasts is beyond the scope of the present work and should be appropriately assessed with realistic benchmarks. It is worth mentioning that some recent studies\ \cite{BuDeE:21,Sevilla-SD-23} showed promising results in this direction, highlighting the robustness of hybrid discretization methods for the treatment of high contrast viscous flow problems, without introducing weighted jump and average operators.
\end{remark}

\subsubsection{Immersed boundary elements}
\label{sec:immElem}
In what follows, $(\cdot,\cdot)_{S}$ and $\langle\cdot,\cdot\rangle_S$ denote the $\eltwo$-product in a generic subset ${S}\subset\Omega$ of dimension $\nsd$ and $\nsd{-}1$, respectively.
For immersed boundary elements, a Nitsche's penalty term\ \cite{Nitsche:71,Arnold-ABCM:02} is used to weakly impose the Dirichlet boundary conditions when $\hGa{D}$ cuts $\Omega_e$, namely, the weak formulation of~\eqref{eq:PSLocal} is as follows: for all $e\inSet{1}{\numel}$ such that $\Omega_e\cap\Su=\emptyset$ and $\Omega_e\cap\partial\Omega\neq\emptyset$, given $\dirData$ on $\Ga{D}{i}$,   $\bt$ on $\Ga{N}{i}$, and $\bhu^i$ on $\Gamma\cap\partial\Omega_e^i$, find $(\bL_e^i, \bu_e^i, p_e^i) \in [\hDiv{\Omega_e^i}]^{\nsd\times\nsd} \times [\sobo(\Omega_e^i)]^{\nsd} \times \sobo(\Omega_e^i)$ that satisfy 
\begin{equation}\label{eq:weakLocalOneFluid}
\left\{\begin{aligned}
  {-}&\bigl( \bG, \bL_e^i \bigr)_{\Omega_e^i}  
  + \bigl( \Div(\sqrt{\mu^i}\bG), \bu_e^i \bigr)_{\Omega_e^i}
  - \langle \bG\,\bn_e, \sqrt{\mu^i}\,\bu_e^i \rangle_{\Ga{N}{i}\cap \overline{\Omega_e}} 
   \\[-0.5ex] &\hspace{150pt}
  = \langle \bG\,\bn_e, \sqrt{\mu^i}\,\dirData \rangle_{\Ga{D}{i}\cap \overline{\Omega_e}} 
  \\[-0.5ex] &\hspace{170pt}
  + \langle \bG\,\bn_e, \sqrt{\mu^i}\,\bhu^i \rangle_{\partial \Omega_e^i \setminus(\Ga{D}{i}\cup\Ga{N}{i})} ,
\\[1ex]
  &\bigl( \bw, \grad{\cdot}( \sqrt{\mu^i} \bL_e^i) \bigr)_{\Omega_e^i}  
  {+} \bigl( \bw, \grad p_e^i \bigr)_{\Omega_e^i} 
  \\[-0.5ex] &\hspace{35pt}
  {-} \bigl\langle \bw, (\sqrt{\mu^i} \bL_e^i {+} p_e^i \Insd)\,\bn_e \bigr\rangle_{\Ga{N}{i}\cap \overline{\Omega_e}} 
  {+} \bigl\langle \bw, \tau\bu_e^i \bigr\rangle_{\partial\Omega_e^i \setminus\Ga{N}{i}}
  {+} \Nitsche h_e^{-1}\langle\bw,\bu_e^i\rangle_{\Ga{D}{i} \cap \Omega_e}
  \\[-0.5ex] &\hspace{150pt}
  = \bigl( \bw, \bm{s} \bigr)_{\Omega_e^i} 
  {+} \bigl\langle \bw, \bt \bigr\rangle_{\Ga{N}{i}\cap \overline{\Omega_e}}
  {+} \bigl\langle \bw, \tau\,\dirData \bigr\rangle_{\partial\Omega_e^i\cap\Ga{D}{i}}
  \\[-0.5ex] &\hspace{170pt}
  {+} \Nitsche h_e^{-1}\langle\bw,\dirData\rangle_{\Ga{D}{i} \cap \Omega_e}
  {+} \bigl\langle \bw, \tau\,\bhu    \bigr\rangle_{\partial\Omega_e^i\setminus(\Ga{D}{i}\cup\Ga{N}{i})} ,
\\[1ex]
    &\bigl( \grad q, \bu_e^i \bigr)_{\Omega_e^i}  
 -  \langle q, \bu_e^i  \cdot \bn_e \rangle_{\Ga{N}{i}\cap \overline{\Omega_e}} 
 = \langle q, \dirData \cdot \bn_e \rangle_{\Ga{D}{i}\cap \overline{\Omega_e}}
 + \langle q, \bhu^i \cdot \bn_e \rangle_{\partial \Omega_e^i \setminus(\Ga{D}{i}\cup\Ga{N}{i}) },\\
\end{aligned}\right.  
 \end{equation}
for all $(\bG, \bw, q) \in [\hDiv{\Omega_e^i}]^{\nsd\times\nsd} \times [\sobo(\Omega_e^i)]^{\nsd} \times \sobo(\Omega_e^i)$, where $i$ is either $1$ or $2$, and where the definition in~\eqref{eq:difflux} of the numerical trace of the diffusive flux has been used. 
Following the literature on interior penalty and fictitious domain methods for the weak imposition of essential boundary conditions\ \cite{Arnold-82,Stenberg-95,Larson-JL-13},  Nitsche's penalty term is set proportional to the inverse of the element size, independently of the employed polynomial degree of approximation. Moreover, 
recall that Nitsche's consistent penalty constant, $\Nitsche$, must be large enough. In all the problems presented in Section\ \ref{sc:numSimu}, $\Nitsche {=} 10$ proved to be sufficient, and, as expected, results are insensitive to the choice of $\Nitsche$.

Note that for those elements with pure Dirichlet boundary conditions, that is, when $\overline{\Omega_e} \cap\Ga{N}{i} = \emptyset$, the problem is closed by adding the pressure condition, 
\begin{equation*}
   \frac{1}{|\Omega_e^i|} (1,p_e^i)_{\Omega_e^i} = \rho_e,
\end{equation*}
with $\rho_e$ an independent global variable representing the mean value of the pressure in $\Omega_e^i$.

%
\subsubsection{Interface elements}
\label{sec:intElem}
For the interface elements, occupied by both immiscible fluids and such that $\Omega_e\cap\partial\Omega=\emptyset$, the interface conditions couple the local variables of $\Omega_e^1$ and $\Omega_e^2$ (velocity, pressure and the mixed variable). Actually, at those interfaces, the average and jump operators determine the identity
\begin{equation}\label{eq:JumpMeanIdentity}
  \sum_{i=1}^2 \bigl\langle \bw^i , \bG^i\bn^i \bigr\rangle_{\Su} 
   = \bigl\langle\jump{\bw\otimes\bn} , \mean{\bG} \bigr\rangle_{\Su} + \bigl\langle \mean{\bw}  , \jump{\bG\bn} \bigr\rangle_{\Su} ,
\end{equation}
see \cite{Arnold-ABCM:02}. Recall that the \emph{jump} operator was already defined in~\eqref{eq:jump} whereas  the \emph{mean}, $\mean{\odot}$, is given by
\begin{equation}
\mean{\odot} := (\odot|_{\Omega^1} + \odot|_{\Omega^2}) /2 ,
\label{eq:HDGmeanOp}
\end{equation}
see for instance \cite{HanHan:04}. Both follow the definition proposed in \cite{AdM-MFH:08}. 

The interface elements present weak local problems as follows: for all $e\inSet{1}{\numel}$ such that $\Omega_e\cap\Su\neq\emptyset$ and $\Omega_e\cap\partial\Omega=\emptyset$, given $\dirData$ on $\hGa{D}{}$,   $\bt$ on $\hGa{N}{}$, $\bhu^1$ on $\Gamma\cap\partial\Omega^1_e$, and $\bhu^2$ on $\Gamma\cap\partial\Omega^2_e$ (i.e., $\bhu$ on $\Gamma\cap \partial\Omega_e$), find $(\bL_e, \bu_e, p_e) \in [\hDiv{\Omega_e\cap\hOm}]^{\nsd\times\nsd} \times [\sobo(\Omega_e\cap\hOm)]^{\nsd} \times \sobo(\Omega_e\cap\hOm)$ that satisfy
\begin{equation}\label{eq:weakLocalTwoFluid}
\left\{\begin{aligned}
  {-}&\bigl( \bG, \bL_e \bigr)_{\Omega_e}  
   \\[-0.5ex] &\hspace{10pt}
  + \bigl( \Div(\sqrt{\mu}\bG), \bu_e \bigr)_{\Omega_e}
  - \langle \bG\,\bn_e, \sqrt{\mu}\,\bu_e \rangle_{\partial\Omega_e\cap\hGa{N}} 
   - \langle\jump{\sqrt{\mu} \bG\bn_e},  \mean{\bu_e}\rangle_{\Su\cap\Omega_e}  
  \\[-0.5ex] &\hspace{105pt}
  = \langle \bG\,\bn_e, \sqrt{\mu}\,\dirData \rangle_{\partial\Omega_e\cap\hGa{D}} 
  + \langle \bG\,\bn_e, \sqrt{\mu}\,\bhu \rangle_{\partial \Omega_e \setminus (\hGa{D}\cup\hGa{N}) } ,
\\[1ex]
  &\bigl\langle \bw, \tau\,\bu_e \bigr\rangle_{\partial\Omega_e\setminus\hGa{N}}
   \\[-0.5ex] &\hspace{10pt}
  {+} \bigl( \bw, \grad{\cdot}( \sqrt{\mu} \bL_e) \bigr)_{\Omega_e}  
  {-} \bigl\langle \bw, \sqrt{\mu} \bL_e \,\bn_e \bigr\rangle_{\partial\Omega_e\cap\hGa{N}}
  {-} \bigl\langle\mean{\bw},  \jump{\sqrt{\mu} \bL_e \,\bn_e} \bigr\rangle_{\Su\cap\Omega_e}  
   \\[-0.5ex] &\hspace{20pt}
  {+} \bigl( \bw, \grad p_e \bigr)_{\Omega_e}
  {-} \bigl\langle \bw, p_e \,\bn_e \bigr\rangle_{\partial\Omega_e\cap\hGa{N}}
  {-} \bigl\langle\mean{\bw},  \jump{ p_e \,\bn_e} \bigr\rangle_{\Su\cap\Omega_e}
  \\[-0.5ex] &\hspace{105pt}
  = \bigl( \bw, \bm{s} \bigr)_{\Omega_e} 
  {+} \bigl\langle \bw, \bt \bigr\rangle_{\partial\Omega_e\cap\hGa{N}} 
  {+} \bigl\langle \bw, \tau\,\dirData \bigr\rangle_{\partial\Omega_e\cap\hGa{D}}
  \\[-0.5ex] &\hspace{135pt}
  {+} \bigl\langle \bw, \tau\,\bhu    \bigr\rangle_{\partial\Omega_e\setminus(\hGa{D}\cup\hGa{N})} 
  {+} \bigl\langle\mean{\bw}, \gamma(\Div\bn^1_e)\bn^1_e \bigr\rangle_{\Su\cap\Omega_e} ,
\\[1ex]
    &\bigl( \grad q, \bu_e \bigr)_{\Omega_e}  
 -  \langle q, \bu_e  \cdot \bn_e \rangle_{\partial \Omega_e \cap \hGa{N} } 
   - \langle \jump{q\,\bn_e} , \mean{\bu_e}\rangle_{\Su\cap\Omega_e} 
  \\[-0.5ex] &\hspace{105pt}
 = \langle q, \dirData \cdot \bn_e \rangle_{\partial \Omega_e \cap \hGa{D}}
 + \langle q, \bhu \cdot \bn_e \rangle_{\partial \Omega_e \setminus (\hGa{D}\cup\hGa{N}) } ,
\end{aligned}\right.  
 \end{equation}
for all $(\bG, \bw, q) \in [\hDiv{\Omega_e\cap\hOm}]^{\nsd\times\nsd} \times [\sobo(\Omega_e\cap\hOm)]^{\nsd} \times \sobo(\Omega_e\cap\hOm)$, where the numerical trace of the diffusive flux defined in~\eqref{eq:difflux} has been used and the interface conditions~\eqref{eq:PSInter} have been imposed such that symmetry is preserved. Actually, it is worth noticing that~\eqref{eq:weakLocalTwoFluid} preserves the symmetry of the Stokes problem.
Moreover, as shown in the numerical examples, since in each element $\abs{\partial\Omega_e{\setminus}\hGa{N}}{>}\abs{\Su}$, there is no need for adding any extra stabilization term along $\Su$.
Note also that the following three particularizations of~\eqref{eq:JumpMeanIdentity} have been exploited:
%
\begin{align*}
   \bigSumE{i}{1}{2} \langle \bG\,\bn^i, \sqrt{\mu^i}\,\bu^i \rangle_{\Su}
  =& \langle \jump{\bu\otimes\bn} , \mean{\sqrt{\mu} \bG}\rangle_{\Su} 
   + \langle \mean{\bu},\jump{\sqrt{\mu} \bG\,\bn}\rangle_{\Su} ,
  \\[-1ex]
    \bigSumE{i}{1}{2} \langle \bw^i, (\sqrt{\mu^i} \bL^i {+} p^i \Insd)\,\bn^i  \rangle_{\Su} 
  =&  \langle \jump{\bw\otimes\bn} , \mean{\sqrt{\mu} \bL {+} p \Insd}\rangle_{\Su} 
  \\[-2ex]
   & {+} \langle \mean{\bw},\jump{(\sqrt{\mu} \bL {+} p \Insd)\bn}\rangle_{\Su} ,
  \\[-1ex]
    \bigSumE{i}{1}{2}  \langle q^i, \bu^i \cdot \bn^i \rangle_{\Su} 
  =&  \langle \jump{q\,\bn} , \mean{\bu} \rangle_{\Su} 
   + \langle \mean{q}, \jump{\bu\cdot\bn} \rangle_{\Su} .
\end{align*}
%
Finally, as before, when $\partial\Omega_e \cap \hGa{N} = \emptyset$, the problem is closed by adding the pressure condition, 
\begin{equation*}
    \frac{1}{|\Omega_e|} \bigSumE{i}{1}{2}  (p^i_e,1)_{\Omega_e \cap \Omega^i} = \rho_{e}.
\end{equation*}

\begin{remark}[HDG with Cauchy stress formulation]
In many problems of engineering interest, the employment of the Cauchy stress formulation of the Stokes equations is preferred. Besides the modeling advantages related to Neumann conditions representing physical tensions (see, e.g.,\ \cite{Donea-Huerta}), the choice of the strain rate tensor as point-wise symmetric HDG mixed variable using Voigt notation can also lead to a reduction of the size of the local problems, by considering only the non-redundant components of $\bL$\ \cite{RS-SGKH-18,MG-GKSH:18,MG-GSH:20,LaSpina-SKGWH-20,JVP-VGSH-21}.
\end{remark}

%
\subsection{HDG global problem}
\label{sc:GlobProblem}
The previous section details element-by-element problems that allow to express the local variables $(\bL_e^i, \bu_e^i, p_e^i)$ in terms of the global unknowns $(\bhu^i,\rho_e)$. The global problem represented by the \emph{transmission conditions}~\eqref{eq:PSTrans} solves for the global unknowns $(\bhu^i,\rho_e)$.  By definition, $\bhu^i$ is uniquely defined on each face of $\Gamma \cap \Omega^i$ thus condition $\jump{\bu^i \otimes \bn^i} = \bm{0}$ is automatically satisfied because of~\eqref{eq:DirichletHybrid}.

The integral formulation resulting from the transmission condition is as follows:
 given $\dirData$ on $\hGa{D}{}$,   $\bt$ on $\hGa{N}{}$, and the expressions of the local variables $(\bL_e^i, \bu_e^i, p_e^i)$ in terms of the global unknowns $(\bhu^i,\rho_e)$ from the local problems, find $(\bhu^i,\rho_e) \in [\sobo[\frac{1}{2}](\Gamma\cap \partial\Omega_e^i)]^{\nsd} \times \RR$ that satisfy
 %
\begin{subequations}\label{eq:GlobalIntegral}
\begin{multline}\label{eq:GlobalIntegralTrans}
   \bigSumE{e}{1}{\numel} \bigSumE{i}{1}{2}
  \bigl\{ \langle \bhw,(\sqrt{\mu^i}\bL_e^i + p_{e}^i \Insd)\,\bn_e \rangle_{\partial \Omega_e^i \setminus (\hGa{D}\cup\hGa{N}) } \bigl. \\[-1ex]
  \bigl.+ \langle \bhw,\tau \bu_e^i \rangle_{\partial \Omega_e^i \setminus (\hGa{D}\cup\hGa{N})}  - \langle\bhw,\tau \bhu^i\rangle_{\partial \Omega_e^i \setminus (\hGa{D}\cup\hGa{N})}\bigl\} = 0,
\end{multline}
for all $\bhw\in [\mathcal{L}_2 (\Gamma\cap \partial\Omega_e^i)]^{\nsd}$. In addition to the transmission conditions, for all elements not affected by the Neumann boundary condition, that is, for $e\inSet{1}{\numel}$ such that $\overline{\Omega_e} \cap \hGa{N} = \emptyset$, the following compatibility condition is considered
\begin{equation}\label{eq:GlobalIntegralCompa}
  \bigSumE{i}{1}{2} \bigl\{\langle \bhu^i\cdot \bn_e     , 1 \rangle_{\partial \Omega_e^i \setminus \hGa{D}} 
                                   +  \langle \dirData \cdot \bn_e , 1 \rangle_{\partial \Omega_e^i \cap \hGa{D} }  \bigl\}= 0.
\end{equation}
\end{subequations}
The above equation is derived from the divergence-free condition, where the continuity of velocity along the interface $\Su$, see~\eqref{eq:PSInter}, has been exploited.
Also, recall that the additional pressure condition, presented in Remark\ \ref{rm:pressureC0}, is required when $\hGa{N} = \emptyset$.

%
The HDG global problem~\eqref{eq:GlobalIntegral} yields a symmetric system at the discrete level, independently of the presence of immersed boundary elements or interface elements.  
Indeed, equation~\eqref{eq:GlobalIntegral} only involves the transmission conditions on the internal mesh skeleton, as in standard HDG methods, and the compatibility condition in the elements not affected by the Neumann boundary condition.  Recall that no additional unknowns are introduced, neither on $\partial\Omega$ nor on $\Su$.
In particular, note that in the case of a two-fluid system with an interface,  condition~\eqref{eq:PSInter} is enforced in the HDG local problem, coupling the local variables $(\bL_e^1, \bu_e^1, p_e^1)$ of $\Omega_e^1$ and $(\bL_e^2, \bu_e^2, p_e^2)$ of $\Omega_e^2$ via the terms along $\Su$. 
This allows to preserve the symmetry of the local system, as detailed in equation~\eqref{eq:weakLocalTwoFluid}, and the resulting global problem~\eqref{eq:GlobalIntegral} is thus not directly affected by the presence of the interface $\Su$.

\subsection{HDG postprocess and degree adaptivity}
\label{sec:degAdapt}
Following \cite{Nguyen-NPC:10,Cockburn-CGNPS:11,GG-GFH:14,MG-GS:19,MG-GSH:20}, for $i=1,2$, the discrete functional space composed of polynomial functions of degree at most $k{+}1$ is defined as
\begin{equation}
    \mathcal{V}_{\star}^h(\hOm) :=\bigl\{ v\in \mathcal{L}_2(\hOm)\mid \; v\vert_{\Omega_e^i}\in \PkPone (\Omega_e^i), \; \forall e \inSet{1}{\numel} , \; \text{and }\forall i=1,2 \bigl\}.
\end{equation}
Then, the accuracy of the velocity field approximation can be enhanced, achieving convergence of order $k{+}2$, by solving an additional, computationally efficient element-by-element problem. For
$e\inSet{1}{\numel}$, find the superconvergent velocity field $\buS\in [\mathcal{V}_{\star}^h(\hOm)]^{\nsd}$ such that $\buSe^i = \buS\vert_{\Omega_e^i}$ for $i=1,2$ and
\begin{equation}\label{eq:postProc}
  \left\{\begin{alignedat}{3}
       -\Div(\sqrt{\mu^i}\grad\buSe^i) &= \Div\bL_e^i &\quad& \text{in }\Omega_e^i,\\
        \sqrt{\mu^i}\grad\buSe^i\, \bn_e &= -\bL_e^i \,\bn_e &\quad& \text{on }\partial\Omega_e^i \cup (\Su\cap \Omega_e),
  \end{alignedat}\right. 
 \end{equation}
and such that the mean superconvergent velocity satisfies
\begin{equation}
    (\buSe^i,1)_{\Omega_e^i} = (\bu_{e}^i,1)_{\Omega_e^i}.
\end{equation}
Hence, the weak problems for each fluid are: for $e\inSet{1}{\numel}$ and $i=1,2$,  solve 
\begin{equation}
\begin{aligned}
     (\grad \bwS^i,\sqrt{\mu^i}\grad\buSe^i)_{\Omega_e^i} &= -(\grad \bwS^i,\bL_e^i)_{\Omega_e^i},\\
       (\buSe^i,1)_{\Omega_e^i} &= (\bu_{e}^i,1)_{\Omega_e^i},
\end{aligned}
\end{equation}
with $\bwS^i$ being a test function of polynomial degree at most $k{+}1$ in $\Omega_e^i$.

Consistently with \cite{SevHu:18}, the following element-wise measure of the error is employed
\begin{equation}
    E_e := \max_{i=1,2}\Bigl\{ \bigl(\buSe^i-\bu_e^i,\buSe^i-\bu_e^i \bigr)_{\Omega_e^i}^{1/2} / 
                                              \bigl(\buSe^i             ,\buSe^i             \bigr)_{\Omega_e^i}^{1/2} \Bigr\}, \qquad \text{for } e \inSet{1}{\numel}.
\end{equation}
Obviously, for values of $\bigl(\buSe^i ,\buSe^i \bigr)_{\Omega_e^i}^{1/2} $ below a prescribed tolerance the absolute error is employed.

The strategy depicted in \cite{SevHu:18} allows to attain the desired accuracy $\varepsilon$ on the error $E_e$ by locally adapting the degree of the HDG approximation.
Problem~\eqref{eq:probstat} is thus solved iteratively increasing or decreasing, for each element $\Omega_e$, the local degree $k_e$ of approximation as
\begin{equation}
    \delta k_e = \Bigl\lceil{ \frac{\log(\varepsilon/E_e)}{\log(h_e)} }\Bigl\rceil,
\end{equation}
with $\lceil \cdot\rceil$ being the ceiling function and $h_e$ the local mesh size of element $\Omega_e$.
The procedure ends when $\delta k_e = 0$, for $e \inSet{1}{\numel}$.

\section{Computational setting for non-conforming meshes}
\label{sec:quadDef}
In this section,  technical details related to the treatment of the unfitted interface, the numerical quadrature in the presence of NURBS and the element extension strategy to ensure robustness with respect to badly cut cells are presented.

\subsection{Identification of the two fluid subdomains}
\label{sec:identifyFluid}
\begin{figure}[!h]
  \centering
  \tikzset{every picture/.style={line width=0.75pt}} 
  \begin{tikzpicture}[x=0.75pt,y=0.75pt,yscale=-1,xscale=1]
\draw [line width=0.75]    (167,102) -- (145.12,80.12) ;
\draw [shift={(143,78)}, rotate = 45] [fill={rgb, 255:red, 0; green, 0; blue, 0 }  ][line width=0.08]  [draw opacity=0] (8.93,-4.29) -- (0,0) -- (8.93,4.29) -- cycle    ;
\draw  [line width=0.75]  (120,32) -- (323,32) -- (323,235) -- (120,235) -- cycle ;
\draw [color={rgb, 255:red, 74; green, 144; blue, 226 }  ,draw opacity=1 ][line width=2.25]    (261,200) .. controls (133,197) and (131,120) .. (180,92) .. controls (221,67) and (341,89) .. (286,116) .. controls (231,143) and (264,187) .. (289,162) .. controls (314,137) and (319,200) .. (261,200) -- cycle ;
\draw  [line width=1.5]  (237.03,75.15) -- (249.76,82.1) -- (236.55,88.07) ;
\draw  [draw opacity=0] (120,32) -- (323,32) -- (323,235) -- (120,235) -- cycle ; \draw   (149,32) -- (149,235)(178,32) -- (178,235)(207,32) -- (207,235)(236,32) -- (236,235)(265,32) -- (265,235)(294,32) -- (294,235) ; \draw   (120,61) -- (323,61)(120,90) -- (323,90)(120,119) -- (323,119)(120,148) -- (323,148)(120,177) -- (323,177)(120,206) -- (323,206) ; \draw   (120,32) -- (323,32) -- (323,235) -- (120,235) -- cycle ;
\draw [line width=0.75]    (169,172) -- (187.08,150.3) ;
\draw [shift={(189,148)}, rotate = 129.81] [fill={rgb, 255:red, 0; green, 0; blue, 0 }  ][line width=0.08]  [draw opacity=0] (8.93,-4.29) -- (0,0) -- (8.93,4.29) -- cycle    ;
\draw  [fill={rgb, 255:red, 208; green, 2; blue, 27 }  ,fill opacity=1 ] (293,158) .. controls (293,154.13) and (296.13,151) .. (300,151) .. controls (303.87,151) and (307,154.13) .. (307,158) .. controls (307,161.87) and (303.87,165) .. (300,165) .. controls (296.13,165) and (293,161.87) .. (293,158) -- cycle ;
\draw  [fill={rgb, 255:red, 208; green, 2; blue, 27 }  ,fill opacity=1 ] (254,200) .. controls (254,196.13) and (257.13,193) .. (261,193) .. controls (264.87,193) and (268,196.13) .. (268,200) .. controls (268,203.87) and (264.87,207) .. (261,207) .. controls (257.13,207) and (254,203.87) .. (254,200) -- cycle ;
\draw  [fill={rgb, 255:red, 208; green, 2; blue, 27 }  ,fill opacity=1 ] (266,125) .. controls (266,121.13) and (269.13,118) .. (273,118) .. controls (276.87,118) and (280,121.13) .. (280,125) .. controls (280,128.87) and (276.87,132) .. (273,132) .. controls (269.13,132) and (266,128.87) .. (266,125) -- cycle ;
\draw  [fill={rgb, 255:red, 208; green, 2; blue, 27 }  ,fill opacity=1 ] (258,83) .. controls (258,79.13) and (261.13,76) .. (265,76) .. controls (268.87,76) and (272,79.13) .. (272,83) .. controls (272,86.87) and (268.87,90) .. (265,90) .. controls (261.13,90) and (258,86.87) .. (258,83) -- cycle ;
\draw  [fill={rgb, 255:red, 208; green, 2; blue, 27 }  ,fill opacity=1 ] (144,134) .. controls (144,130.13) and (147.13,127) .. (151,127) .. controls (154.87,127) and (158,130.13) .. (158,134) .. controls (158,137.87) and (154.87,141) .. (151,141) .. controls (147.13,141) and (144,137.87) .. (144,134) -- cycle ;
\draw  [line width=1.5]  (294.65,84.7) -- (298.66,98.63) -- (285.12,93.44) ;
\draw  [line width=1.5]  (258.31,149.75) -- (261.3,163.94) -- (248.18,157.78) ;
\draw  [line width=1.5]  (301.39,190.7) -- (287.27,193.98) -- (293.16,180.73) ;
\draw  [line width=1.5]  (202.36,200.77) -- (194.77,188.42) -- (209.2,189.8) ;

\draw (209,122) node [anchor=north west][inner sep=0.75pt]  [font=\large] [align=left] {$\displaystyle \Omega ^{1}$};
\draw (122,209) node [anchor=north west][inner sep=0.75pt]  [font=\large] [align=left] {$\displaystyle \Omega ^{2}$};
\draw (189,63) node [anchor=north west][inner sep=0.75pt]  [font=\large,color={rgb, 255:red, 0; green, 0; blue, 0 }  ,opacity=1 ] [align=left] {$\displaystyle \Su $};
\draw (325,35) node [anchor=north west][inner sep=0.75pt]  [font=\large] [align=left] {$\displaystyle \meshDomain $};
\draw (151,64) node [anchor=north west][inner sep=0.75pt]  [font=\large] [align=left] {$\displaystyle n^{1}$};
\draw (180,151) node [anchor=north west][inner sep=0.75pt]  [font=\large] [align=left] {$\displaystyle n^{2}$};

\end{tikzpicture}
  \caption{Schematic representation of an unfitted NURBS interface $\Su$ in the computational domain $\meshDomain$.}
    \label{fig:domainExamples}
\end{figure}
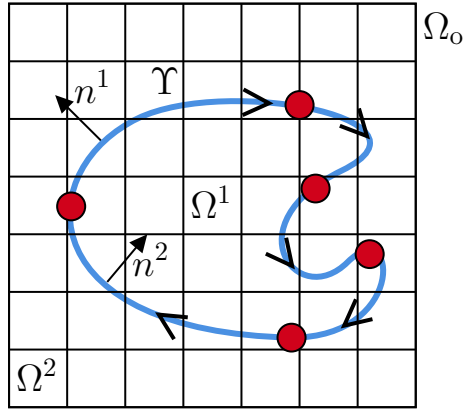

Given the unfitted NURBS interface, see \eqref{eq:nurbsDef},  the two fluid subdomains $\Omega^i$, $i=1,2$ need to be characterized. Without loss of generality, only the interface $\Su$ splitting domains $\Omega^1$ and $\Omega^2$ is considered (i.e., assume that the external boundary $\partial\Omega$ is aligned with the mesh boundary $\partial\meshDomain$). Moreover, the interface $\Su$ is assumed closed, orientable, and not touching the boundary $\partial\Omega$. Hence, the interior subdomain is denoted by $\Omega^1$ whereas the exterior one is described as $\Omega^2$. 
A schematic representation of the setup is reported in Figure\ \ref{fig:domainExamples}.  The interface $\Su$ is clockwise oriented and is composed of $5$ NURBS curves separated by red bullets. The outward unit normal to $\Omega^i$ is denoted by $\bn^i$, for $i=1,2$. 
%

To identify the two fluid domains,  this work follows the procedure presented in\ \cite{Navarro:21}.
The idea is to decompose the problem element-wise by identifying, for each element, the regions $\Omega_e^i$ for $i=1,2$ and then obtain each fluid domain as $\Omega^i = \bigSum{e}{1}{\numel} \Omega_e^i$.

First, the elements cut by the NURBS are identified and the intersection points between the internal skeleton $\Gamma$ and $\Su$ are located and ordered according to the orientation of $\Su$ (see Appendix\ \ref{ap:intersect} for the details of this step).
Consequently, $\Su_e$ is defined as the portion of $\Su$ inside element $\Omega_e$, i.e., $\Su_e := \Su \cap \Omega_e$. Note that, by definition, $\Su_e$ is a continuous curve and it can be composed of several NURBS curves, that is,
\begin{equation}
    \Su_e = \bigcup_{j\in\mathcal{J}_e} \Su_{e}^{j} \quad \text{ with $\Su_{e}^{j} :=\pmb{C}^j_{\Su}([\lambda_{e,1}^j,\lambda_{e,2}^j])$,}
\end{equation}
where $\mathcal{J}_e$ is the set of coefficients $j$ of the NURBS curves $\mathcal{\pmb{C}}^j_{\Su}$ inside $\Omega_e$, while $\lambda_{e,1}^j$ and $\lambda_{e,2}^j$ are the extrema of $\pmb{C}^j_{\Su}$ in $\Omega_e$. 

\begin{figure}[!h]
    \centering
\subfloat[Identification of $\Omega_e^1$ and $\Omega_e^2$.]{

\tikzset{every picture/.style={line width=0.75pt}} 

\begin{tikzpicture}[x=0.75pt,y=0.75pt,yscale=-0.8,xscale=0.8]

\draw [color={rgb, 255:red, 74; green, 144; blue, 226 }  ,draw opacity=1 ][line width=2.25]    (191,96) .. controls (192,97) and (248,207) .. (335,143) .. controls (422,79) and (466,173) .. (463,172) ;
\draw  [line width=1.5]  (216,58) -- (449,58) -- (449,281) -- (216,281) -- cycle ;
\draw  [fill={rgb, 255:red, 208; green, 2; blue, 27 }  ,fill opacity=1 ] (442,151) .. controls (442,147.13) and (445.13,144) .. (449,144) .. controls (452.87,144) and (456,147.13) .. (456,151) .. controls (456,154.87) and (452.87,158) .. (449,158) .. controls (445.13,158) and (442,154.87) .. (442,151) -- cycle ;
\draw  [fill={rgb, 255:red, 184; green, 233; blue, 134 }  ,fill opacity=1 ] (328,143) .. controls (328,139.13) and (331.13,136) .. (335,136) .. controls (338.87,136) and (342,139.13) .. (342,143) .. controls (342,146.87) and (338.87,150) .. (335,150) .. controls (331.13,150) and (328,146.87) .. (328,143) -- cycle ;
\draw  [fill={rgb, 255:red, 74; green, 144; blue, 226 }  ,fill opacity=1 ] (442,58) .. controls (442,54.13) and (445.13,51) .. (449,51) .. controls (452.87,51) and (456,54.13) .. (456,58) .. controls (456,61.87) and (452.87,65) .. (449,65) .. controls (445.13,65) and (442,61.87) .. (442,58) -- cycle ;
\draw  [fill={rgb, 255:red, 208; green, 2; blue, 27 }  ,fill opacity=1 ] (209,126) .. controls (209,122.13) and (212.13,119) .. (216,119) .. controls (219.87,119) and (223,122.13) .. (223,126) .. controls (223,129.87) and (219.87,133) .. (216,133) .. controls (212.13,133) and (209,129.87) .. (209,126) -- cycle ;
\draw  [fill={rgb, 255:red, 74; green, 144; blue, 226 }  ,fill opacity=1 ] (442,281) .. controls (442,277.13) and (445.13,274) .. (449,274) .. controls (452.87,274) and (456,277.13) .. (456,281) .. controls (456,284.87) and (452.87,288) .. (449,288) .. controls (445.13,288) and (442,284.87) .. (442,281) -- cycle ;
\draw  [fill={rgb, 255:red, 74; green, 144; blue, 226 }  ,fill opacity=1 ] (209,281) .. controls (209,277.13) and (212.13,274) .. (216,274) .. controls (219.87,274) and (223,277.13) .. (223,281) .. controls (223,284.87) and (219.87,288) .. (216,288) .. controls (212.13,288) and (209,284.87) .. (209,281) -- cycle ;
\draw  [line width=1.5]  (280.47,168.6) -- (267.79,161.56) -- (281.04,155.68) ;
\draw  [line width=1.5]  (407.47,127.6) -- (394.79,120.56) -- (408.04,114.68) ;
\draw  [line width=1.5]  (209.81,98.34) -- (216.27,85.36) -- (222.74,98.34) ;
\draw  [line width=1.5]  (455.24,101.36) -- (448.77,114.34) -- (442.31,101.36) ;
\draw  [line width=1.5]  (319.9,52.27) -- (332.76,58.96) -- (319.67,65.2) ;
\draw  [line width=1.5]  (319.73,274.39) -- (332.71,280.85) -- (319.73,287.32) ;
\draw  [fill={rgb, 255:red, 74; green, 144; blue, 226 }  ,fill opacity=1 ] (209,58) .. controls (209,54.13) and (212.13,51) .. (216,51) .. controls (219.87,51) and (223,54.13) .. (223,58) .. controls (223,61.87) and (219.87,65) .. (216,65) .. controls (212.13,65) and (209,61.87) .. (209,58) -- cycle ;
\draw  [line width=1.5]  (222.74,221.36) -- (216.27,234.34) -- (209.81,221.36) ;
\draw  [line width=1.5]  (442.31,229.34) -- (448.77,216.36) -- (455.24,229.34) ;
\draw  [fill={rgb, 255:red, 255; green, 255; blue, 255 }  ,fill opacity=1 ] (385.5,101.75) .. controls (385.5,93.6) and (364.01,87) .. (337.5,87) -- (337.5,72) .. controls (364.01,72) and (385.5,78.6) .. (385.5,86.75) ;\draw  [fill={rgb, 255:red, 255; green, 255; blue, 255 }  ,fill opacity=1 ] (385.5,86.75) .. controls (385.5,92.8) and (373.65,98) .. (356.7,100.27) -- (356.7,95.77) -- (337.5,109) -- (356.7,119.77) -- (356.7,115.27) .. controls (373.65,113) and (385.5,107.8) .. (385.5,101.75)(385.5,86.75) -- (385.5,101.75) ;
\draw  [fill={rgb, 255:red, 255; green, 255; blue, 255 }  ,fill opacity=1 ] (288.5,86.25) .. controls (288.5,94.4) and (309.99,101) .. (336.5,101) -- (336.5,116) .. controls (309.99,116) and (288.5,109.4) .. (288.5,101.25) ;\draw  [fill={rgb, 255:red, 255; green, 255; blue, 255 }  ,fill opacity=1 ] (288.5,101.25) .. controls (288.5,95.2) and (300.35,90) .. (317.3,87.73) -- (317.3,92.23) -- (336.5,79) -- (317.3,68.23) -- (317.3,72.73) .. controls (300.35,75) and (288.5,80.2) .. (288.5,86.25)(288.5,101.25) -- (288.5,86.25) ;
\draw  [fill={rgb, 255:red, 255; green, 255; blue, 255 }  ,fill opacity=1 ] (385.5,215.25) .. controls (385.5,223.4) and (364.01,230) .. (337.5,230) -- (337.5,245) .. controls (364.01,245) and (385.5,238.4) .. (385.5,230.25) ;\draw  [fill={rgb, 255:red, 255; green, 255; blue, 255 }  ,fill opacity=1 ] (385.5,230.25) .. controls (385.5,224.2) and (373.65,219) .. (356.7,216.73) -- (356.7,221.23) -- (337.5,208) -- (356.7,197.23) -- (356.7,201.73) .. controls (373.65,204) and (385.5,209.2) .. (385.5,215.25)(385.5,230.25) -- (385.5,215.25) ;
\draw  [fill={rgb, 255:red, 255; green, 255; blue, 255 }  ,fill opacity=1 ] (289.5,230.75) .. controls (289.5,222.6) and (310.99,216) .. (337.5,216) -- (337.5,201) .. controls (310.99,201) and (289.5,207.6) .. (289.5,215.75) ;\draw  [fill={rgb, 255:red, 255; green, 255; blue, 255 }  ,fill opacity=1 ] (289.5,215.75) .. controls (289.5,221.8) and (301.35,227) .. (318.3,229.27) -- (318.3,224.77) -- (337.5,238) -- (318.3,248.77) -- (318.3,244.27) .. controls (301.35,242) and (289.5,236.8) .. (289.5,230.75)(289.5,215.75) -- (289.5,230.75) ;

\draw (451,68) node [anchor=north west][inner sep=0.75pt]   [align=left] {$\displaystyle \Omega _{e}$};
\draw (254,86) node [anchor=north west][inner sep=0.75pt]   [align=left] {$\displaystyle \Omega _{e}^{1}$};
\draw (408,241) node [anchor=north west][inner sep=0.75pt]   [align=left] {$\displaystyle \Omega _{e}^{2}$};
\draw (310,160) node [anchor=north west][inner sep=0.75pt]   [align=left] {$\displaystyle \Su _{e}$};

\end{tikzpicture}

\label{fig:decompositionOmegaeA}
}
\quad
\subfloat[Triangulation for quadrature.]{

\tikzset{every picture/.style={line width=0.75pt}} 

\begin{tikzpicture}[x=0.75pt,y=0.75pt,yscale=-0.8,xscale=0.8]

\draw    (449,281) -- (335,143) ;
\draw    (276,280) -- (216,126) ;
\draw    (276,280) -- (335,143) ;
\draw    (449,58) -- (335,143) ;
\draw [color={rgb, 255:red, 74; green, 144; blue, 226 }  ,draw opacity=1 ][line width=2.25]    (191,96) .. controls (192,97) and (248,207) .. (335,143) .. controls (422,79) and (466,173) .. (463,172) ;
\draw  [line width=1.5]  (216,58) -- (449,58) -- (449,281) -- (216,281) -- cycle ;
\draw  [fill={rgb, 255:red, 208; green, 2; blue, 27 }  ,fill opacity=1 ] (442,151) .. controls (442,147.13) and (445.13,144) .. (449,144) .. controls (452.87,144) and (456,147.13) .. (456,151) .. controls (456,154.87) and (452.87,158) .. (449,158) .. controls (445.13,158) and (442,154.87) .. (442,151) -- cycle ;
\draw    (216,126) -- (449,58) ;
\draw  [fill={rgb, 255:red, 184; green, 233; blue, 134 }  ,fill opacity=1 ] (328,143) .. controls (328,139.13) and (331.13,136) .. (335,136) .. controls (338.87,136) and (342,139.13) .. (342,143) .. controls (342,146.87) and (338.87,150) .. (335,150) .. controls (331.13,150) and (328,146.87) .. (328,143) -- cycle ;
\draw  [fill={rgb, 255:red, 74; green, 144; blue, 226 }  ,fill opacity=1 ] (442,58) .. controls (442,54.13) and (445.13,51) .. (449,51) .. controls (452.87,51) and (456,54.13) .. (456,58) .. controls (456,61.87) and (452.87,65) .. (449,65) .. controls (445.13,65) and (442,61.87) .. (442,58) -- cycle ;
\draw  [fill={rgb, 255:red, 208; green, 2; blue, 27 }  ,fill opacity=1 ] (209,126) .. controls (209,122.13) and (212.13,119) .. (216,119) .. controls (219.87,119) and (223,122.13) .. (223,126) .. controls (223,129.87) and (219.87,133) .. (216,133) .. controls (212.13,133) and (209,129.87) .. (209,126) -- cycle ;
\draw  [fill={rgb, 255:red, 74; green, 144; blue, 226 }  ,fill opacity=1 ] (442,281) .. controls (442,277.13) and (445.13,274) .. (449,274) .. controls (452.87,274) and (456,277.13) .. (456,281) .. controls (456,284.87) and (452.87,288) .. (449,288) .. controls (445.13,288) and (442,284.87) .. (442,281) -- cycle ;
\draw  [fill={rgb, 255:red, 74; green, 144; blue, 226 }  ,fill opacity=1 ] (269,280) .. controls (269,276.13) and (272.13,273) .. (276,273) .. controls (279.87,273) and (283,276.13) .. (283,280) .. controls (283,283.87) and (279.87,287) .. (276,287) .. controls (272.13,287) and (269,283.87) .. (269,280) -- cycle ;
\draw  [fill={rgb, 255:red, 74; green, 144; blue, 226 }  ,fill opacity=1 ] (209,281) .. controls (209,277.13) and (212.13,274) .. (216,274) .. controls (219.87,274) and (223,277.13) .. (223,281) .. controls (223,284.87) and (219.87,288) .. (216,288) .. controls (212.13,288) and (209,284.87) .. (209,281) -- cycle ;
\draw  [fill={rgb, 255:red, 74; green, 144; blue, 226 }  ,fill opacity=1 ] (209,58) .. controls (209,54.13) and (212.13,51) .. (216,51) .. controls (219.87,51) and (223,54.13) .. (223,58) .. controls (223,61.87) and (219.87,65) .. (216,65) .. controls (212.13,65) and (209,61.87) .. (209,58) -- cycle ;

\draw (451,68) node [anchor=north west][inner sep=0.75pt]   [align=left] {$\displaystyle \Omega _{e}$};
\draw (410,95) node [anchor=north west][inner sep=0.75pt]   [align=left] {$\displaystyle T_{e}^{1,1}$};
\draw (251,71) node [anchor=north west][inner sep=0.75pt]   [align=left] {$\displaystyle T_{e}^{1,3}$};
\draw (293,114) node [anchor=north west][inner sep=0.75pt]   [align=left] {$\displaystyle T_{e}^{1,2}$};
\draw (397,162) node [anchor=north west][inner sep=0.75pt]   [align=left] {$\displaystyle T_{e}^{2,1}$};
\draw (269,179) node [anchor=north west][inner sep=0.75pt]   [align=left] {$\displaystyle T_{e}^{2,3}$};
\draw (340,211) node [anchor=north west][inner sep=0.75pt]   [align=left] {$\displaystyle T_{e}^{2,2}$};
\draw (223,234) node [anchor=north west][inner sep=0.75pt]   [align=left] {$\displaystyle T_{e}^{2,4}$};

\end{tikzpicture}

\label{fig:decompositionOmegaeB}
}
    \caption{Element $\Omega_e$ cut in two regions by an interface $\Su_e$.
    (a) Identification of $\Omega_e^1$ and $\Omega_e^2$.
    (b) Triangulation for quadrature. 
    The interface is composed of two NURBS curves connected at the green bullet.  The red bullets denote the intersections between $\Su$ and $\partial\Omega_e$ and the blue ones are the vertices of $\Omega_e$ and of the triangulation for qudrature.}
    \label{fig:decompositionOmegae}
\end{figure}
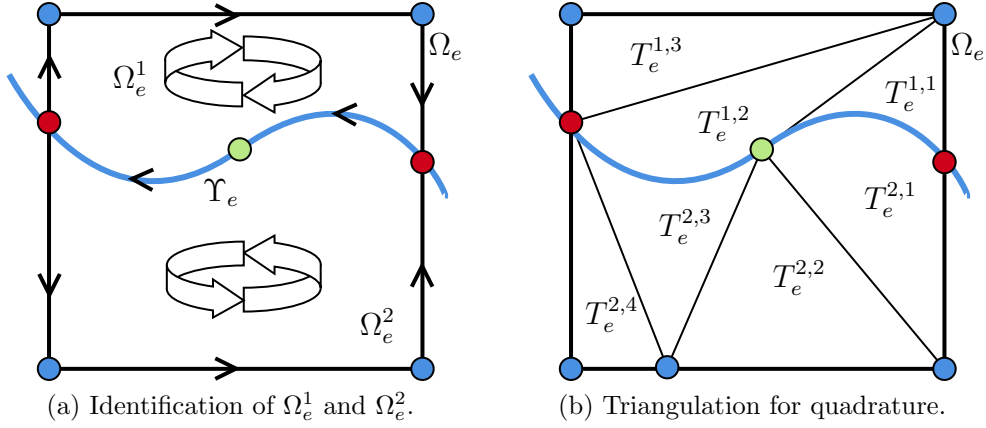

For each cut element $\Omega_e$ (see, e.g., the green element in Figure\ \ref{fg:cut_cell_ex}), the two regions, $\Omega_e^1$ and $\Omega_e^2$, are constructed connecting the mesh vertices and the ordered intersection points. 
Starting from the first intersection point between $\Gamma$ and $\Su_e$ and travelling along $\Su_e$ towards the left or towards the right, two sides of $\Su_e$ are identified. Depending on the orientation of $\Su_e$, the fluid index $i$ of each adjacent region is determined: when $\Su_e$ is clockwise oriented, the right region is $\Omega_e^1$,  whereas the left one is denoted by $\Omega_e^2$, and vice versa when $\Su_e$ is counter-clockwise oriented. 
This procedure is sketched in Figure\ \ref{fig:decompositionOmegaeA} and it can also be generalized to the case of unfitted boundary $\partial\Omega$.

The assignment of the appropriate fluid index to the remaining, uncut elements is performed iteratively. Starting from a neighboring element, with a common face, for which such information has already been determined,  the fluid index of the uncut element is inherited across the shared face.

\subsection{Quadrature in elements and faces cut by NURBS}
Specific numerical techniques are used to integrate along cut faces (i.e., along NURBS curves) and in cut elements, whereas uncut faces and uncut elements employ standard HDG quadratures. 

Quadratures along cut faces only require the proper mapping between each segment of the face and the reference domain of the chosen quadrature. Typically Gauss-Legendre quadratures are used along each segment.

The quadrature along NURBS curves is treated element-wise. For a generic scalar function $f$, the quadrature along the portion of NURBS interface $\Su_e$ is defined as
\begin{equation*}\label{eq:exact_lineInt}
   \int_{\Su_e} \! f\, \mathrm{d} l = \sum_{j\in\mathcal{J}_e}
                                                 \int_{\lambda^{j}_1}^{\lambda^{j}_2} f\bigl(\pmb{C}^j_{\Su}\!(\lambda)\bigr)\, \norm{J_{\pmb{C}^j_{\Su}}\!(\lambda)}\, \mathrm{d}\lambda,
\end{equation*}
where $\norm{J_{{\pmb{C}}^j_{\Su}}\!}$ is the norm of the differential of the parametrization of the NURBS curve $\pmb{C}^j_{\Su}$ describing $\Su_e^j$. 
The integrals are approximated using the $1$D Gauss--Legendre quadrature defined in the parametric space $[\lambda_{e,1}^j,\lambda_{e,2}^j]$ and mapped along $\pmb{C}^j_{\Su}$. The same procedure can be employed also for unfitted external boundaries. 

The quadrature within a generic region $\Omega_e^i$ involves multiple steps. Firstly, the region $\Omega_e^i$ undergoes a partitioning process, guided by Lee's visibility algorithm\ \cite{Joe:87}, which is essential to ensure the accuracy of the quadrature, as emphasized in\ \cite{Sevilla-SRH-16}. 
The procedure, presented in Appendix\ \ref{ap:quadTri}, constructs a triangular partition  (see Figure\ \ref{fig:decompositionOmegaeB}) of each subregion $\Omega_e^i$ such that
\begin{equation*}
    \Omega_e^i = \bigSum{l}{1}{\numelOmega^i} T_e^{i,l},
\end{equation*}
with $\numelOmega^i$ the total number of triangles composing $\Omega_e^i$.  
The resulting triangulation consists of two types of elements: affine and curved. 
The former are handled using standard isoparametric mapping in finite elements.
The latter, which require a special treatment, are constructed with at most one curved edge, described by one NURBS curve.

Given a generic scalar function $f$, the quadrature in $\Omega_e^i$ is obtained as the sum of the quadrature in each triangle identified according to the previous procedure, that is, 
\begin{equation*}
  \int_{\Omega_e^i}\! f(x,y)\, \mathrm{d}x \mathrm{d}y = \bigSumE{l}{1}{\numelOmega^i}\int_{T_e^{i,l}} f(x,y)\, \mathrm{d}x \mathrm{d}y.
\end{equation*}
When $T_e^{i,l}$ is an affine triangle, the standard Gauss--Legendre quadrature is employed, whereas the quadrature procedure proposed in\ \cite{Sevilla-SFH-08-IJNME} for NEFEM is employed in the case of a curved element $T_e^{i,l}$. In this case, the Gauss--Legendre quadrature is constructed in the rectangle $R{:=}[\lambda_{e,1}^j,\lambda_{e,2}^j]{\times}[0, 1]$ which is transformed to $T_e^{i,l}$ by means of the affine mapping 
\begin{equation*}\label{eq:varphi}
        \bm{\psi}(\lambda,\theta) :=
       (1-\theta) \mathcal{\pmb{C}}^j\!(\lambda) + \theta\, \pmb{x}_3, \quad \forall \lambda \in [\lambda_{e,1}^j,\lambda_{e,2}^j], \quad \forall \theta \in [0, 1],
    \end{equation*}
with $\bx_3$ being the vertex of the triangle not lying on the NURBS curve. Then, it follows that
\begin{equation*}\label{eq:elem_quad}
  \int_{T_e^{i,l}}\! f(x,y) \,\mathrm{d}x \mathrm{d}y = \int_{R}\! f\bigl(\bm{\psi}(\lambda,\theta) \bigr) \abs{J_{\pmb{\psi}}(\lambda,\theta)} \, \mathrm{d}\lambda \mathrm{d}\theta,
\end{equation*}
where $|J_{\bm{\psi}}|$ is the determinant of the Jacobian of the transformation $\bm{\psi}$.

%
\subsection{Element extension for badly cut cells}
\label{sc:strategiesEE_SC}
It is well-known that unfitted methods suffer in the presence of elements, cut by the interface or the boundary, such that the portion of the computational volume is small with respect to the volume of the entire element. 
To measure such cuts, the parameter 
\begin{equation}
    \alpha_e^i := \frac{\abs{\Omega_e^i}}{\abs{\Omega_e}}, \qquad  e\inSet{1}{\numel}, \qquad  i=1,2,
    \label{eq:alpha_cutElem}
\end{equation} 
is defined as the ratio between the region occupied by fluid $i$, i.e., $\Omega_e^i$, and the entire element $\Omega_e$. The smallest cut in a computational domain $\meshDomain$ is denoted by $\alpha := \min_e \min_i \alpha_e^i$. 
An element is said to be badly cut when $\alpha_e^i<\alpha_{\min}$, with $\alpha_{\min}$ a user-defined value, usually set to $\alpha_{\min} {=} 0.3$ consistently with the literature, see, e.g.\ \cite{BuCDE:21,BuDeE:21}. 

\begin{figure}[!h]
    \centering
    \includegraphics[width=0.75\linewidth]{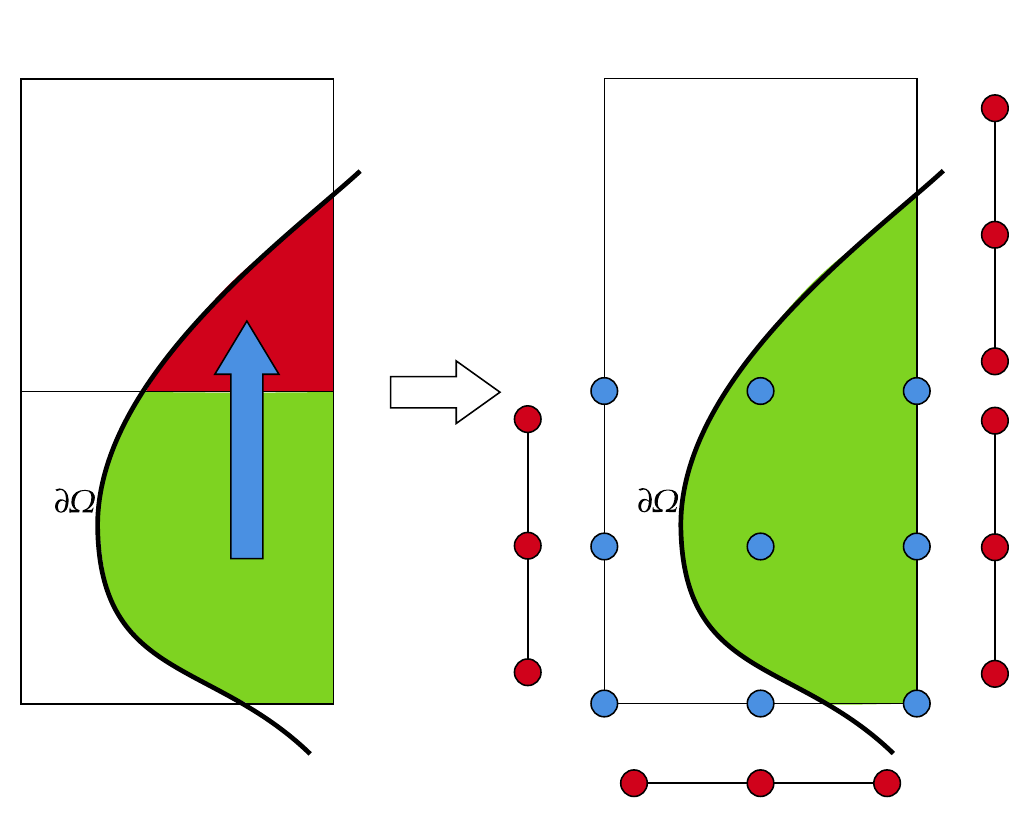}
    \caption[Element extension strategy.]{Element extension strategy. Left: initial configuration with badly cut cell (red) and well-cut cell (green). Right: extended element. Blue bullets are the element-based unknowns and red bullets the face-based unknowns of the resulting extended element.}
    \label{fig:exampleEE}
\end{figure}

To address badly cut cells, the element extension method proposed in \cite{Navarro:21} is used.
If $\alpha_e^i<\alpha_{\min}$ for $i=1,2$ , the element extension strategy, depicted in Figure\ \ref{fig:exampleEE}, performs the integration using the shape functions and elemental nodal distribution of a well-cut, neighboring element. More precisely, the procedure executes the following operations:
\begin{enumerate}[label=(\arabic*)]
\item It erases the element-based unknowns of the badly cut cell.
\item It extends the element-based approximation (and, hence, the quadrature) from a selected neighboring element to the badly cut cell.
\item It erases the face-based unknowns along the face shared by the two elements.
\item It preserves the face-based unknowns of the badly cut cell.
\end{enumerate}

For each badly cut cell, potential extensions are assessed based on three factors: (i) to maximize the combined area; (ii) to minimize the Euclidean distance between element centers; (iii) to penalize previously extended elements. The element with the highest combined score of these factors is extended onto the badly cut one. In two-fluid problems, a stronger penalty is used to discourage multiple extensions per element owing to inter-region coupling, thus minimizing the local problem size.

%
\section{Numerical results}
\label{sc:numSimu}
In this section,  the unfitted HDG method with exact treatment of NURBS geometry is tested for a series of 2D problems featuring a single fluid and two immiscible fluids. 
After assessing the robustness of the method in the presence of badly cut faces and badly cut cells in Section\ \ref{sc:conditioning}, numerical experiments show the optimal convergence and high-fidelity accuracy achieved by the  scheme in the case of an unfitted external boundary (Section\ \ref{sc:immersedBoundaryResults}) and an unfitted interface (Section\ \ref{sc:immersedInterfaceResults}) described by means of NURBS.
Finally, the suitability of the proposed approach to simulate microfluidic systems featuring complex geometries is showcased in Sections\ \ref{sc:microMix} and\ \ref{sc:emulsion} for one and two-fluid problems.  

\begin{remark}[Extension to three dimensional problems]
Whilst the proposed formulation is seamlessly applicable to problems in two and three dimensions, the treatment of unfitted 3D NURBS surfaces entails some technical difficulties.  On the one hand, to identify the intersections between the NURBS surfaces and the faces/edges of the mesh hexahedra (e.g., via the marching cubes algorithm\ \cite{Lorensen-LC-98}). On the other hand,  to subdivide the resulting three-dimensional cut cells into tetrahedra, as explained in\ \cite{Sevilla-MSZRT-15,Rodenas-MRNT-17}.  Finally, once tetrahedralization of the regions is completed,  the numerical integration in the elements and on the faces cut by the NURBS surfaces can be performed using the 3D version of the NEFEM rationale\ \cite{Sevilla-SFH-11}.
\end{remark}

%
\subsection{Robustness with respect to badly cut faces and cells}
\label{sc:conditioning} 
As noted earlier, it is well know that unfitted methods suffer from badly cut cells. However, in HDG face/skeleton unknowns and integrals play a prominent role. Here, a synthetic problem with manufactured solution is used to study the effect of badly cut faces and badly cut cells on the discretization error and on the condition number of the matrices arising from the HDG formulation.

Consider a one-fluid problem on the domain $\Omega$, with Dirichlet boundary conditions on the unfitted boundary $\partial\Omega$. The detailed definition of the domain $\Omega$ is provided in the following subsections, whereas the computational domain is set to $\meshDomain {=} (0,1)^2$.
Given a viscosity $\mu{=}1$,  the Dirichlet datum $\dirData$ and the source term $\bm{s}$ are analytically obtained from the expressions of the Stokes velocity and pressure fields, see\ \cite{Donea-Huerta}, given by
\begin{equation}\label{eq:sol_testCase}
    \bu^{\text{ref}} = \begin{pmatrix}x^2(1-x)^2(2y-6y^2+4y^3)\\-y^2(1-y)^2(2x-6x^2+4x^3)\end{pmatrix}, \qquad p^{\text{ref}} = x(1-x) \,.
\end{equation}

The quantities of interest in the present study are the condition numbers $\kappaLoc$ and $\kappaGlob$, computed using the Euclidean norm,  of the HDG local and global matrices, respectively.
Consistently with the definition of the element ratio $\alpha_e^i$, see \eqref{eq:alpha_cutElem}, for each face $\Gamma_{\!\! f} \subset \Gamma$, the ratio between the computational face $\Gamma_{\!\! f} \cap \Omega^i$ and the entire face $\Gamma_{\!\! f}$ is defined as
%
\begin{equation}
    \beta^i_{\! f} := \frac{\abs{\Gamma_{\!\! f} \cap \Omega^i}}{\abs{ \Gamma_{\!\! f}}} , \qquad  f\inSet{1}{\numelfc}, \quad i \inSetTwo{1}{2} \, ,
\end{equation}
and $\beta := \min_f \min_i \beta^i_{\! f}$ denotes the smallest face cut in the computational domain.

%
\subsubsection{Badly cut faces}
\label{sc:badlyfaces}
\begin{figure}[!h]
 \centering
 {\includegraphics[width=0.35\textwidth]{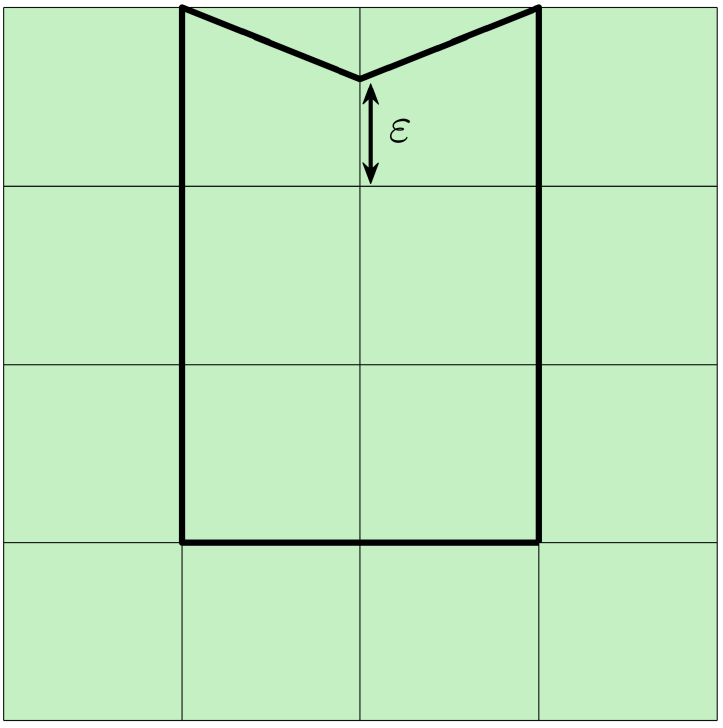}}
 \caption{M-shaped domain with $\varepsilon = 1.5 \times 10^{-1}$ ($\beta = 60\%$ of the entire face).}
 \label{fg:figMdomain}
\end{figure}

Consider the M-shaped domain (see Figure\ \ref{fg:figMdomain}),  defined as $\Omega := \Omega_{\textrm{R}} \setminus \Omega_{\textrm{T}}$, with $\Omega_{\textrm{R}} {=} (0.25,0.75)\times(0.25,1)$ and $\Omega_{\textrm{T}}$ being the triangular region with vertices $(0.25,1)$, $(0.75,1)$ and $(0.5,0.75+\varepsilon)$, $\varepsilon {>}0$.
The computational domain is the square $\meshDomain {=} (0,1)^2$, composed of a $4 {\times} 4$ mesh of square elements.

By construction, all the elements and faces of the M-shaped domain are well-cut, except for the vertical face of extrema $(0.5,0.75)$ and $(0.5,0.75+\varepsilon)$. 
The influence on the global condition number of the decreasing length of this face is studied for $\varepsilon\in \{ 5 \times 10^{-3}, 5 \times 10^{-2},10^{-1},1.5 \times 10^{-1} \}$, leading to $\beta \in \{2\%,20\%,40\%,60\% \}$ of the entire face, while always maintaining well-cut elements with $\alpha > 50\%$.
%

\begin{figure}[!h]
 \centering
 \includegraphics[width=\textwidth]{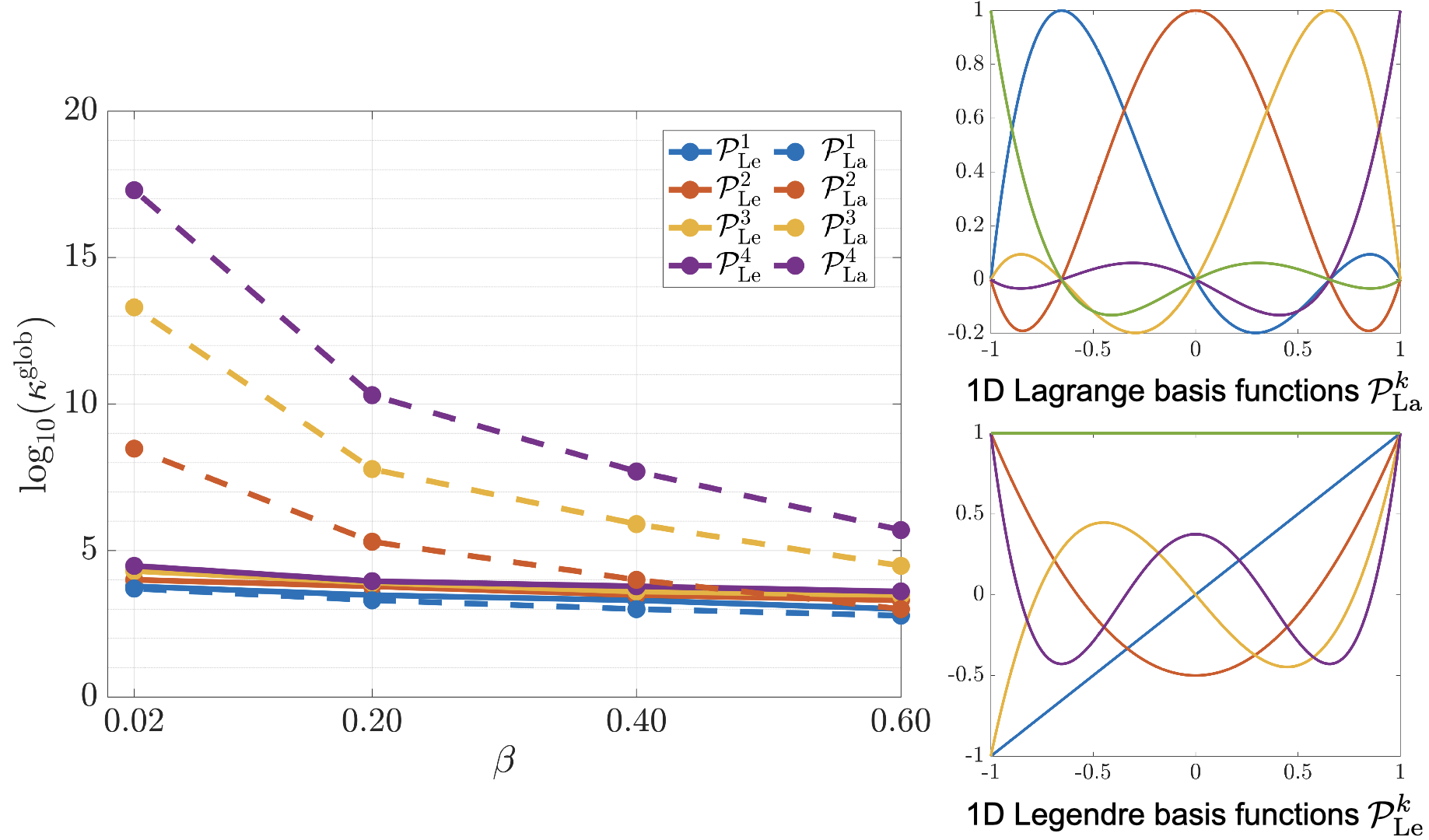}
 \caption{M-shaped domain.  Comparison of the global condition number $\log_{10}(\kappaGlob)$ as a function of the smallest face cut $\beta$ using Lagrange $\Pk_{\text{La}}$ (dashed lines) and Legendre $\Pk_{\text{Le}}$ (continuous lines) basis functions of degree $k\inSet{1}{4}$ for the face-based variable. }
  \label{fg:Kglob_LagVSLeg}%
\end{figure}

Figure\ \ref{fg:Kglob_LagVSLeg} shows the global condition number $\kappaGlob$ as a function of the smallest face cut $\beta$, using  
nodal and modal basis functions of polynomial degree $k\inSet{1}{4}$ for the hybrid variable. More precisely,  Lagrange basis functions $\Pk_{\text{La}}$ are compared with Legendre basis functions $\Pk_{\text{Le}}$. For the element-based unknowns, 2D Lagrange tensor basis functions with Fekete nodal distribution are employed, guaranteeing the optimal condition number and accuracy for standard, uncut, elements.

The results show that, whilst modal basis functions yield a condition number $\kappaGlob$ independent of $\beta$ for any polynomial degree,  nodal basis functions tend to suffer when small cut faces appear in the domain and this phenomenon becomes even more severe for high-order approximations. 
Note that this is consistent with the recent results on the hybrid high order method (a close relative of HDG) presented in\ \cite{Droniou-BDY-22}.
Indeed, by construction, Lagrange basis functions are equal to zero at all nodes except for the one where they are centered. Moreover, the contribution of each basis function progressively concentrates in the vicinity of its center as the polynomial degree $k$ increases (see Figure\ \ref{fg:Kglob_LagVSLeg} for the Lagrange basis functions of degree $4$ in the interval $[-1,1]$).
Thus, when integrating over a portion of a face,  the influence of the Lagrange functions centered at nodes lying outside of the domain significantly decreases, leading to small contributions that tend to zero as the polynomial degree $k$ increases and, consequently, to an ill-conditioning of the matrix.
On the contrary,  Legendre basis functions feature a larger support where each polynomial function is different than zero (see Figure\ \ref{fg:Kglob_LagVSLeg}),  reducing the effect of small face cuts during quadrature. 
Thus, Legendre basis functions are  employed for the face-based variables in all the numerical results presented henceforth.

%
\subsubsection{Badly cut cells}
The adverse effect of badly cut cells on the condition number of immersed boundary methods is well documented in the literature, see, e.g., \cite{Badia-PVBLB-23}.
The objective of this section is to assess the suitability of a standard technique to address this numerical issue (viz. the element extension discussed in \cite{Navarro:21}) in the context of the unfitted HDG method.  Note that alterative approaches such as ghost penalty or cell agglomeration could also be considered\ \cite{Burman:10,JohLa:13}.

\begin{figure}[!h]
 \centering
 \subfloat[Smoothed square domain.] {\includegraphics[width=0.35\textwidth]{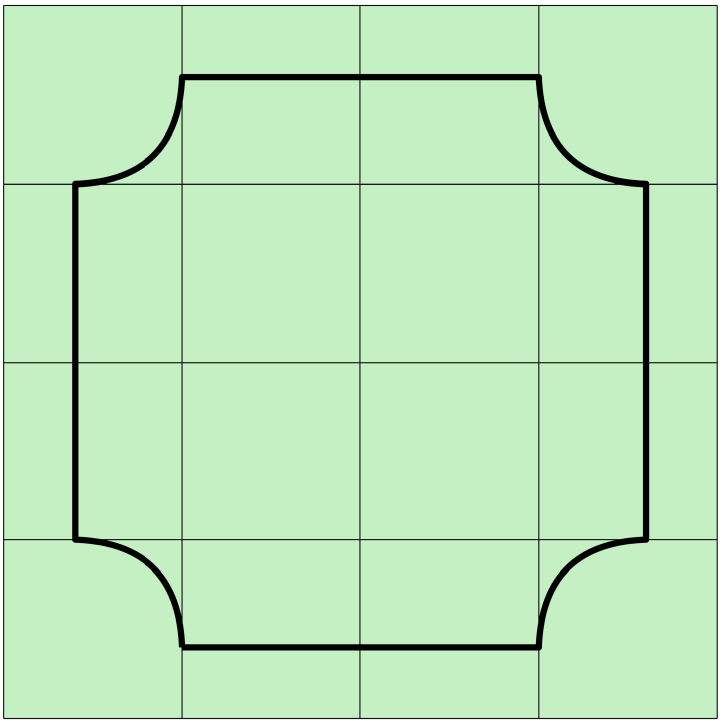}}
 \quad
 \subfloat[Extended elements. \label{fg:smoothSquareDomainExtended}]{\includegraphics[width=0.35\textwidth]{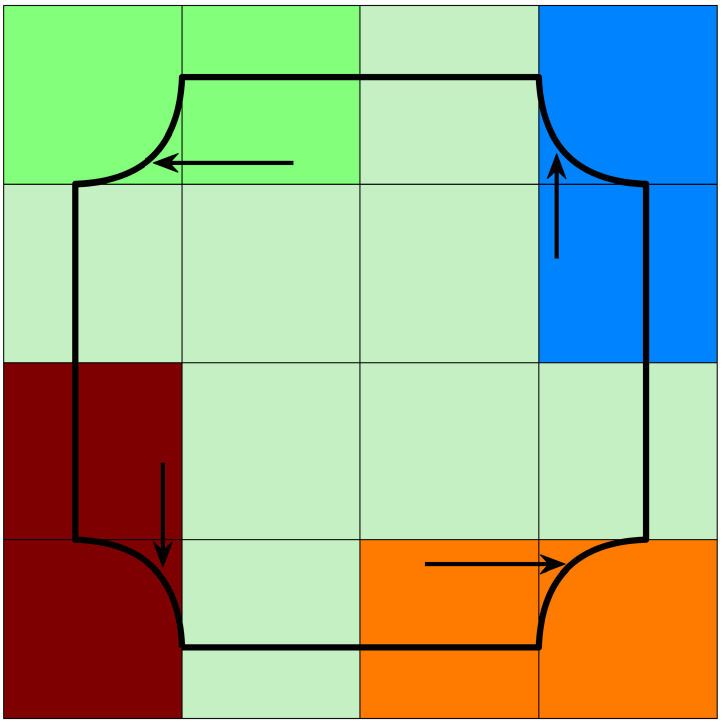}}
 \caption{Smoothed square domain and extension procedure to treat badly cut cells.}
  \label{fg:smoothSquareDomain}%
\end{figure}

Let $\Omega := \Omega_{\textrm{S}} \setminus \bigcup_{i=1}^4 \Omega_{\textrm{A}}^i$ be the domain in Figure\ \ref{fg:smoothSquareDomain}, obtained by subtracting from the square region $\Omega_{\textrm{S}} {=} (0.1,0.9)^2$ the four circular arcs $\Omega_{\textrm{A}}^i, \ i=1,\ldots,4$ of radius $0.15$ and centered in $(0.1,0.1)$,  $(0.9,0.1)$, $(0.9,0.9)$ and $(0.1,0.9)$.
%
%
The computational domain is the square $\meshDomain {=} (0,1)^2$, partitioned in a mesh of $4\times 4$ square elements.
In this configuration, all faces are well-cut, with $\beta > 60\%$.  Nonetheless, whilst most cut elements feature a portion of occupied area equal to $60\%$, the four ones at the corners are badly cut, with $\alpha$ dropping to $6 \%$. 
The extended elements obtained using the strategy described in \cite{Navarro:21} and in Section\ \ref{sc:strategiesEE_SC} are displayed in Figure\ \ref{fg:smoothSquareDomainExtended}.

%

It is worth recalling that the HDG gobal problem, see~\eqref{eq:GlobalIntegral}, only involves face integrals. On the contrary, HDG local problems (see Sections\ \ref{sec:immElem} and\ \ref{sec:intElem}) also feature element integrals and can experience ill-conditioning issues when small elements appear in the computational mesh.
Of course, since the HDG local problems are employed to construct the hybrid system of the HDG global problem, an indirect effect of badly cut cells on the condition number $\kappaGlob$ is also present but this can be controlled by the element extension procedure as reported in Section\ \ref{sc:immersedBoundaryResults}.
Thus,  the following study focuses on the local systems only. 

\begin{figure}[!h]
 \centering
  \subfloat[With element extension.\label{fg:badlyCellsCondResEE}]{\includegraphics[width=0.4\textwidth]{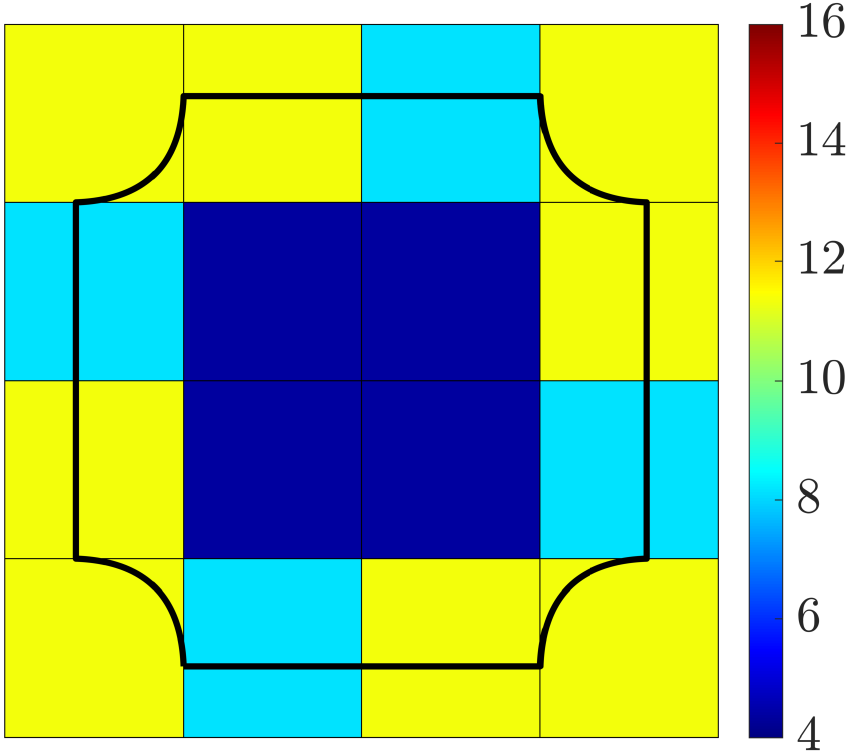}}    
 \quad 
 \subfloat[Without element extension.\label{fg:badlyCellsCondResStd}] {\includegraphics[width=0.4\textwidth]{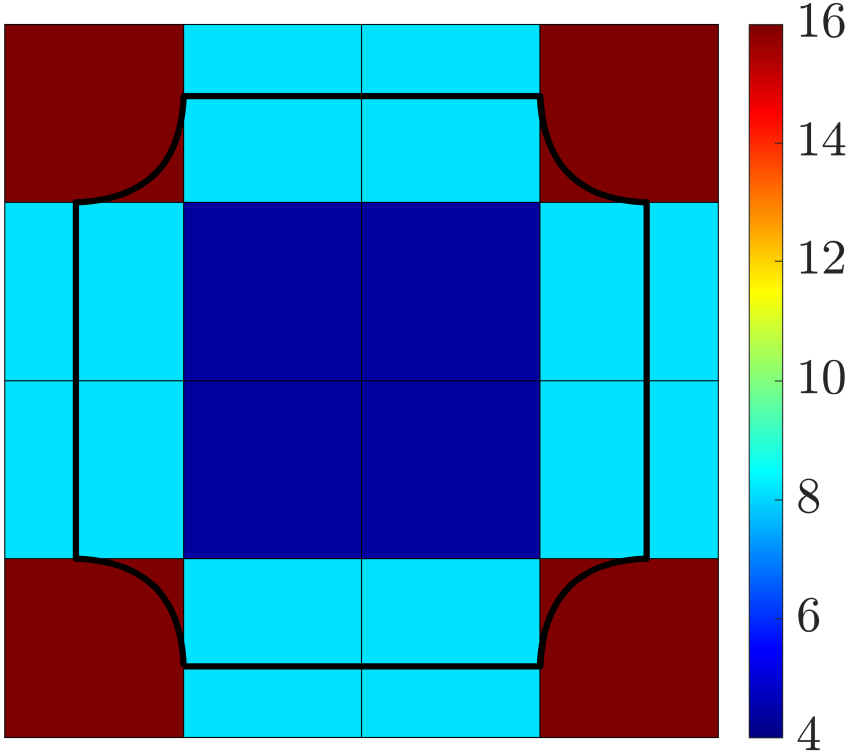}}
 \caption{Smoothed square domain. Element-by-element condition number $\log_{10}(\kappaLoc)$ of the HDG local problem (a) with and (b) without element extension.}
  \label{fg:badlyCellsCondResults}
\end{figure}

Figure\ \ref{fg:badlyCellsCondResults} displays the condition number $\kappaLoc$ for the HDG local problems, with and without element extensions.
As expected,  the badly cut cells are responsible for the local condition number to quickly deteriorate, achieving unmanageable values (Figure\ \ref{fg:badlyCellsCondResStd}).
By means of element extension,  the negative effects of badly cut cells are mitigated, reducing the local condition number of up to $5$ orders of magnitude, as reported in Figure\ \ref{fg:badlyCellsCondResEE}. 
Of course, this procedure also affects the extended, previously well cut, elements, which experience a degradation of the condition number of the local problem: for instance,  an extrapolation over $6\%$ of the element area yields an increase of $\kappaLoc$ of about $2$ orders of magnitude compared to its initial value. 
Nonetheless, as it will be detailed in Section\ \ref{sc:immersedBoundaryResults}, the condition number of the HDG global problem remains bounded even when element extension is performed, thus guaranteeing the robustness of the method and the accuracy of the global solution.

\begin{figure}[!h]
 \centering
 \subfloat[With element extension.\label{fg:badlyCellsVelResEE}]{
    \includegraphics[width=0.4\textwidth]{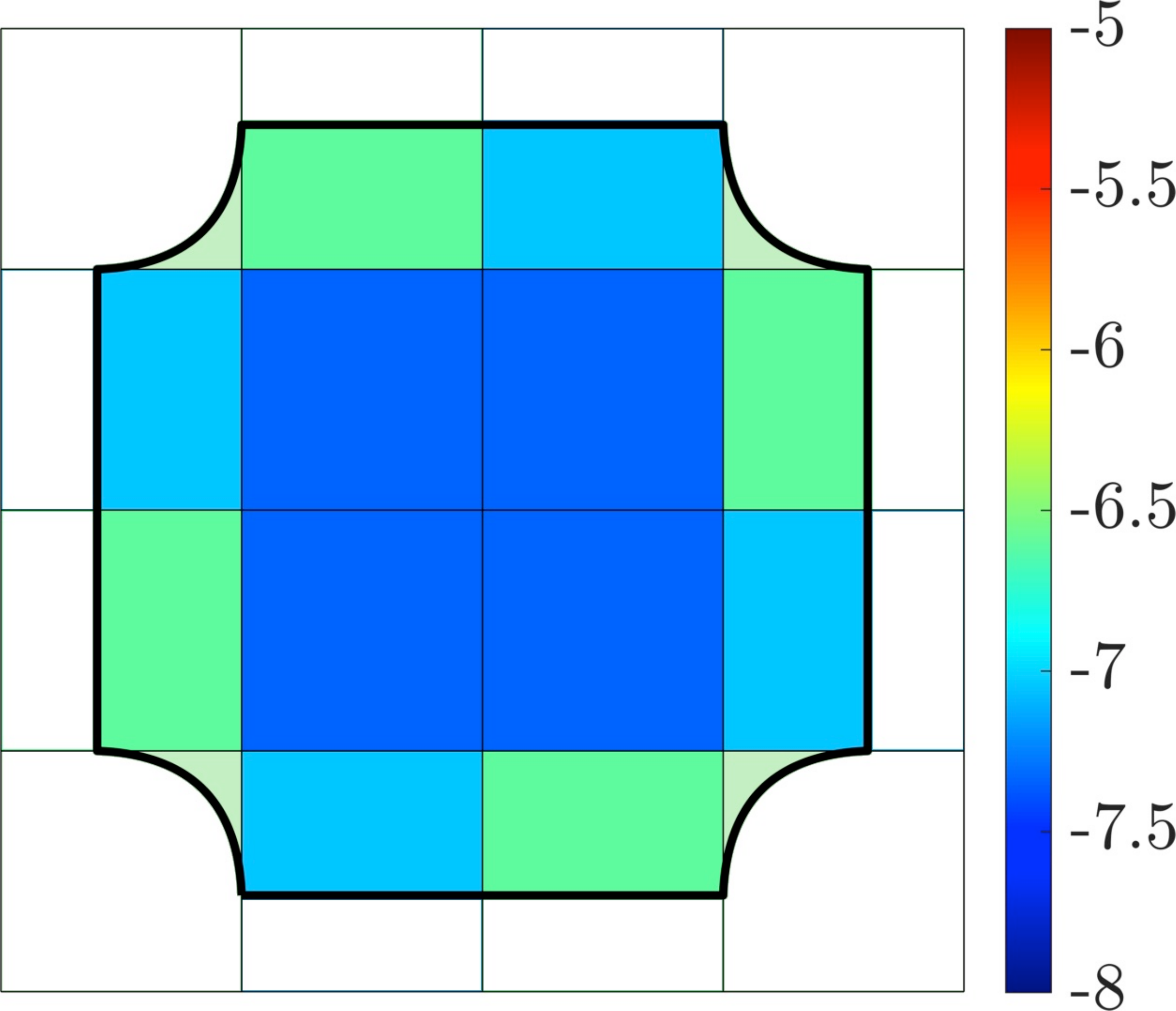}} %
    \quad    
     \subfloat[Without element extension.\label{fg:badlyCellsVelResStd}] {\includegraphics[width=0.4\textwidth]{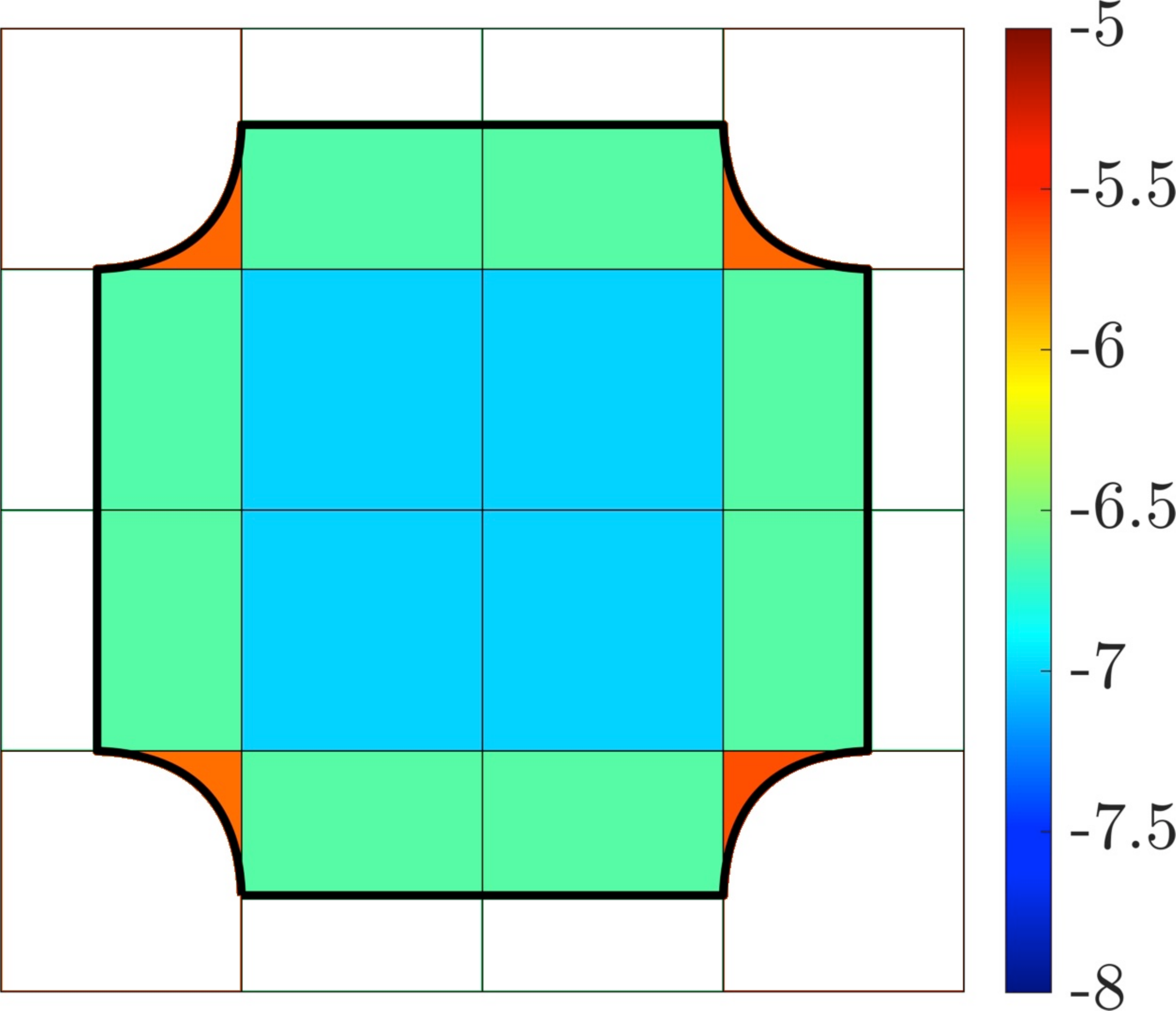}}
 \caption{Smoothed square domain.  Element-by-element map of the error $\log_{10}(\|\bu-\bu^{\text{ref}}\|_{\eltwo(\Omega_e)})$ (a) with and (b) without element extension.}
  \label{fg:badlyCellsVelResults}%
\end{figure}

Similarly, the element-by-element map of the velocity error, measured in the $\eltwo(\Omega_e)$ norm, is presented in Figure\ \ref{fg:badlyCellsVelResults} using a logarithmic scale. 
The results confirm the conclusions from Figure\ \ref{fg:badlyCellsCondResults}, confirming the capability of element extension to significantly reduce the error of the local variable in badly cut cells, while slightly worsening the approximation in the extended region.

The above results are obtained with a discretization employing Lagrange polynomial functions of degree $4$ for the element-based unknowns and a Legendre basis of degree $4$ for the face-based variable.  
Alternative solutions, not reported here for brevity, were explored for the treatment of the element-based unknowns but no significant gain in accuracy nor robustness was observed. 
Thus, Lagrange basis functions with Fekete nodes are henceforth employed to approximate the element-based variables in the following simulations.

%
\subsection{Taylor-Couette flow with unfitted boundary}
\label{sc:immersedBoundaryResults} 
A two-dimensional, one-fluid, coaxial Taylor-Couette flow is considered to numerically assess the accuracy and optimal convergence property of the unfitted HDG method, in the presence of an exact boundary represented via NURBS.
Let us define $(x_c,y_c) {=} (0.5,0.5)$ and $r {=} ((x-x_c)^2+(y-y_c)^2)^{1/2}$.  The domain $\Omega = \{ (x,y) \in \RR^2 \ | \ \Rint \leq r(x,y) \leq \Rext \}$ consists of two coaxial circles of radii $\Rint {=} 1/6$ and $\Rext {=} 1/3$,  centered in $(x_c,y_c)$, and rotating with angular velocities $\OmInt {=} 0$ and $\OmExt {=} 1$, respectively. 
%
%
Setting $\mu {=} 1$,  the Taylor-Couette velocity and pressure fields are defined as
\begin{equation*}
     \bu^{\text{ref}} = (A + B / r^2)\begin{pmatrix} -y+y_c \\[1ex]
                                    \phantom{-} x-x_c \end{pmatrix}, \qquad p^{\text{ref}} = 1,
\end{equation*}
with
\begin{equation*}
    A :=\frac{ \OmExt \Rext^2 - \OmInt \Rint^2}{\Rext^2 -\Rint^2}, \quad B := \frac{(\OmInt  - \OmExt)\Rext^2\Rint^2}{\Rext^2 -\Rint^2} \,  ,
\end{equation*}
and verify the one-fluid Stokes equations with zero body forces and Dirichlet conditions $\dirData$ (obtained from the analytical velocity $\bu^{\text{ref}}$) on the entire boundary $\partial\Omega$, see, e.g.,\ \cite{GiaSeH:21}.

\begin{figure}[!h]
 \centering
 \subfloat[Taylor-Couette domain.\label{fg:TayCouMesh}]{\includegraphics[width=0.31\textwidth]{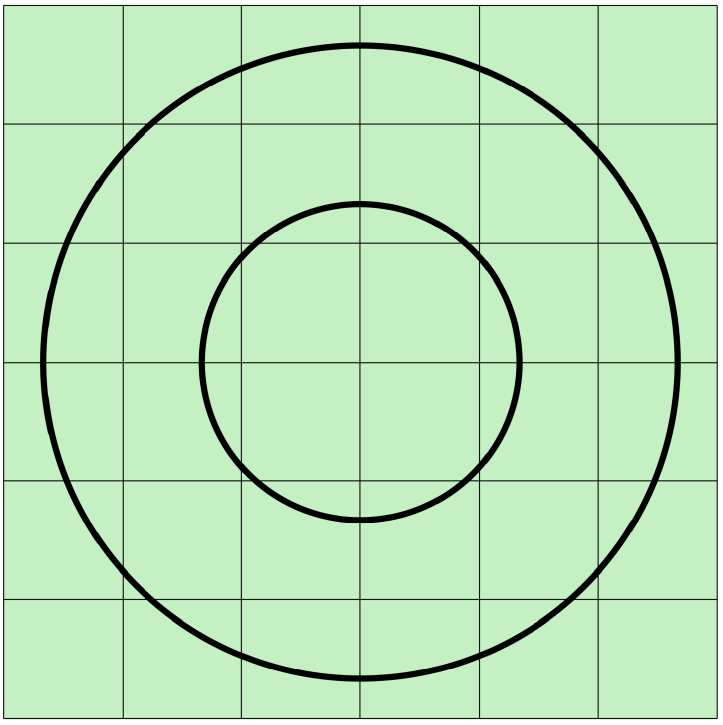}}
 \hfill
 \subfloat[Extended elements.\label{fg:TayCouEE}]{\includegraphics[width=0.31\textwidth]{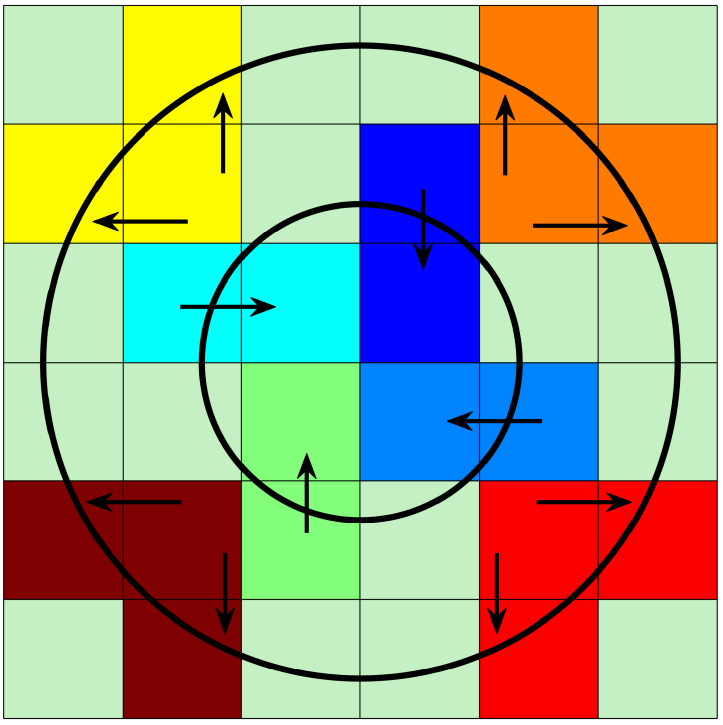}}
\hfill
\subfloat[Module of the velocity.\label{fg:TayCouVel}]{\includegraphics[width=0.36\textwidth]{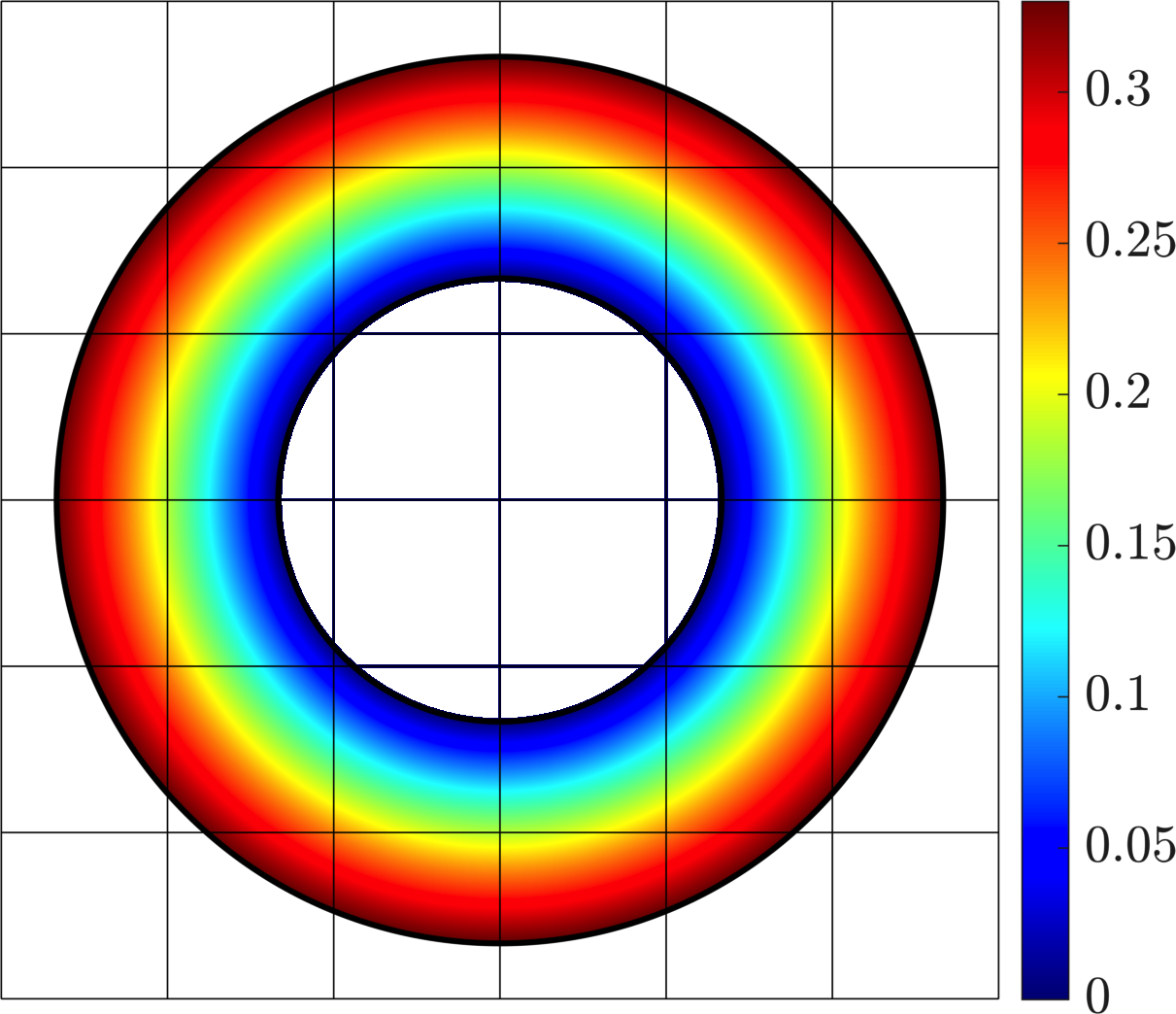}}
 \caption{Taylor-Couette domain. (a) Detail of the computational domain in the region $(0.125,0.875)^2$. (b) Element extension to avoid badly cut cells.  (c) Module of the velocity approximated using polynomial functions of degree $k=4$.}
  \label{fg:couetteDomain}%
\end{figure}

The domain $\Omega$ is immersed in a computational domain $\meshDomain {=} (0,1)^2$, featuring a Cartesian mesh of square elements. 
A detail of the $8\times 8$ mesh, in the region $(0.125,0.875)^2$, is displayed in Figure\ \ref{fg:TayCouMesh}, along with an example of element extension (Figure\ \ref{fg:TayCouEE}).  The module of the velocity field computed on the same mesh using polynomial approximation of degree $4$ is reported in Figure\ \ref{fg:TayCouVel}.

A mesh convergence study is then performed using a set of Cartesian structured meshes composed of $2\times 2$, $4\times 4$, $8\times 8$, $16\times 16$, and $32\times32$ square elements, for polynomial degree approximations ranging from $k=1$ up to $k=4$.

\begin{figure}[!h]
 \centering
 \subfloat[Global problem.]{
    \includegraphics[width=0.4\textwidth]{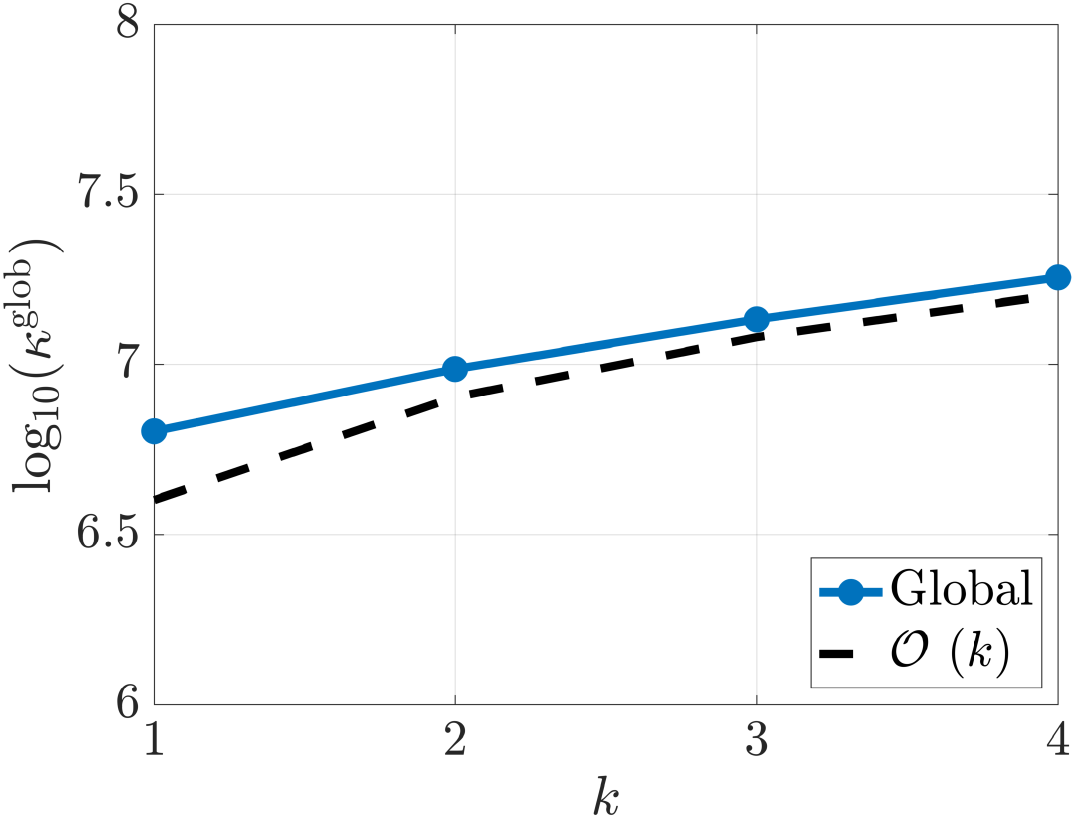}\label{fg:convCouetteKwrt_ka}}
    \qquad
    \subfloat[Local problem.]{
    \includegraphics[width=0.4\textwidth]{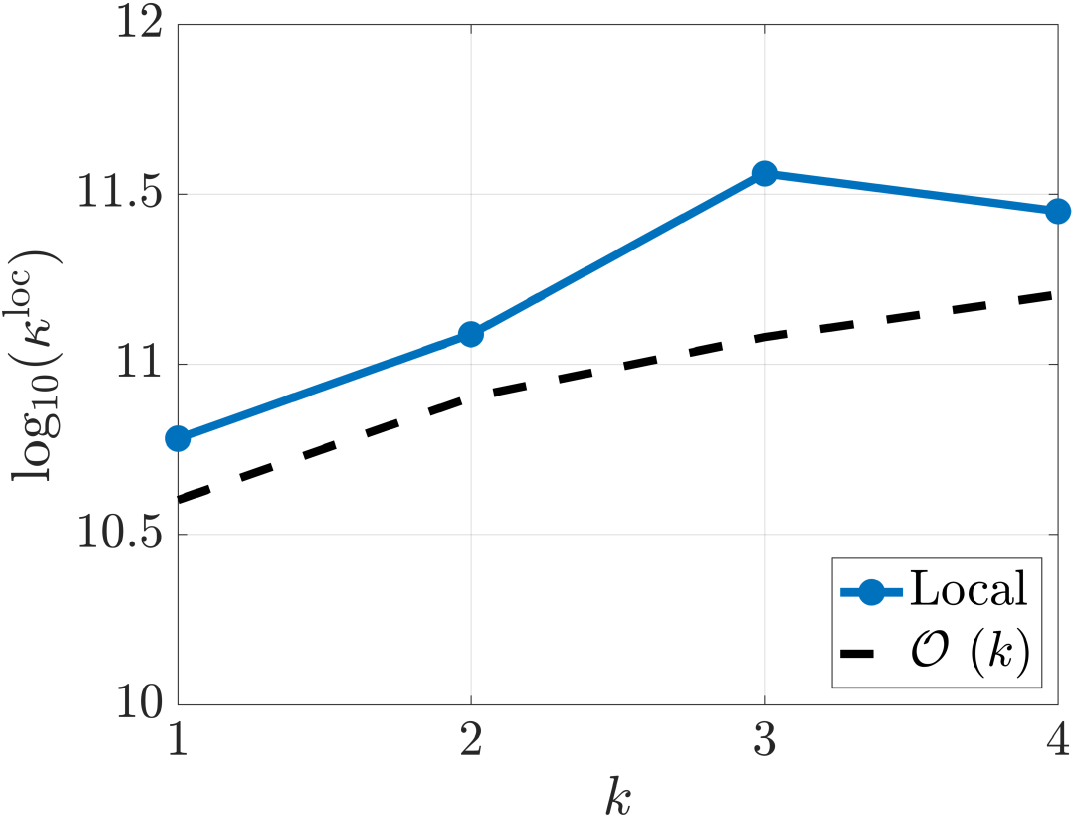}\label{fg:convCouetteKwrt_kb}}
 \caption{Taylor-Couette domain. Maximum value, over all the tested meshes, of the condition number of (a) the global problem and (b) the local problem as a function of the polynomial degree of approximation.}
  \label{fg:convCouetteKwrt_k}%
\end{figure}

First, the evolution of the  condition numbers for the HDG global and local problems is reported, for different mesh sizes and different polynomial degrees of approximation.
Figure\ \ref{fg:convCouetteKwrt_k} displays the maximum condition numbers $\kappaGlob$ and $\kappaLoc$ of the global and local problems, respectively, over the set of meshes under analysis. In both cases, a linear growth $\mathcal{O}(k)$ of the condition numbers as a function of the polynomial degree of approximation $k$ is observed.
It is worth remarking that, although the condition number of the local problem increases for high-order approximations and achieves values of approximately $10^{11}$ for $k=4$,  the size of such problems is usually fairly small (see Table\ \ref{tb:neqLocal}) and they can thus be efficiently treated using direct solvers. On the contrary, owing to the choice of Legendre basis functions for the face-based unknowns, the global condition number is maintained at lower values, also allowing for iterative solvers to be employed. 
Of course, suitable preconditioners tailored for matrices with HDG structure play a significant role in this context, see, e.g.,\ \cite{Wells-RW-18,BuiThanh-MBS-20} and should be further studied in the framework of unfitted discretizations.
\begin{table}[!h]
\caption{Dimension of the local problems using square and hexahedra elements.}
\centering
\begin{tabular}{lrrrrrrr}\hline
 Order of interpolation &   $1$ &     $2$ &     $3$ &          $4$ &          $6$ &           $8$ &            $10$ \\ \hline
 \multicolumn{8}{l}{Standard HDG element / Immersed boundary element} \\ \hline
 2D                              & $29$ &   $64$ &   $113$ &        $176$ &      $344$ &     $568$  &        $848$  \\
 3D                              & $105$ & $352$ & $833$ &    $1\ 626$ &  $4\ 460$ & $9\ 478$  &  $17\ 304$  \\ \hline\hline
 \multicolumn{8}{l}{Interface element} \\ \hline
 2D                              & $57$ &   $127$ &     $225$ &     $351$ &     $687$ &   $1\ 135$  &   $1\ 695$  \\
 3D                              & $209$ & $703$ & $1\ 665$ & $3\ 251$ & $8\ 919$ & $18\ 954$  & $34\ 607$  \\ \hline
\end{tabular}	
\label{tb:neqLocal}
\end{table}

\begin{figure}[!h]
 \centering
{\includegraphics[width=0.4\linewidth]{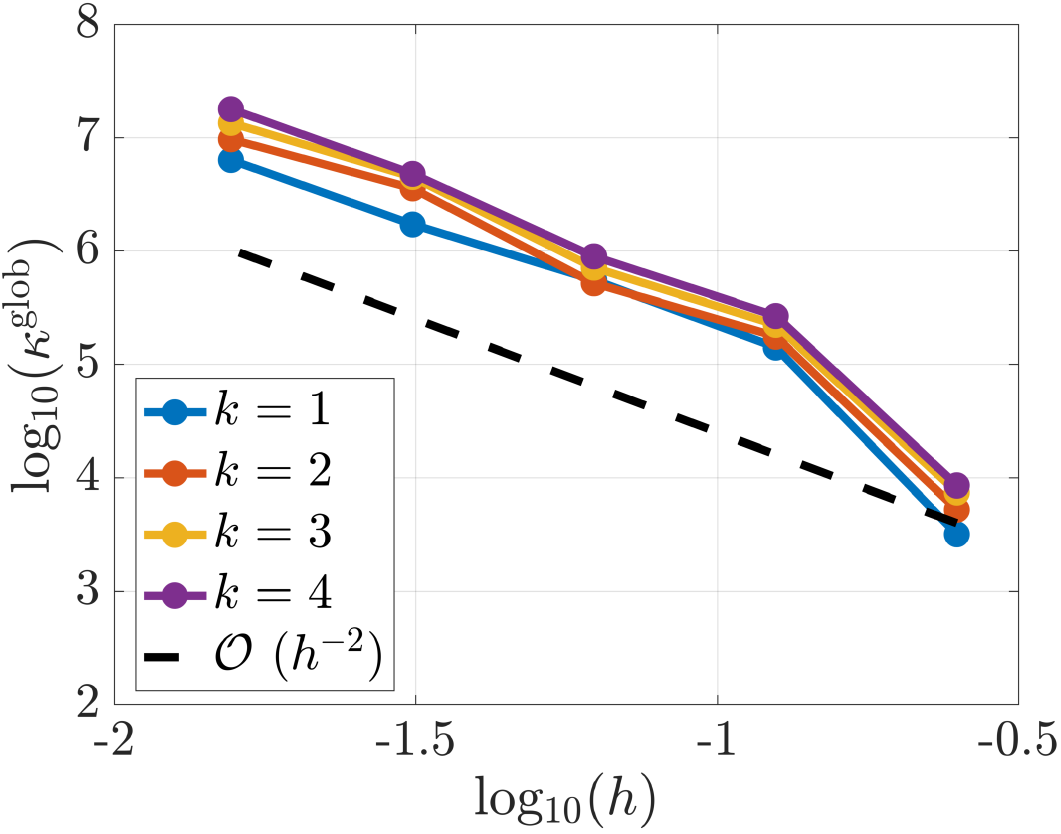}} 
 \caption{Taylor-Couette domain. Condition number of the global problem as a function of the local mesh size, for different polynomial degrees of approximation.}
  \label{fg:convCouetteKwrt_h}%
\end{figure}

The evolution of the condition number $\kappaGlob$ of the HDG global problem in terms of the local mesh size $h$ is presented in Figure\ \ref{fg:convCouetteKwrt_h}, for different polynomial degrees of approximations. The results display that the optimal rate $\mathcal{O}(h^{-2})$ for second-order elliptic differential operators is achieved  for $h\rightarrow 0$.

\begin{figure}[!h]
 \centering
 \subfloat[Pressure and mixed variable.] {\includegraphics[width=0.33\textwidth]{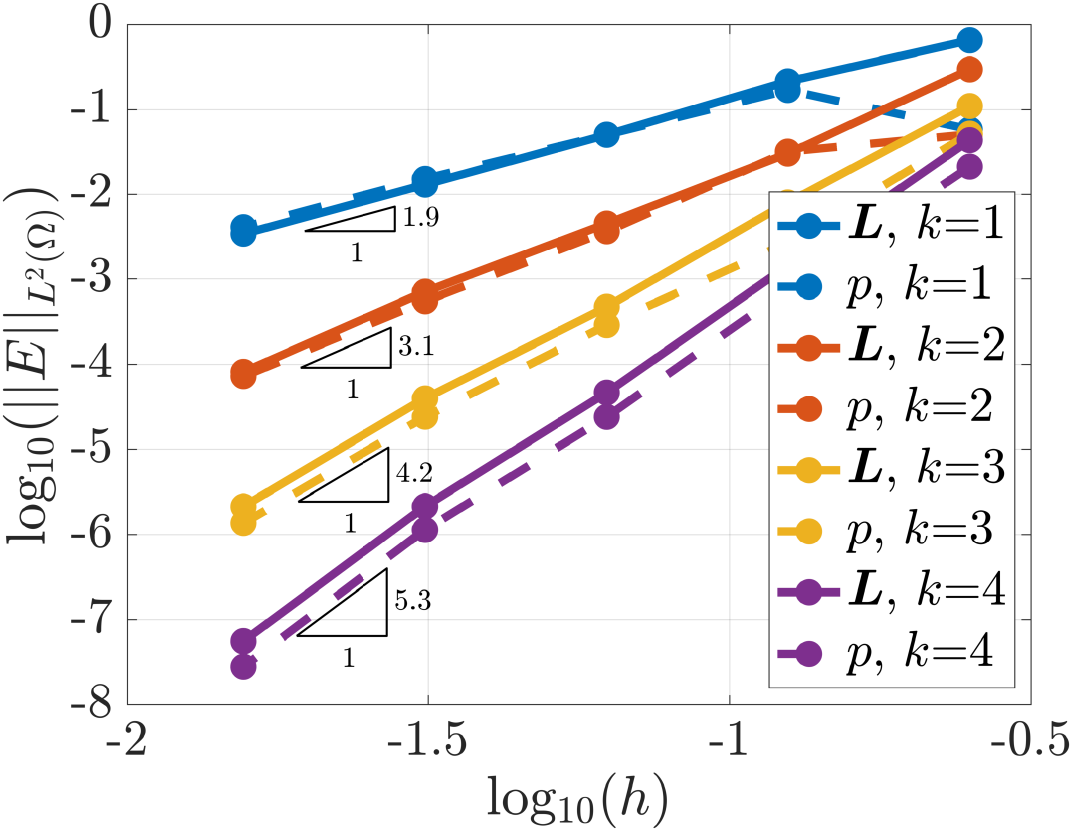}\label{fg:errCouetteConvLP}}
\subfloat[Velocity.] {\includegraphics[width=0.33\textwidth]{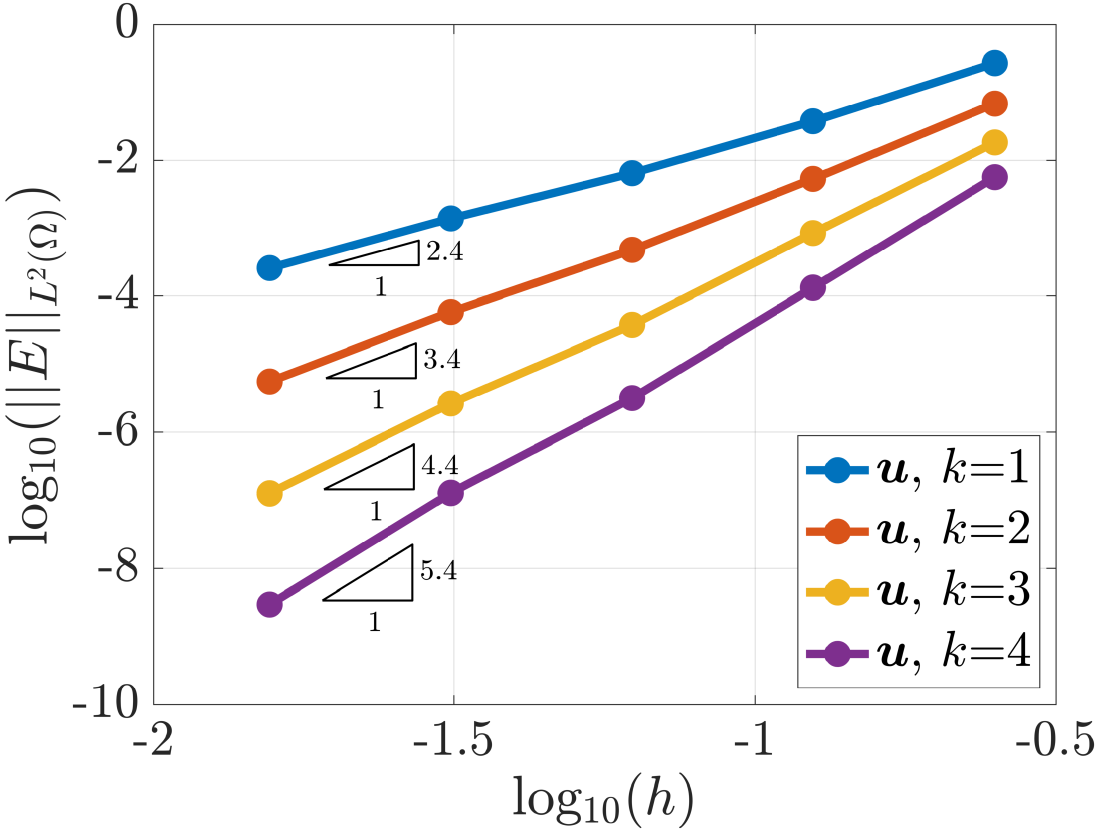}\label{fg:errCouetteConvUa}}
 \subfloat[Postprocessed velocity.]{
\includegraphics[width=0.33\textwidth]{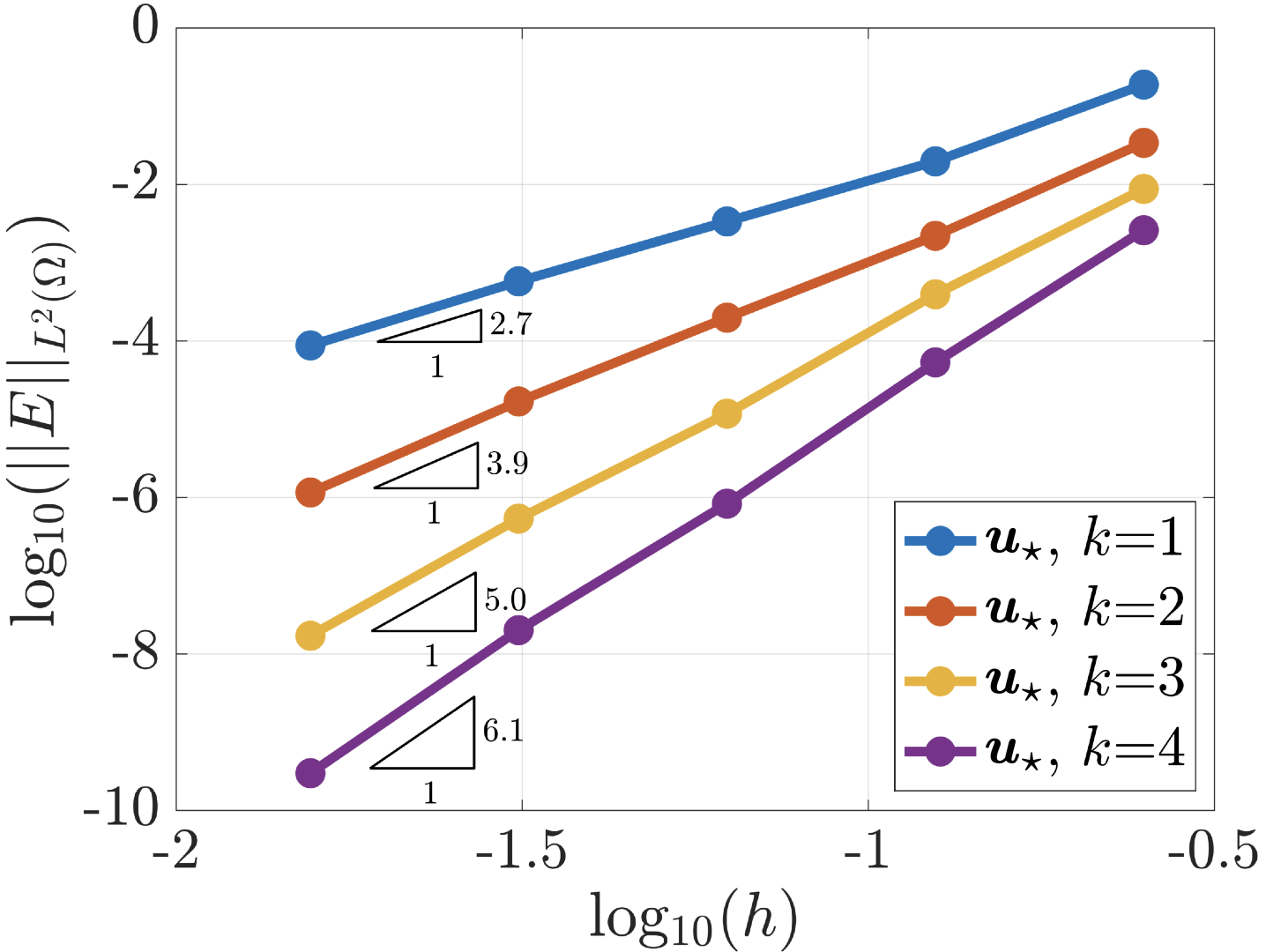}\label{fg:errCouetteConvUb}}
 \caption{Taylor-Couette domain.  Convergence of the error, measured in the $\eltwo(\Omega)$ norm, as a function of the local mesh size $h$ for (a) pressure (dashed lines) and mixed variable (continuous line), (b) velocity and (c) postprocessed velocity for different degrees of approximation $k=1,\ldots,4$.}
  \label{fg:errCouetteConv}
\end{figure}

Figure\ \ref{fg:errCouetteConv} displays the convergence of the error, measured in the $\eltwo(\Omega)$ norm,  for pressure, gradient of velocity, primal and postprocessed velocities for different mesh sizes $h$ and different polynomial degrees $k$. The results show optimal convergence of order $k+1$ for $p$, $\bL$ and $\bu$, as well as superconvergence of order $k+2$ for $\buS$.
It is worth noticing that the obtained optimal convergence rates confirm the robustness of the proposed unfitted HDG methodology coupling the exact treatment of the boundary geometry according to the NEFEM rationale with the previously discussed element extension procedure, leading to a high-order immersed boundary method, without any loss of accuracy.

\begin{figure}[!h]
 \centering
 \subfloat[Cut cell contour.]{\includegraphics[width=0.35\textwidth]{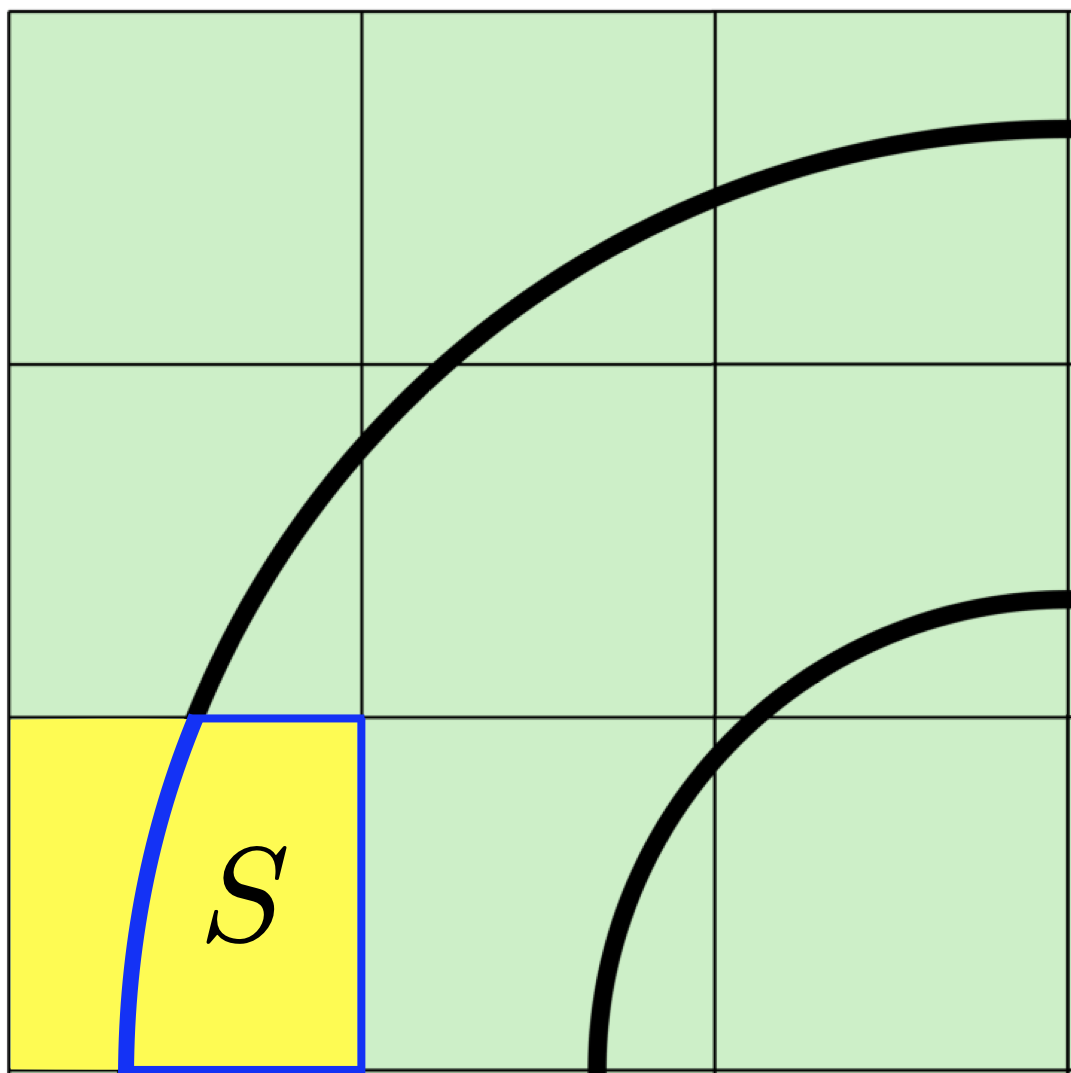}}
 \quad
 \subfloat[Extended element contour.]{\includegraphics[width=0.35\textwidth]{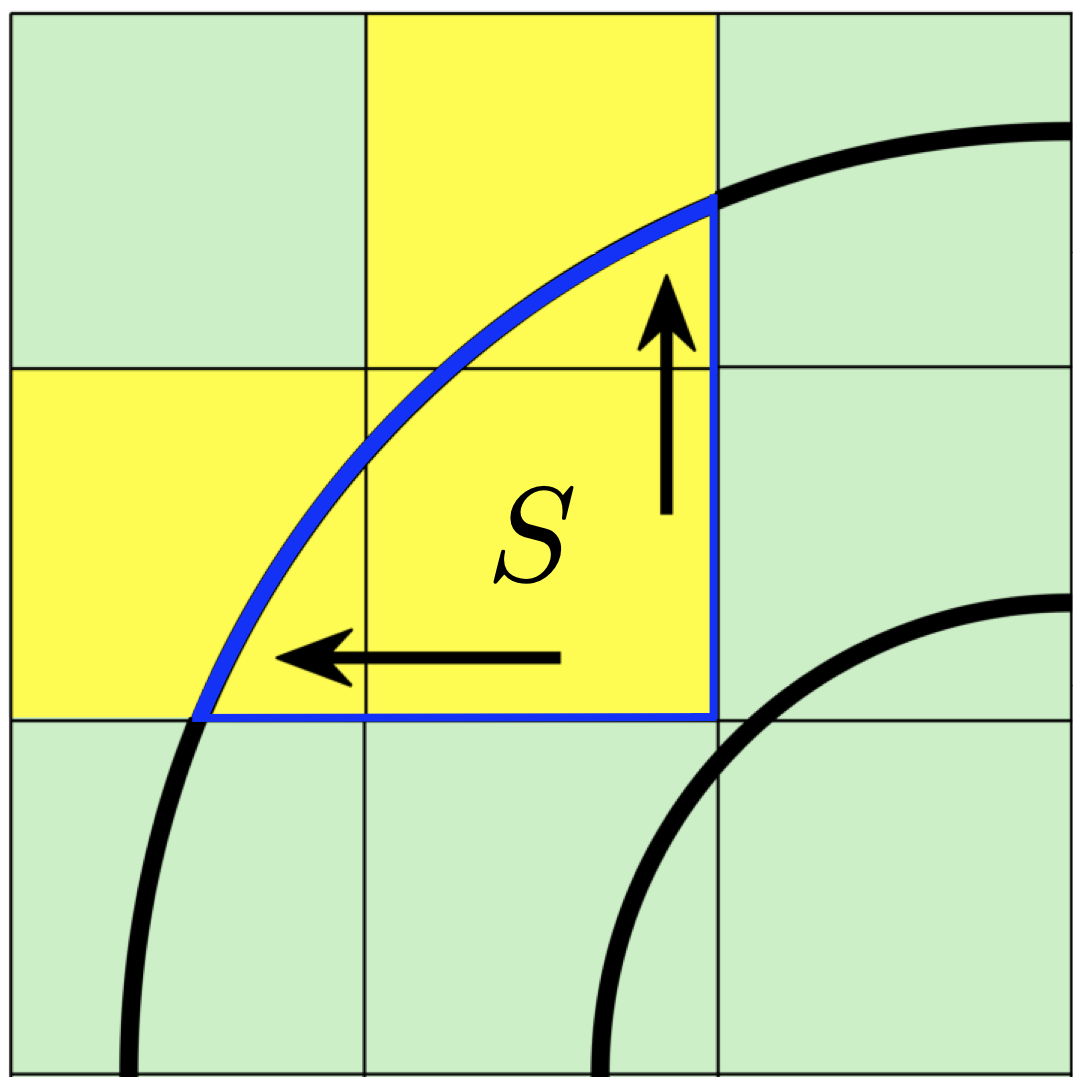}}
 \caption{Definition of the contour (in blue) to compute the mass flux entering/exiting a subset $S$ for (a) a cut cell and (b) an extended element.}
  \label{fg:fluxDef}%
\end{figure}

Finally, the error introduced by cut cells and element extension in the fulfillment of the divergence-free condition of velocity is assessed. 
Let us denote by $S \subset \Omega$ a generic subset representing a cut cell or an extended element, with $S^i = S \cap \Omega^i$ being the portion of $S$ occupied by fluid $i$. For each subset, the total mass flux entering/exiting $S$ is defined as
\begin{equation}\label{eq:elementFluxCut}
%
\fluxE = \langle \dirData \cdot \bn_e , 1 \rangle_{\Ga{D}{i}\cap \overline{S}} +  \langle \bu_e^i  \cdot \bn_e , 1 \rangle_{\Ga{N}{i}\cap \overline{S}} + \langle \bhu^i \cdot \bn_e , 1 \rangle_{\partial S^i \setminus(\Ga{D}{i}\cup\Ga{N}{i})} .
\end{equation}
For a cut cell, the subset $S$ is defined as the corresponding element $\Omega_e$ of the computational mesh,  whereas, for an extended element, $S$ is given by the union of the cells belonging to the extension patch. Figure\ \ref{fg:fluxDef} shows $S$ and the boundary employed to compute the flux $\fluxE$ in the two cases.

\begin{figure}[!h]
 \centering
 \subfloat[$k=1$.]{\includegraphics[width=0.35\textwidth]{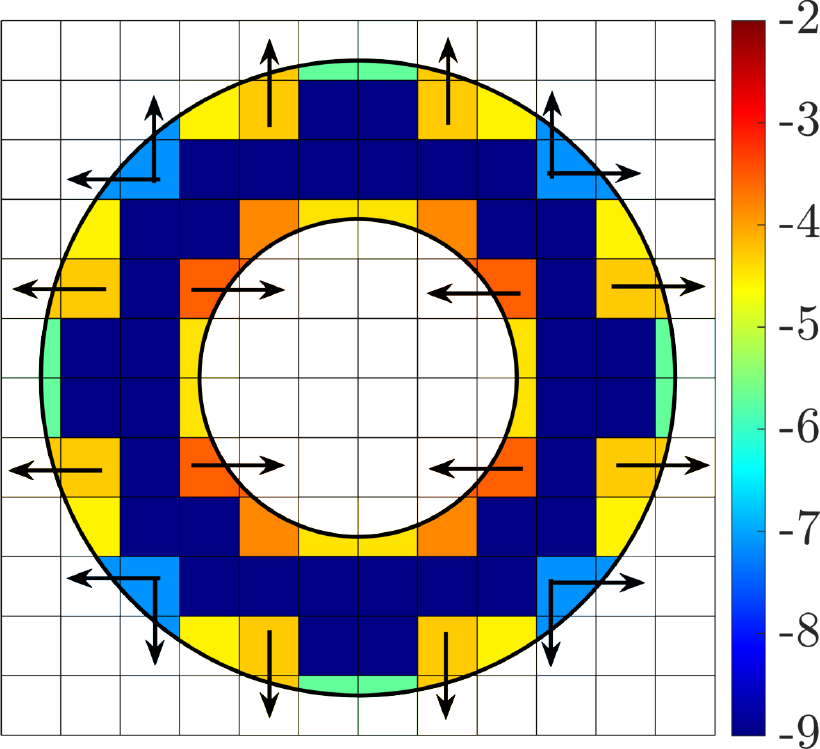}}
 \quad
 \subfloat[$k=2$.]{\includegraphics[width=0.35\textwidth]{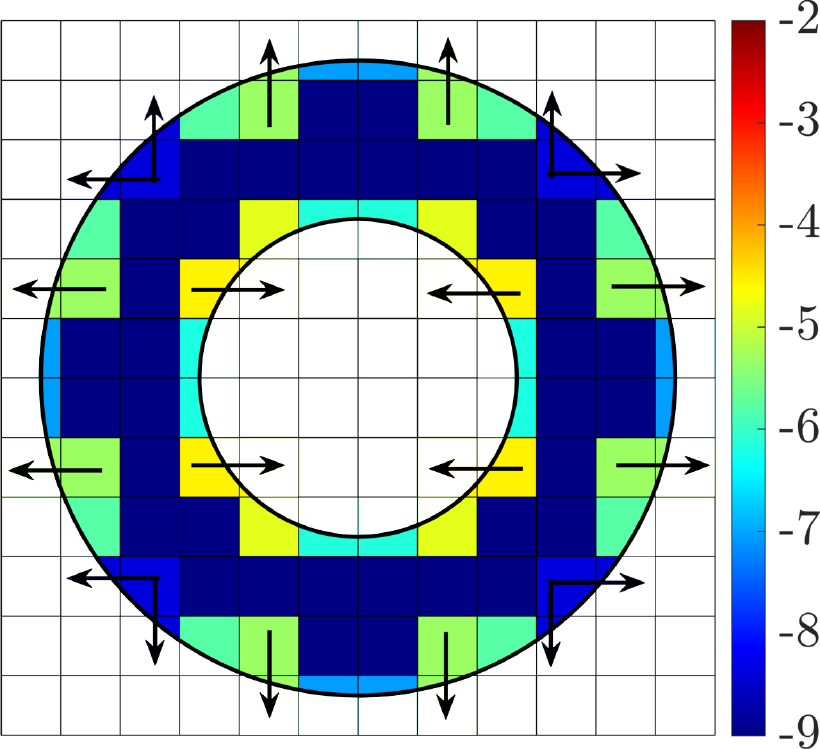}}
 
  \subfloat[$k=3$.]{\includegraphics[width=0.35\textwidth]{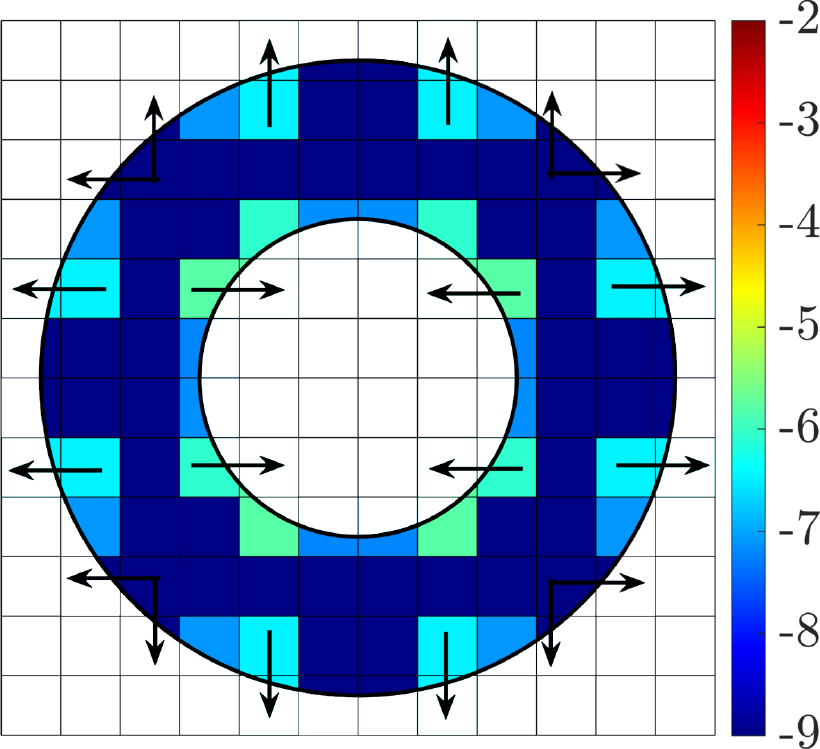}}
 \quad
 \subfloat[$k=4$.]{\includegraphics[width=0.35\textwidth]{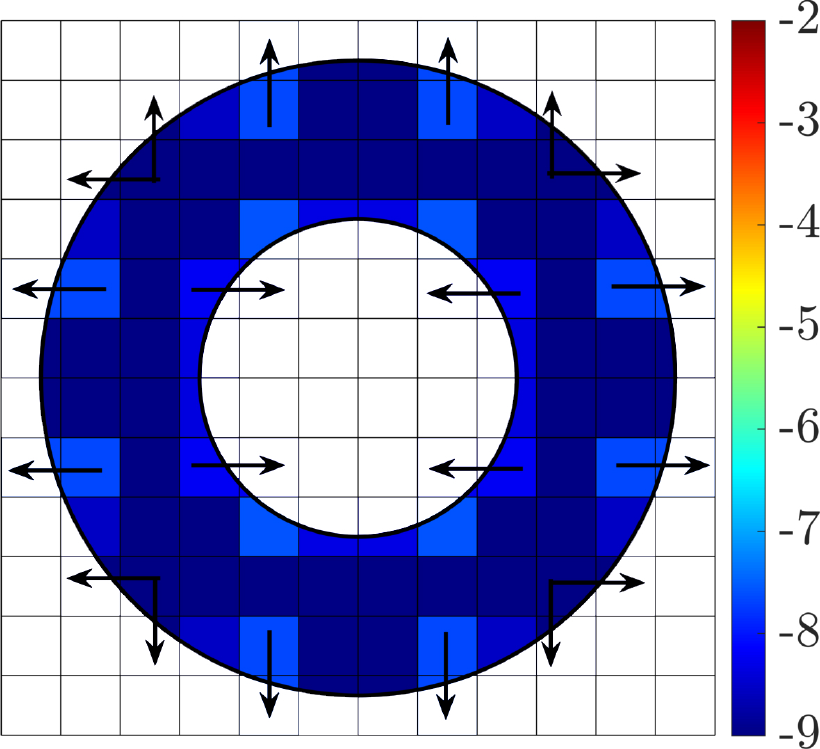}}
 \caption{Taylor-Couette domain.  Element-by-element map of $\log_{10} \abs{\fluxE}$ on the mesh composed by $16 \times 16$ elements, for different degrees of approximation $k=1,\ldots,4$.}
  \label{fg:fluxErr}%
\end{figure}

Figure\ \ref{fg:fluxErr} displays the logarithm of the absolute value of the mass flux $\fluxE$ for uncut, cut, and extended cells on the mesh composed by $16 \times 16$ elements,  for different polynomial degrees of approximation.  The choice of this specific grid stems from it being the coarsest mesh featuring all three types of elements: uncut, cut, and extended.
Cut cells and extended elements showcase a numerical error in the mass flux that decreases when the polynomial degree of approximation increases. In particular, the maximum mass flux error reported in Table\ \ref{tb:errMass} ranges from $10^{-4}$ for $k=1$ to $10^{-8}$ for $k=4$, while the corresponding best approximation errors of order $h^{k+1}$ for the mesh under analysis are always one order of magnitude above such values.
Note that the element extension procedure only affects the accuracy of mass conservation in a negligible way with respect to non-extended, cut cells.
\begin{table}[!h]
\caption{Taylor-Couette domain. Maximum flux error in cut cells and extended elements on the mesh composed by $16 \times 16$ elements.}
\centering
\begin{tabular}{lcccc}\hline
 Approximation degree $k$   &  $1$                              &   $2$                             &   $3$                              &   $4$ \\ \hline
 $h^{k+1}$  &  $3.9 \times 10^{-3}$ &   $2.4 \times 10^{-4}$ &   $1.5 \times 10^{-5}$ &   $9.5 \times 10^{-7}$  \\ \hline
 Cut cell                      &  $1.5 \times 10^{-4}$ &   $1.5 \times 10^{-5}$ &   $8.7 \times 10^{-7}$ &   $2.8 \times 10^{-8}$  \\
 Extended element    &  $2.8 \times 10^{-4}$ &   $2.8 \times 10^{-5}$ &   $1.6 \times 10^{-6}$ &   $2.2 \times 10^{-8}$   \\ \hline\hline
\end{tabular}	
\label{tb:errMass}
\end{table}

In uncut cells, mass is preserved almost exactly, with errors at most of order $10^{-9}$ (see Figure\ \ref{fg:fluxErr}). Whilst element-by-element mass conservation is fulfilled exactly in standard fitted HDG methods, the error in the hybrid variable introduced by the presence of cut cells and extended elements marginally propagates also to uncut cells. Although this prevents the error in the mass flux to achieve machine precision,  the approximation of the divergence-free condition in uncut cells attains results at least one order of magnitude more accurate than the corresponding cases featuring cut cells and extended elements.

Finally,  it is worth noticing that the total mass flux obtained by summing the contributions $\fluxE$ for all uncut, cut and extended elements in the domain achieves machine precision independently of the polynomial degree of approximation, confirming the capability of the proposed formulation to globally preserve the mass in the domain.

%
\subsection{Unfitted interface of a circular bubble at equilibrium}
\label{sc:immersedInterfaceResults}
\begin{figure}[!h]
 \centering
 \subfloat[Equilibrium bubble domain.] {\includegraphics[width=0.35\textwidth]{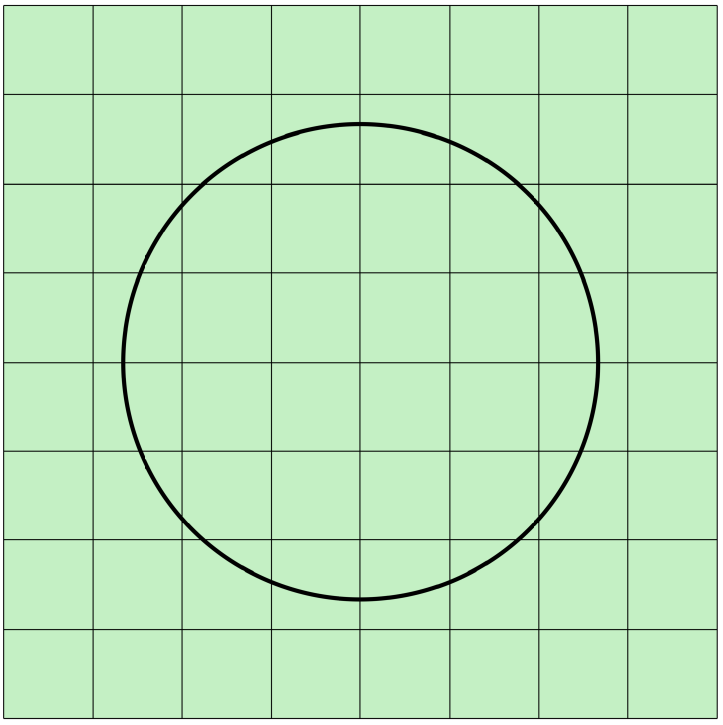}}
 \quad
 \subfloat[Extended elements.]{
\includegraphics[width=0.35\textwidth]{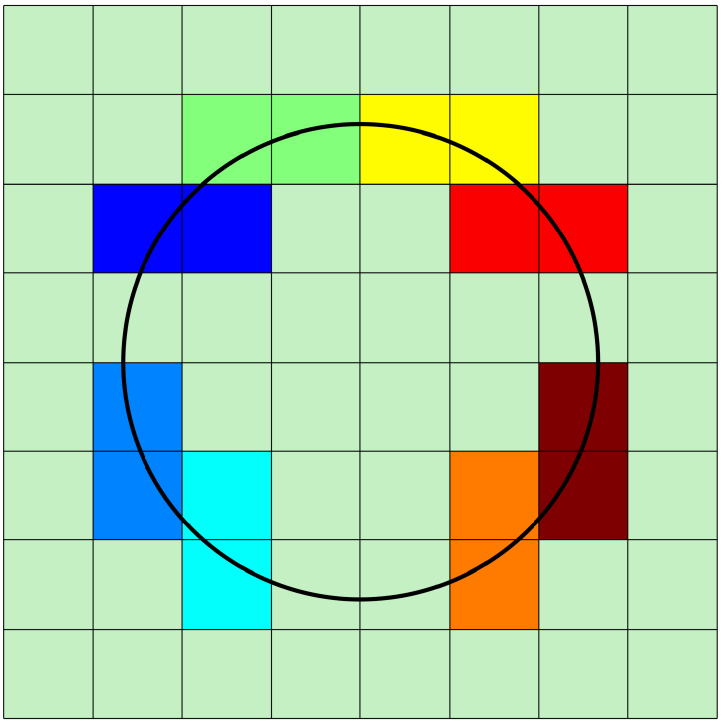}}
 \caption{Equilibrium bubble domain and extension procedure to treat badly cut cells.}
  \label{fg:bubbleDomain}%
\end{figure}

In this section, a two-fluid problem is considered, with an unfitted interface describing a circular bubble $\Omega^1$,  centered in $(x_c,y_c) {=} (0.5,0.5)$ and with radius $\R {=} 1/3$, immersed in an external fluid $\Omega^2$.  The two fluids are assumed to be immiscible, with viscosity $\mu^1 {=} 10$ and $\mu^2 {=} 1$, respectively.
The computational domain is defined as $\meshDomain {=} (0,1)^2$ and an $8 \times 8$ mesh of square elements is constructed, with an appropriate element extension procedure to handle small cut cells (Figure\ \ref{fg:bubbleDomain}).

\begin{figure}[!h]
 \centering
 \subfloat[$\log_{10}(\|\bu - \bu^{\text{ref}} \|)$.] {\includegraphics[width=0.4\linewidth]{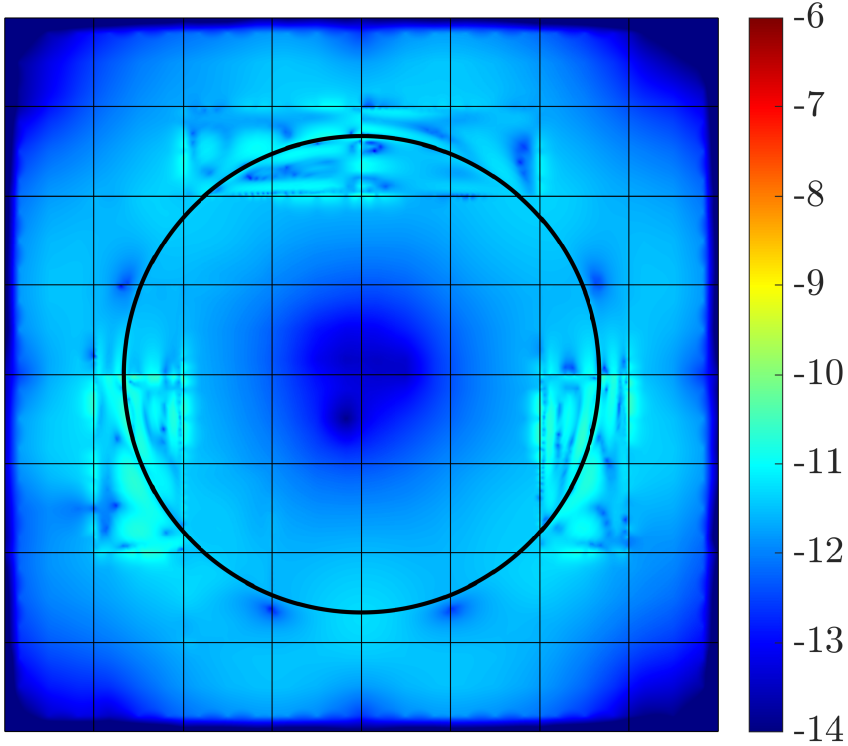}}
 \quad
 \subfloat[$\log_{10}(|p - p^{\text{ref}} |)$.]{
\includegraphics[width=0.4\linewidth]{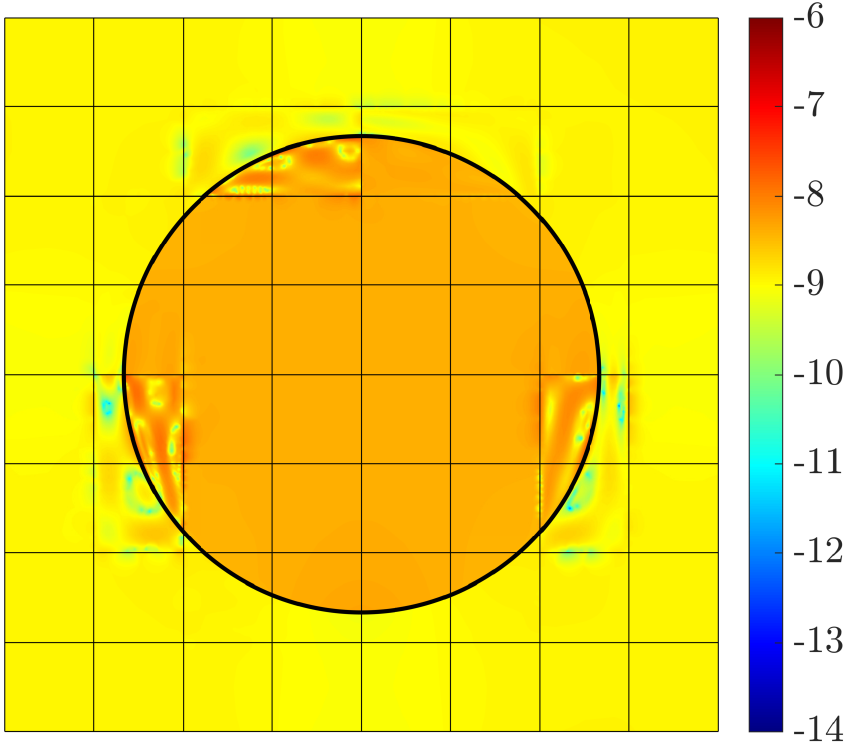}}
 \caption{Equilibrium bubble domain. Logarithm of the point-wise absolute error of (a) the module of the velocity and (b) the pressure computed with polynomial functions of degree $k=1$.}
  \label{fg:bubbleErr}
\end{figure}

Under the assumption of zero body forces (i.e.,  $\bm{s} {=} \bm{0}$) and no-slip boundary conditions $ \dirData {=} \bm{0}$ on $\partial\meshDomain$, the problem admits an analytical solution 
\begin{equation*}
\bu^{\text{ref}} = \bm{0} \, ,
\qquad
 p^{\text{ref}} = \begin{cases}
     \pi \R \gamma - \displaystyle \frac{\gamma}{\R},  & \text{if } r < \R \, ,\\
  \pi \R  \gamma, & \text{elsewhere} \, ,
\end{cases}
\end{equation*}
with $r {=} ((x-x_c)^2+(y-y_c)^2)^{1/2}$ and $\gamma {=} 1$.
It is worth noticing that the target velocity solution is constant in $\Omega^1 \cup \Omega^2$, whereas the target pressure is a discontinuous function, piecewise constant in $\Omega^1$ and $\Omega^2$.
Hence, the numerical results obtained using an HDG approximation of degree $k=1$ yield a discretization error achieving machine precision for all variables, with errors measured in the $\eltwo(\meshDomain)$ norm of order $10^{-12}$ for velocity,  $10^{-9}$ for pressure and $10^{-10}$ for the gradient of velocity. 
Figure\ \ref{fg:bubbleErr} reports a point-wise error map of the logarithm of the absolute error for the module of the velocity and the pressure, confirming the ability of the method to accurately approximate the solution in the two fluid domains, while guaranteeing the equilibrium via the imposition of the transmission condition across the interface, geometrically described using a NURBS.

%
\subsection{Microfluidic channel with complex unfitted obstacles}
\label{sc:microMix}
\begin{figure}[t]
 \centering
 \subfloat[Problem setup.\label{fg:channelGeo}]{\quad\includegraphics[width=0.245\textwidth]{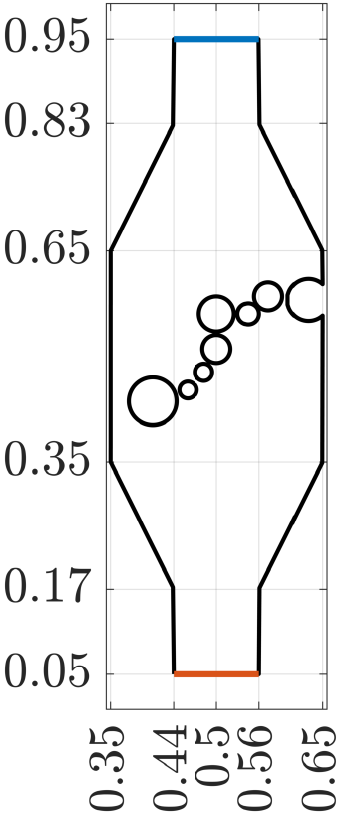}\quad} 
\subfloat[Microchannel domain.\label{fg:channelMesh}] {\qquad\includegraphics[width=0.2\textwidth]{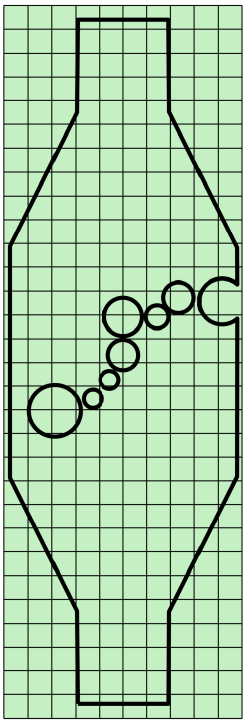}\qquad} 
 \subfloat[Extended elements.\label{fg:channelEE}]{\quad\includegraphics[width=0.2\textwidth]{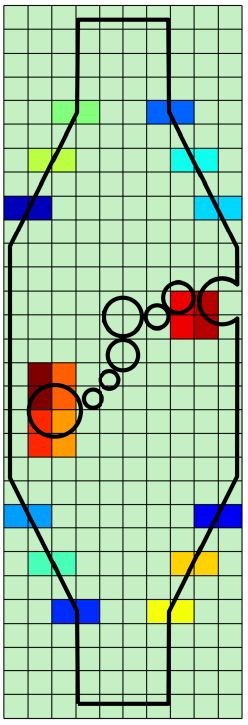}\quad}
 \caption{Microchannel domain.  (a) Problem setup with physical walls (black), inlet (blue) and outlet (red) boundaries.  (b) Detail of the computational domain in the region $(0.34375,0.65625)\times(0,1)$. (c) Element extension to avoid badly cut cells.}
  \label{fg:channelDomain}%
\end{figure}

A microchannel inspired by the realistic design of a passive microfluidic mixer\ \cite{Wang:19} is studied here to showcase the suitability of the unfitted HDG-NEFEM method to perform high-fidelity simulations of microfluidic flows in complex geometric systems described by NURBS.
The microchannel $\Omega$, displayed in Figure\ \ref{fg:channelGeo},  consists of a divergent-convergent structure with a series of eight circular obstacles located in the central portion of the domain.  The geometric details of each inclusion (i.e., coordinates $(x_c^i,y_c^i)$ of the centers and radii $\R^i$ of the circles) are reported in Table\ \ref{tb:mixerData}. 
%
A fluid with viscosity $\mu =1$ enters from the top boundary (in blue), where the Dirichlet datum $\dirData {=} (0,-1)$ is imposed,  and exits from the bottom contour (in red), where a homogeneous Neumann condition is applied. No-slip conditions are enforced on the remaining boundaries and on the obstacles to model physical walls. Moreover, body forces are set to zero. 
\begin{table}[!h]
\centering
\caption{Geometric data (coordinates $(x_c^i,y_c^i)$ of the centers and radii $\R^i$) of the eight circular obstacles,  ordered from left to right. }
\resizebox{\textwidth}{!}{%
\begin{tabular}{lccccccccc}
\hline
Obstacle $i$ & & $1$ & $2$ & $3$ & $4$ & $5$ & $6$ & $7$ & $8$ \\ \hline\hline
Center coordinate $x_c^i$ & [$\times 10^{-1}$] & $4.10$ & $4.60$ & $4.82$ & $5.00$ & $5.00$ & $5.45$ & $5.73$ & $6.30$ \\
Center  coordinate $y_c^i$ & [$\times 10^{-1}$] & $4.38$ & $4.52$ & $4.77$ & $5.10$ & $5.60$ & $5.60$ & $5.85$ & $5.80$ \\  \hline
Radius $\R^i$ & [$\times 10^{-2}$]  & $3.4$ & $1.2$ & $1.2$ & $2.0$ & $2.5$ & $1.5$ & $2.0$ & $3.0$ \\
\hline
\end{tabular}
}%
\label{tb:mixerData}
\end{table}

As classic in unfitted methods, the computational domain needs to be selected to guarantee that the boundaries of the physical domain $\Omega$ are contained in $\overline{\meshDomain}$. In this case, no \emph{a priori} knowledge on the shape of $\Omega$ is employed and a square computational domain $\meshDomain {=} (0,1)^2$ is selected to ensure that both the inlet and outlet surfaces of the microchannel lie in the interior of $\meshDomain$. The computational domain is thus subdivided into $32 \times 32$ square elements. 
%
%
It is worth noticing that given the rectangular shape of the physical domain $\Omega$ under analysis, approximately two thirds of the mesh elements are inactive lying in $\meshDomain \setminus \Omega$. 
Of course, a smaller computational domain tailored for the problem under analysis could be devised. Indeed, Figure\ \ref{fg:channelMesh} displays the portion of $\meshDomain$ comprised in the region $(0.34375,0.65625)\times(0,1)$, where $\Omega$ is located, which features only $10\times 32$ square elements. Nonetheless, the computational complexity of the proposed method is independent of the number of inactive cells since these are identified once, during the preprocessing stage, and subsequently omitted from the routines to construct the HDG local and global problems, which are only assembled for the active cells.
%
%
In addition, the extended elements employed to avoid ill conditioning issues due to badly cut cells are reported in Figure\ \ref{fg:channelEE}.

\begin{figure}[!h]
 \centering
  \subfloat[Pressure.]{
\includegraphics[height=0.53\textwidth]{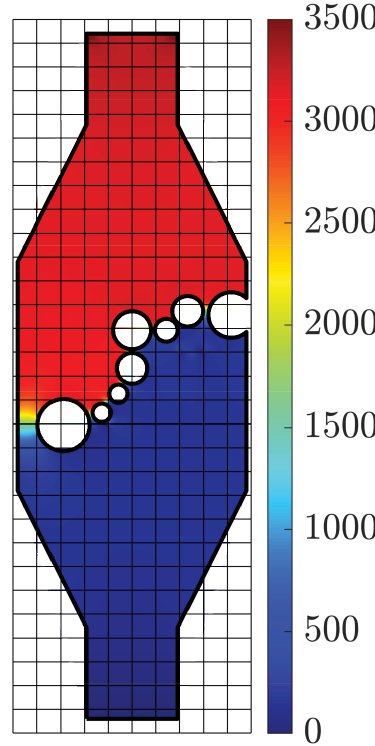}\label{fg:channelSolP}} 
 \;
 \subfloat[Module of velocity.] {\includegraphics[height=0.53\textwidth]{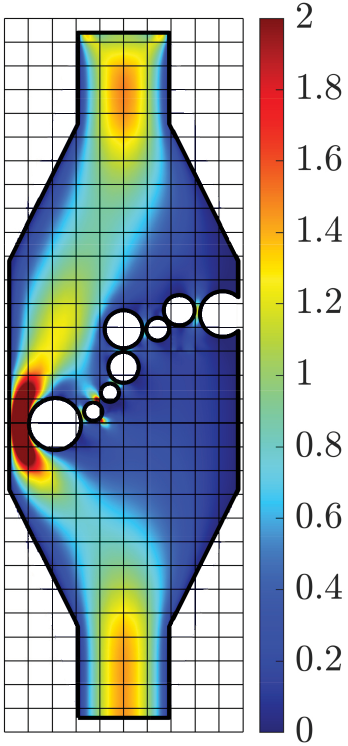}\label{fg:channelSolVel}}
 \;
 \subfloat[Streamlines.]{\includegraphics[height=0.52\textwidth]{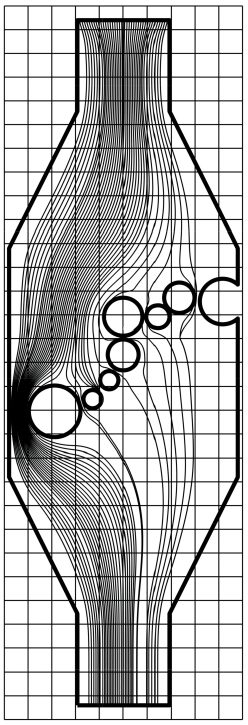}\label{fg:channelSolStream}}
\; 
 \subfloat[Adapted degree.]{
\includegraphics[height=0.53\textwidth]{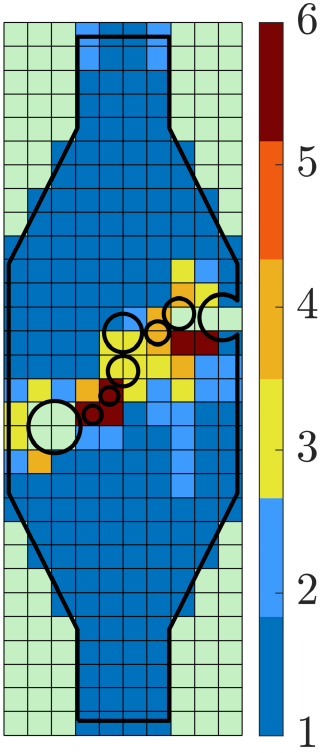}\label{fg:channelSolAdapt}}
 \caption{Microchannel domain.  (a) Pressure. (b) Module of velocity. (c) Streamlines. (d) Map of the adapted polynomial degree of approximation.}
  \label{fg:channelSolution}
\end{figure}

Figure\ \ref{fg:channelSolution} shows pressure and velocity fields computed using the unfitted HDG-NEFEM method.  While both fields appear smooth, the presence of the obstacles is responsible for a significant deviation of the flow, with a high pressure gradient (Figure\ \ref{fg:channelSolP}) and an acceleration of the flow (Figure\ \ref{fg:channelSolVel}) in the vicinity of the left-most obstacle.  In addition, the streamlines in Figure\ \ref{fg:channelSolStream} clearly highlight an additional difficulty of this problem, with localized flow features appearing when two obstacles are not perfectly in contact.

\begin{figure}[!h]
 \centering
 \subfloat[Adapted degree.]{
\includegraphics[width=0.4\textwidth]{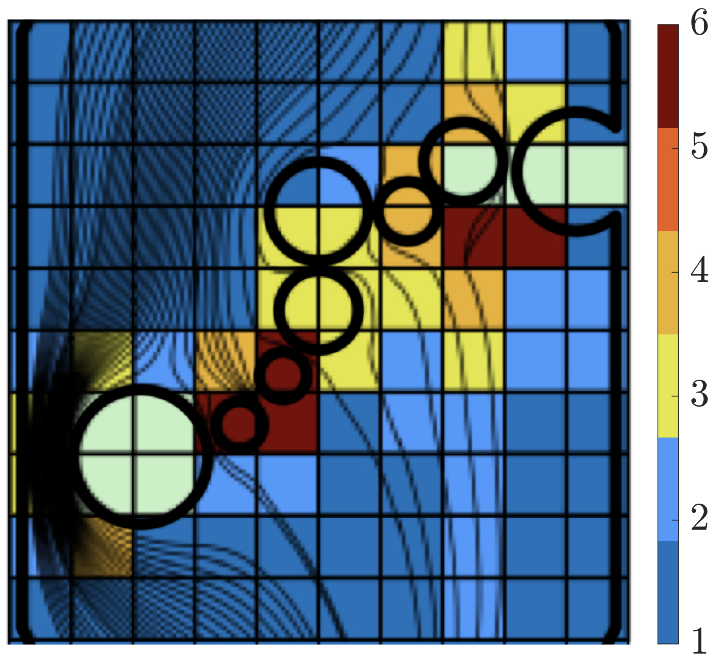}\label{fg:channelAdaptZoom}}
\quad 
 \subfloat[Extended elements.]{\includegraphics[width=0.36\textwidth]{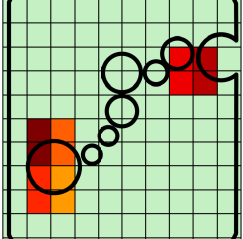}\label{fg:channelEEzoom}}
 \caption{Microchannel domain.  Detail of (a) the adapted degree superimposed with the streamlines and (b) the extended elements in the vicinity of the obstacles.}
  \label{fg:ChannelZoom}
\end{figure}

To accurately describe the flow variations observed in\ \cite{Wang:19} in the region featuring the obstacles,  the degree adaptive procedure described in Section\ \ref{sec:degAdapt} is employed. Exploiting the superconvergence properties of HDG approximations in elliptic problems,  the map of non-uniform polynomial degrees in Figure\ \ref{fg:channelSolAdapt} is devised to ensure an accuracy of two significant digits in the velocity field.
More precisely, the methodology automatically selects the most appropriate polynomial degree to locally approximate the flow, achieving degree $6$ in the regions where very small features appear,  while maintaining degree $1$ and $2$ where the flow is less complex. The map of the adapted degree in the vicinity of the obstacles,  together with the extended elements employed for computation are reported in Figures\ \ref{fg:channelAdaptZoom} and\ \ref{fg:channelEEzoom}, respectively.

This confirms the ability of the method to provide highly accurate results, going high-order with the polynomial approximations, while limiting the increase of the computational cost. In addition, the adaptive procedure is capable of identifying the regions where physically relevant phenomena take place and automatically devising the precision required by the discretization. Finally, it is worth remarking that the combined use of unfitted meshes and exact description of the geometry via NURBS allows to achieve these results on relatively coarse meshes,  whose design is not constrained by the geometric features of the domain. Moreover, the presented framework is robust also in the presence of elements cut in multiple subregions.

%
\subsection{Multiple unfitted interfaces in an emulsion of two fluids}
\label{sc:emulsion}
The last example features a two-fluid system describing an emulsion in a porous medium, with a downward flow due to the gravitational force as presented in\ \cite{Avraam:95,LuRiv:15,Janetti:17}.
More precisely,  a dispersed phase referred to as fluid $1$ (e.g., oil droplets) is immersed in a continuous phase of the second fluid (e.g., water) that permeates the pores of a homogeneous porous medium.

\begin{figure}[!h]
 \centering
 \subfloat[Problem setup. \label{fg:dropPb}]{\includegraphics[width=0.3\textwidth]{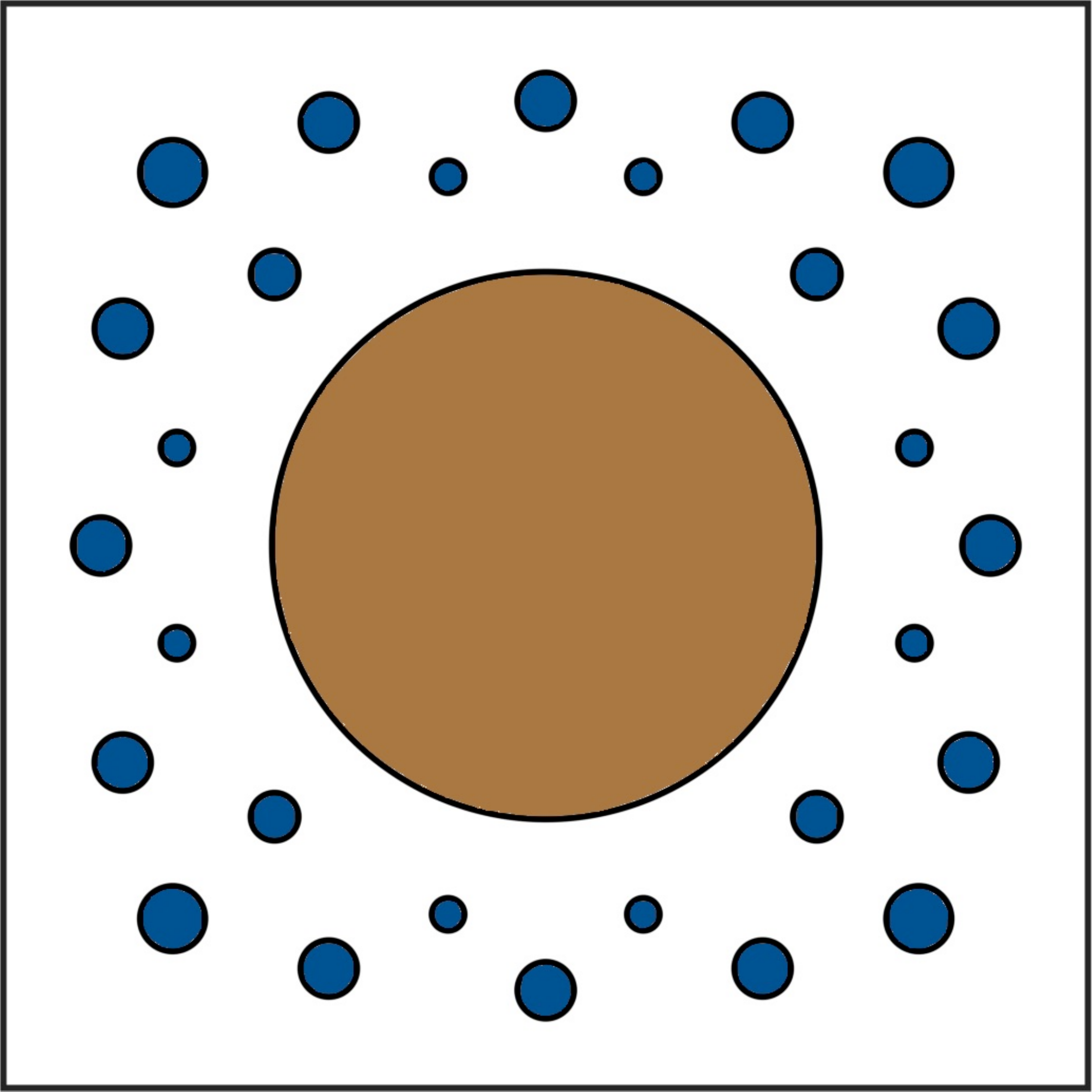}}
 \quad
 \subfloat[Emulsion domain.\label{fg:dropMesh}] {\includegraphics[width=0.3\textwidth]{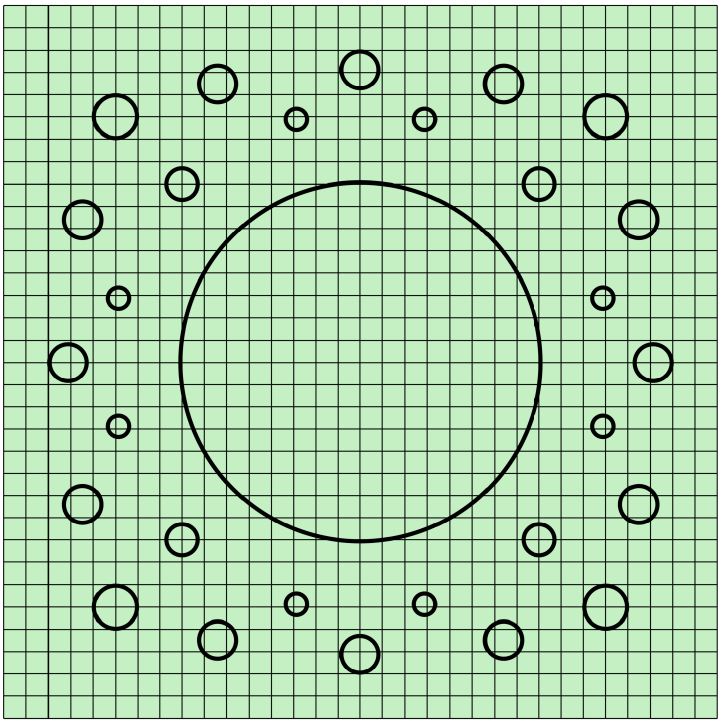}}
 \caption{Emulsion domain. (a) Problem setup with fluid $1$ (blue), fluid $2$ (white) and solid pore (brown). (b) Computational domain. }
  \label{fg:dropDomain}%
\end{figure}

Under the assumption of periodic porous medium, a microscopic unit cell $\meshDomain {=} (0,1)^2$, illustrated in Figure\ \ref{fg:dropPb}, is considered as computational domain,  with a mesh of $32 \times 32$ square elements (Figure\ \ref{fg:dropMesh}). 
Upon a non-dimensional analysis of the physical data reported in\ \cite{Avraam:95,LuRiv:15,Janetti:17},  the following setup is employed.
The solid matrix of the porous medium is defined as the brown circular region centered in $(x_c,y_c) {=} (0.5,0.5)$ with radius $\Rpore {=} 0.25231$. This is assumed to be rigid and impermeable, leading to a material with porosity of approximately $0.8$.
%
The domain $\Omega^1$ of fluid $1$ is subdivided into $28$ circular droplets,  each defined by closed NURBS curves, whose geometric data are detailed in Table\ \ref{tb:dropData}.
Fluid $2$ is defined in the domain $\Omega^2 := \meshDomain \setminus \Omega^1$.
A viscosity ratio $\mu^1/\mu^2 {=} 10$ is considered, with the viscosities of the two fluids being $\mu^1 {=} 40$ and $\mu^2 {=} 4$, respectively.
To model the effect of gravity,  the body force is set to $\bm{s} {=} (0,-613.125)$, whereas the surface tension coefficient $\gamma {=} 2.4\times 10^5$ is employed to enforce the interface conditions between the two fluids.
Finally, no-slip boundary conditions are imposed on the contour of the central pore $\{ (x,y) \in \meshDomain \ : \  (x-x_c)^2 + (y-y_c)^2 = \Rpore^2 \}$ and periodic conditions are prescribed on the top and bottom, as well as left and right boundaries of $\partial\meshDomain$.
%
%
\begin{table}[!h]
\centering
\caption{Geometric data (center coordinates $(x_c^i,y_c^i)$ and radius $\R^i$) of the $28$ circular droplets.}
\resizebox{\textwidth}{!}{%
\begin{tabular}{lccccccccccc}
\hline
Droplet $i$ & & $1$ & $2$ & $3$ & $4$ & $5$ & $6$ & $7$ & $8$ & $9$ & $10$ \\ \hline
Coord. $x_c^i$ & [$\times 10^{-1}$] & $1.56$ & $1.56$ & $8.44$ & $8.44$ & $5.0$ & $1.0$ & $5.0$ & $9.0$ & $7.0$ & $3.0$ \\
Coord. $y_c^i$ & [$\times 10^{-1}$] & $1.56$ & $8.44$ & $8.44$ & $1.56$ & $1.0$ & $5.0$ & $9.0$ & $5.0$ & $1.1$ & $1.1$ \\
\hline
Radius $\R^i$ & [$\times 10^{-2}$] & $3.0$ & $3.0$ & $3.0$ & $3.0$ & $2.6$ & $2.6$ & $2.6$ & $2.6$  & $2.6$ & $2.6$ \\
\hline\hline
Droplet $i$ &  & $11$ & $12$ & $13$ & $14$ & $15$ & $16$ & $17$ & $18$ & $19$ & $20$ \\ \hline
Coord. $x_c^i$ & [$\times 10^{-1}$] & $1.1$ & $1.1$ & $3.0$ & $7.0$ & $8.9$ & $8.9$ & $7.5$ & $2.5$ & $2.5$ & $7.5$ \\
Coord. $y_c^i$ & [$\times 10^{-1}$] & $3.0$ & $7.0$ & $8.9$ & $8.9$ & $7.0$ & $3.0$ & $2.5$ & $2.5$ & $7.5$ & $7.5$ \\
\hline
Radius $\R^i$ & [$\times 10^{-2}$] & $2.6$ & $2.6$ & $2.6$ & $2.6$ & $2.6$ & $2.6$ & $2.2$ & $2.2$ & $2.2$ & $2.2$ \\
\hline\hline
Droplet $i$ & & $21$ & $22$ & $23$ & $24$ & $25$ & $26$ & $27$ & $28$ & & \\ \hline
Coord. $x_c^i$ & [$\times 10^{-1}$] & $5.9$ & $4.1$ & $1.6$ & $1.6$ & $4.1$ & $5.9$ & $8.4$ & $8.4$ & & \\
Coord. $y_c^i$ & [$\times 10^{-1}$] & $1.6$ & $1.6$ & $4.1$ & $5.9$ & $8.4$ & $8.4$ & $4.1$ & $5.9$ & & \\
\hline
Radius $\R^i$ & [$\times 10^{-2}$] & $1.5$ & $1.5$ & $1.5$ & $1.5$  & $1.5$ & $1.5$ & $1.5$ & $1.5$ & & \\
\hline
\end{tabular}
}%
\label{tb:dropData}
\end{table}

\begin{remark}[Periodic boundary conditions in HDG]
Let $\Ga{L}{}$ and $\Ga{R}{}$ respecively denote the left and right vertical boundaries where periodic conditions are prescribed. 
Assuming that the vertical coordinate of the nodes along $\Ga{L}{}$ and $\Ga{R}{}$ coincides,  the imposition of periodic conditions in HDG methods is straightforward.
By eliminating the hybrid unknown defined along $\Ga{R}{}$,  the faces along $\Ga{L}{}$ can be considered as internal faces, shared by the left-most and right-most elements.
Thus,  the matching of the solution along $\Ga{L}{}$ and $\Ga{R}{}$ follows automatically from the uniqueness of the hybrid variable along $\Ga{L}{}$.
The equality of the fluxes across $\Ga{L}{}$ and $\Ga{R}{}$ is naturally enforced by means of the HDG global problem since the faces along $\Ga{L}{}$ are now considered as internal faces.
Of course, the above considerations can be extended to any pair of appropriately oriented boundaries $\Ga{I}{}$ and $\Ga{J}{}$.
\end{remark}

\begin{figure}[!h]
 \centering
 \subfloat[Module of velocity.] {\includegraphics[height=0.3\textwidth]{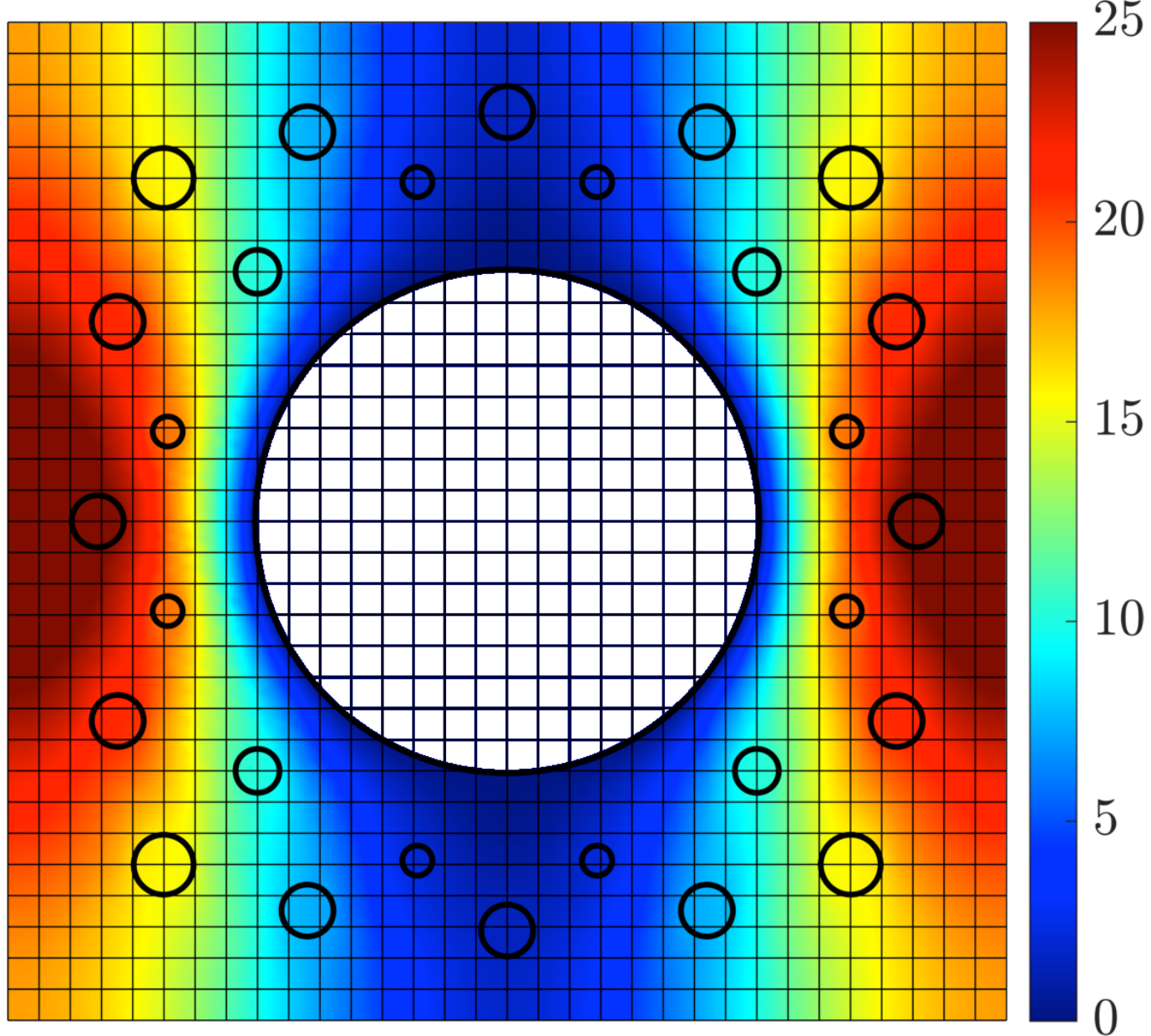}}
 \quad
 \subfloat[Streamlines.]{\includegraphics[height=0.3\textwidth]{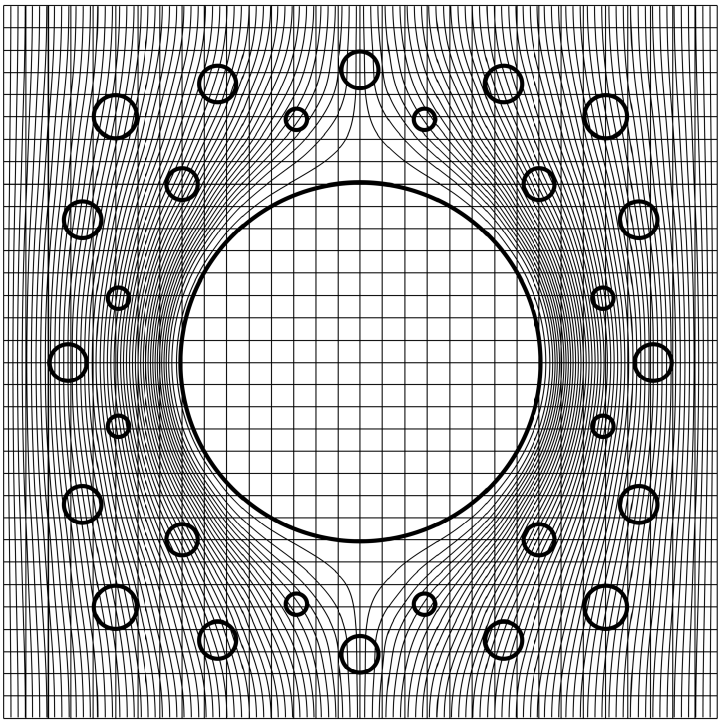}}

 \subfloat[Pressure. \label{fg:dropP}]{\quad\quad\includegraphics[height=0.3\textwidth]{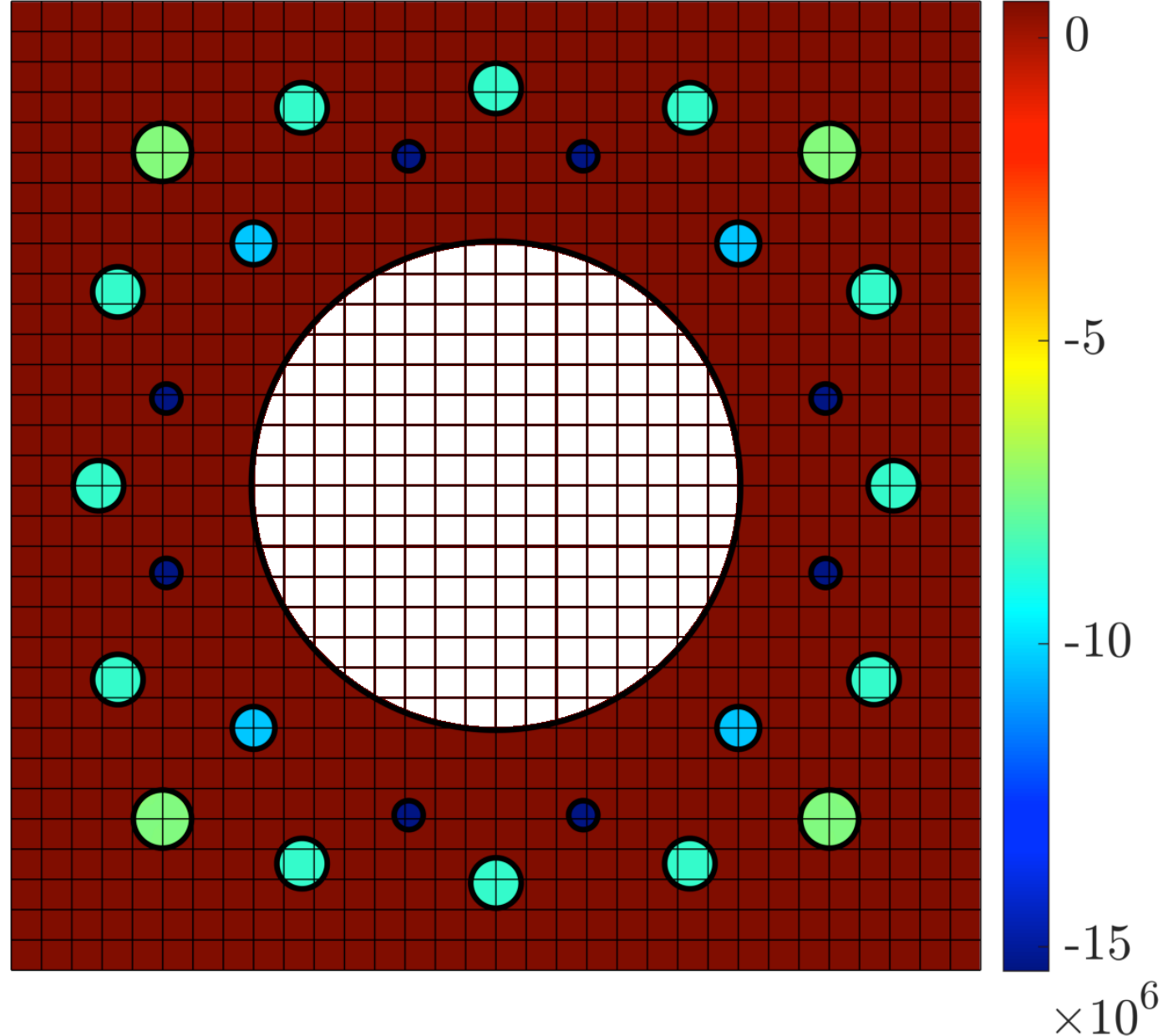}}
 \caption{Emulsion domain.  (a) Module of velocity. (b) Streamlines. (c) Pressure. }
  \label{fg:dropSol}
\end{figure}

\begin{figure}[!h]
 \centering
 \quad
  \subfloat[Extended elements.\label{fg:dropEE}]{\includegraphics[height=0.3\textwidth]{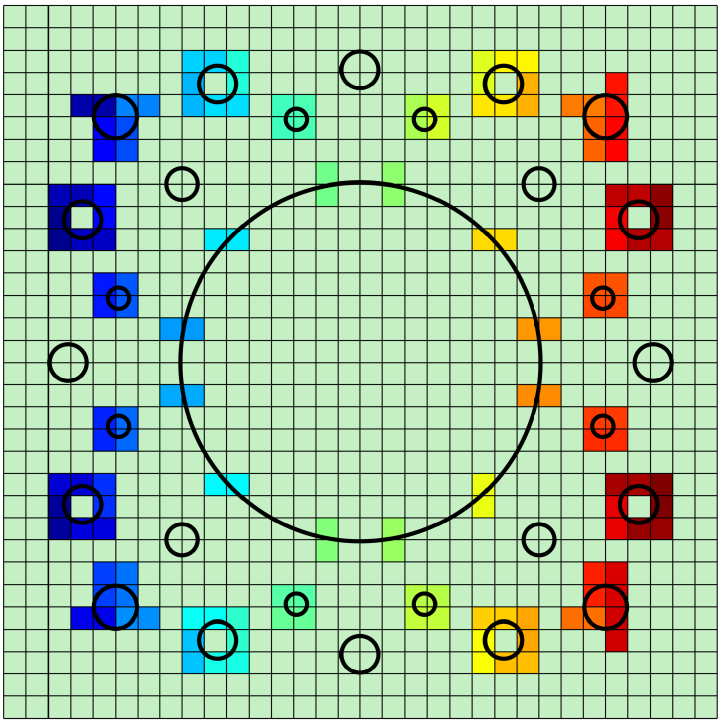}}
 \quad
  \subfloat[Adapted degree. \label{fg:dropKadaptDeg}]{\includegraphics[height=0.31\textwidth]{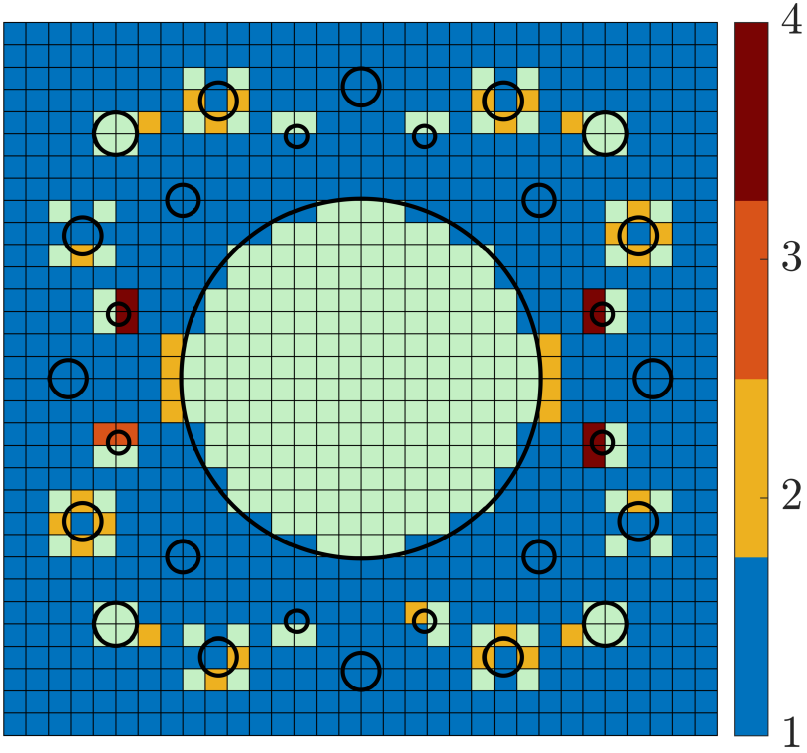}}
 \caption{Emulsion domain.  (a) Element extension to avoid badly cut cells. (b) Map of the adapted polynomial degree of approximation.}
  \label{fg:dropKadapt}%
\end{figure}

Figure\ \ref{fg:dropSol} displays the pressure and velocity fields computed using the unfitted HDG approach, with a non-uniform polynomial degree approximation.
The results show localized pressure jumps at the interfaces between fluid $1$ and fluid $2$. These are due to the effect of the surface tension, which depends upon the curvature of the droplets in the dispersed phase.  More precisely, the pressure jump increases proportionally to the inverse of the radius of the droplets (Figure\ \ref{fg:dropP}).
%
As commented for the previous cases, the element extension procedure is performed to avoid ill conditioning issues due to badly cut cells (Figure\ \ref{fg:dropEE}).  
Moreover, the degree-adaptive procedure is employed to ensure two significant digits in the velocity approximation, yielding the map of adapted degrees in Figure\ \ref{fg:dropKadaptDeg}.  The results show that the highest polynomial degree $k=4$ is achieved in the vicinity of the central pore and, specifically, for the smallest droplets, whereas low-order polynomial functions are sufficient to accurately capture the solution away from the interfaces between fluid $1$ and fluid $2$.
Hence, the combination of high-order functional approximation and  exact description of the geometry via NURBS is also applicable to multi-fluid systems, with the unfitted framework allowing to bypass the difficulties of generating high-order curved meshes.

%
\section{Concluding remarks}
\label{sc:Conclusions}
In this article, an unfitted high-order HDG method capable of treating exact geometries was presented.
The method was successfully applied to incompressible one-fluid and two-fluid systems, with boundaries and interfaces described by means of NURBS, showcasing optimal high-order convergence and robustness even in the presence of badly cut cells or faces.

The non-conforming nature of this approach allows to achieve high-fidelity solutions, circumventing the generation of fitted, high-order, curved meshes. To this end, the method embeds the NURBS in a fixed Cartesian grid and an algorithm to intersect the NURBS with the mesh skeleton is executed. Note that this is trivial in 2D, whereas extension to 3D is more cumbersome, but doable using standard CAD tools\ \cite{Sevilla-MSZRT-15,Rodenas-MRNT-17}. 

Once the sets of cut elements and cut faces are identified,  the unfitted HDG formulation only requires the appropriate integration of the equations in the physical domains for both element-based and face-based unknowns.
On the one hand, uncut elements are treated as classic elements in standard HDG formulations. On the other hand, \emph{immersed} elements cut by an unfitted boundary and \emph{interface} elements featuring two fluids require special treatment: integration along immersed NURBS and integration in elements cut by a set of NURBS curves are accurately performed following the NEFEM rationale.
In addition, to appropriately handle interface conditions, the method duplicates the unknowns in the mesh elements and on the mesh faces cut by the interface, thus enriching the approximation space to achieve a more precise discretization.
Note that whilst all the discussed two-fluid numerical tests are static, the presented methodology can be extended to account for moving interfaces and deformable boundaries. This can be achieved following the idea to deform NURBS curves via a suitable movement of the control points, as discussed in the context of NURBS-based treatment of free-surface flows and fluid-structure interaction systems\ \cite{Elgeti-SKEH-16,Behr-HHSBE-18}.

Hence, by combining non-conforming meshes, exact NURBS geometry and high-order degree-adaptive approximations, the presented unfitted HDG method provides accurate results on coarse meshes, which are independent of the geometric features of the domain.
Moreover,  the use of coarse meshes minimizes the number of intersections between the mesh skeleton and the boundaries/interfaces, consequently reducing the computational overhead due to cut elements.
Finally, the resulting unfitted HDG-NEFEM scheme does not introduce any new unknown on non-matching boundaries or interfaces,  always solving for the global unknowns defined on the grid skeleton and preserving the favorable matrix structure of standard HDG methods.

\section*{Acknowledgements}
This work was supported by the Spanish Ministry of Science, Innovation and Universities and the Spanish State Research Agency MICIU/AEI/10.13039/501100011033 (PID2020-113463RB-C33 to MG, PID2020-113463RB-C32 to AH, CEX2018-000797-S to MG and AH) and the Generalitat de Catalunya (2021-SGR-01049). MG is Fellow of the Serra H\'unter Programme of the Generalitat de Catalunya.

\bibliographystyle{unsrt}
\bibliography{Ref-Immersed}

\appendix

\section{Identification of the intersections of NURBS curves with the computational mesh}
\label{ap:intersect}

In this appendix, the strategy to identify the intersection of a NURBS curve with the computational mesh $\meshDomain$ is briefly presented.

First,  for each NURBS curve $\mathcal{\pmb{C}}^j_{\Su}$, $j\inSet{1}{\nSu}$ composing the interface $\Su$, a \emph{sufficiently large} number (e.g., $n_j = 100$) of sampling points $\{\bx^j_k\}_{k\inSet{1}{n_j}}$ is introduced to represent the NURBS. The points $\bx_k^j$ are set equispaced in the parametric NURBS domain $[0,1]$, that is, $\bx^j_k = \mathcal{\pmb{C}}^j_{\Su}(\lambda^j_k)$, with $\lambda^j_k = k/n_j$, and $k\inSet{0}{n_j}$.

For all $k\inSet{0}{n_j}$ and all $j\inSet{1}{\nSu}$,  the mesh element $\Omega_e$ containing $\bx^j_k$ is identified, thus providing the list of \emph{cut elements}.  Note that this strategy can be easily optimized, via parallelization or by means of efficient sequential algorithms searching for cut elements only in the neighborhood of already defined cut elements.

Once cut elements are identified, cut faces need to be determined. To this end, for each cut element $\Omega_e$ and all $k\inSet{1}{n_j-1}$,  the procedure seeks the points $\bx^j_m, \ m\inSet{k-1}{k+1}$,  neighbors of $\bx^j_k$, and lying outside of $\Omega_e$.
Given a point $\bx^j_m$ outside of $\Omega_e$, owing to the continuity of $\mathcal{\pmb{C}}^j$, the intersection point $\overline{\bx}^j_{e,n} = \mathcal{\pmb{C}}^j_{\Su}(\lambda^j_{n})$ between $\Su$ and $\partial\Omega_e$ is determined via a dichotomy algorithm, with $\lambda^j_{n} \in [\lambda^j_{k-1},\lambda^j_{k+1}]$.

Hence, the portion of the interface lying within the cut element is denoted by $\Su_e = \Su \cap \Omega_e$ and its intersections with $\partial \Omega_e$ are given by the points $\overline{\bx}^j_{e,n}$, with $n=1,2$.  
Under the assumption that the interface does not intersect an element through a vertex,  for any intersection points $\overline{\bx}^j_{e,n}$, with $j\inSet{1}{\nSu}$ and $n\inSetTwo{1}{2}$,  between $\Su$ and $\partial\Omega_e$,  there exists a corresponding intersection point  $\overline{\bx}^j_{k,n} = \overline{\bx}^j_{e,n}$, lying in a neighboring element $\Omega_k$ which shares a face with $\Omega_e$.
If this condition is not satisfied, a local refinement is performed to look for missing cut elements. 
  
\begin{remark}[NURBS sampling]
The previously mentioned approximation of $\Su$ using $n_j$ sampling points is only employed to identify the intersections of the NURBS with the mesh skeleton. 
Note that alternative strategies, beyond the scope of the present work, could be employed to improve the efficiency of the algorithm to identify the cut elements, e.g., by exploiting \emph{a priori} information on the local size and structure of the computational mesh. 
It is worth recalling that these operations are performed once, during a preprocessing stage, and do not introduce a significant computational burden in the algorithm.  
On the contrary, during computation, the information of the \emph{true} boundary/interface, described by the exact NURBS curve, is exploited to perform numerical integration using the NEFEM rationale\ \cite{Sevilla-SFH-08-IJNME}.
%
\end{remark}
  
\begin{remark}[Multiple unfitted regions]
The above considerations for the identification of the intersections rely on the assumption that a cut element is divided, at most, in two regions. 
Consider now an element $\Omega_e$ with $n_e^j$ intersection points between $\Su$ and $\partial\Omega_e$, that is, $\overline{\bx}^j_{e,n}$, with $n \inSet{1}{n_e^j}$. 
Hence, the element $\Omega_e$ is split into $n_e^j/2 +1$ disjoint regions and $\Su_e$ is composed of $n_e^j-1$ disjoint curves $\Su_{e,n}$ such that
\begin{equation*}
    \Su_e = \bigSum{\ell}{1}{n_e^j-1}\Su_{e,\ell},
\end{equation*}
where each curve $\Su_{e,\ell}$ is continuous and can be composed of multiple NURBS curves  $\pmb{C}^j_{\Su}$.
This leads to
\begin{equation*}
   \Su_{e,\ell} = \bigcup_{j\in\mathcal{J}_{e,\ell}} \pmb{C}^j_{\Su}(\lambda), \qquad \forall \lambda\in[0,1] ,
\end{equation*}
%
with $\mathcal{J}_{e,\ell}$ being the set of coefficients $j$ of the NURBS curves $\pmb{C}^j_{\Su}$ describing the portion $\Su_{e,\ell}$ of the interface.
\end{remark}

\section{Triangulation procedure for quadrature involving an element cut by NURBS curves}
\label{ap:quadTri}
\begin{figure}[!h]
  \centering
  \subfloat[Cut element.]{\includegraphics[width=0.38\textwidth]{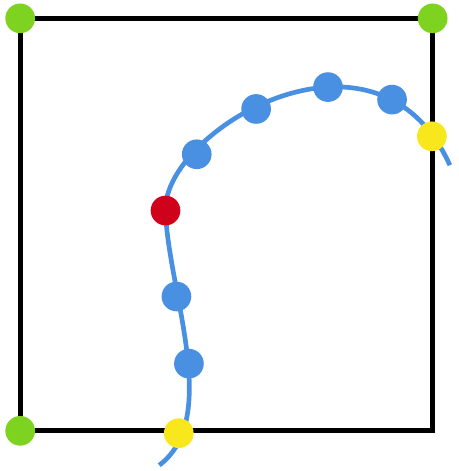}\label{fig:triangulationA}}
  \quad
  \subfloat[Visibility of the orange vertex.]{\includegraphics[width=0.38\textwidth]{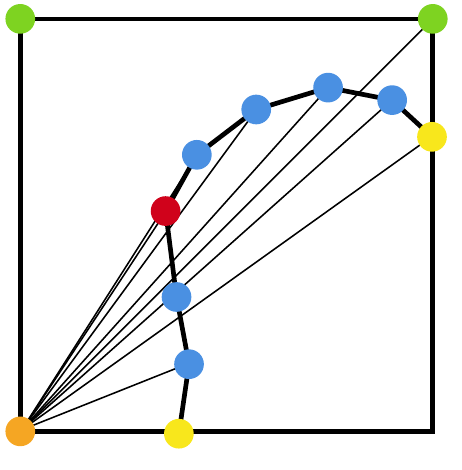}\label{fig:triangulationB}}
  
  \subfloat[Visibility of the orange point.]{\includegraphics[width=0.38\textwidth]{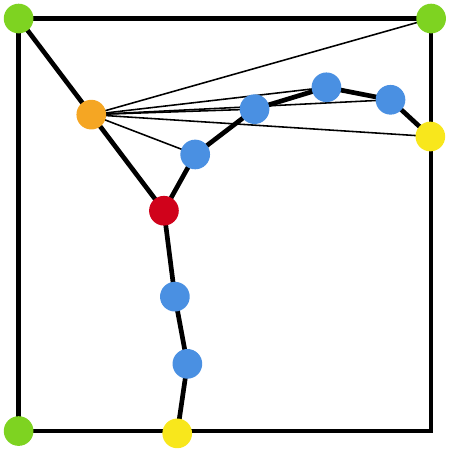}\label{fig:triangulationC}}
  \quad
  \subfloat[Triangulation for quadrature.]{\includegraphics[width=0.38\textwidth]{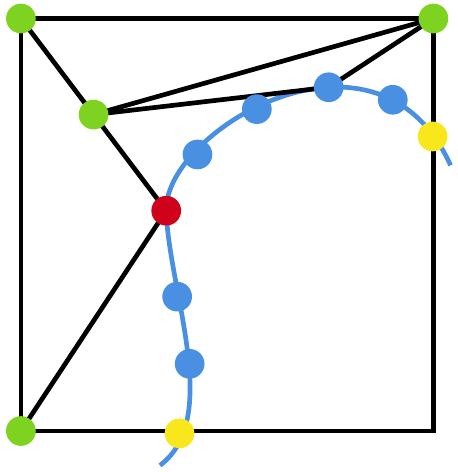}\label{fig:triangulationD}}
  \caption{Triangulation procedure for a region $\Omega_e^i$, within an element $\Omega_e$ cut by a NURBS curve $\Su_e$.}
  \label{fig:triangulation}
\end{figure}

In this appendix, the quadrature procedure for an element $\Omega_e$, cut by a NURBS curve $\Su_e$, is presented. Given a cut element,  quadrature is performed following the NEFEM rationale\ \cite{Sevilla-SRH-16}. More precisely,  for each region $\Omega_e^i$, a triangulation is constructed, as sketched in Figure\ \ref{fig:triangulation}.

First, the region $\Omega_e^i$ is identified (Figure\ \ref{fig:triangulationA}) according to the procedure presented in Section\ \ref{sec:identifyFluid}. This region is described by the vertices of $\Omega_e$ lying in $\Omega_e^i$ (green vertices),  the extrema of $\Su_e^j$, with $j\inSet{1}{\nSu}$, intersecting $\partial\Omega_e$ (yellow points) and a set of points sampling the NURBS curves $\Su_e^j$.  More precisely, in Figure\ \ref{fig:triangulationA}, the interface $\Su_e$ is composed of two NURBS curves, split by the red point, and sampled at the locations of the six blue points. As a general rule, the number of sampling points along the NURBS is selected proportional to the inverse of the curvature of $\Su_e$, to accurately describe the curve.

Then, an initial point $\bx_0 \in \Omega_e^i$ is selected (e.g., a vertex or the barycenter of a face belonging to $\partial\Omega_e^i$) to execute Lee's visibility algorithm\ \cite{Joe:87}.  In Figure\ \ref{fig:triangulationB}, the starting point is the bottom-left orange vertex. The goal of this procedure is to identify the subregions of $\Omega_e^i$ visible from $\bx_0$,  that is,  to determine all the sampling points that can be connected to $\bx_0$, without intersecting the NURBS curves $\Su_e$ more than once. The visibility region is thus defined as the polygon delimited by $\bx_0$ and all the visible points. In Figure\ \ref{fig:triangulationB},  the visibility region is identified by the orange vertex, the yellow point at the bottom and the red point.

The procedure is then repeated for the portion of the NURBS not yet explored. Following\ \cite{Sevilla-SRH-16}, the second starting point is selected as the mid-point between the last visible point (red point) and a vertex (top-left green point).  Figure\ \ref{fig:triangulationC} assesses the visibility region starting from the orange point.

The procedure is iteratively repeated until $\Omega_e^i$ is completely partitioned into visible subregions, as reported in Figure\ \ref{fig:triangulationD}.  The resulting triangles feature, at most,  one curved edge, exactly describing the geometry of $\Su_e$ by means of one NURBS curve. Generalizations of these results are discussed in\ \cite{Sevilla-SRH-16}.
Finally, it is worth recalling that this triangulation is only constructed to perform numerical quadrature on $\Omega_e^i$, whereas the functional approximation is maintained in the original nodes of the underlying Cartesian grid.

\end{document}